\documentclass[final]{siamltex}

\usepackage{amsmath} 
\usepackage{amsfonts}

\usepackage{graphicx}
\usepackage{subcaption}
\usepackage{url}
\usepackage{color}
\usepackage{bm}
\usepackage{multirow}

\newtheorem{remark}{Remark}

\title{Parameter-Robust Preconditioning for Oseen Iteration Applied to Stationary and Instationary Navier--Stokes Control}

\author{Santolo Leveque\thanks{School of Mathematics, The University of Edinburgh, James Clerk Maxwell Building, The King's Buildings, Peter Guthrie Tait Road, Edinburgh, EH9 3FD, United Kingdom ({\tt S.Leveque@sms.ed.ac.uk})} \and John W. Pearson\thanks{School of Mathematics, The University of Edinburgh, James Clerk Maxwell Building, The King's Buildings, Peter Guthrie Tait Road, Edinburgh, EH9 3FD, United Kingdom ({\tt j.pearson@ed.ac.uk})}}

\begin{document}
\maketitle

\begin{abstract}
We derive novel, fast, and parameter-robust preconditioned iterative methods for steady and time-dependent Navier--Stokes control problems. Our approach may be applied to time-dependent problems which are discretized using backward Euler or Crank--Nicolson, and is also a valuable candidate for Stokes control problems discretized using Crank--Nicolson. The key ingredients of the solver are a saddle-point type approximation for the linear systems, an inner iteration for the $(1,1)$-block accelerated by a preconditioner for convection--diffusion control, and an approximation to the Schur complement based on a potent commutator argument applied to an appropriate block matrix. A range of numerical experiments validate the effectiveness of our new approach.
\end{abstract}



\pagestyle{myheadings}
\thispagestyle{plain}
\markboth{S. LEVEQUE AND J. W. PEARSON}{ROBUST PRECONDITIONING FOR NAVIER--STOKES CONTROL}

\section{Introduction}
\label{1}
Optimal control problems with PDE constraints have received increasing research interest of late, due to their applicability to scientific and industrial problems, and also due to the difficulties arising in their numerical solution (see \cite{Hinze, Troltzsch} for excellent overviews of the field). An example of a highly challenging problem attracting significant attention is the (distributed) control of incompressible viscous fluid flow problems. Here, the constraints may be the (non-linear) incompressible Navier--Stokes equations or, in the limiting case of viscous flow, the (linear) incompressible Stokes equations. The control of the Navier--Stokes equations is of particular interest: due to the non-linearity involved, to find a solution linearizations of the constrained problem need to be repeatedly solved until a sufficient reduction on the non-linear residual is achieved \cite{Hintermuller_Hinze, Hinze_Koster_Turek, Posta_Roubicek}. This has motivated researchers to devise solvers for this type of problem which exhibit robustness with the respect to all the parameters involved; see \cite{Hinze_Koster_Turek} for a robust multigrid method applied to Newton iteration for instationary Navier--Stokes control, for instance. Despite the recent development of parameter-robust preconditioners for the control of the (stationary and instationary time-periodic) Stokes equations \cite{Axelsson_Farouq_Neytcheva, Krendl_Simoncini_Zulehner,Rees_Wathen, Zulehner}, to our knowledge no such preconditioner has proved to be completely robust when applied to the Navier--Stokes control problem considered here. We also point out \cite{Heidel_Wathen} for a preconditioned iterative solver for Stokes and Navier--Stokes boundary control problems, and \cite{Qiu_vanGijzen_vanWingerden_Verhaegen_Vuik} for an efficient and robust preconditioning technique for in-domain Navier--Stokes control.

A popular preconditioner for the \emph{Oseen linearization} of the forward stationary Navier--Stokes equations combines saddle-point theory with a commutator argument for approximating the Schur complement \cite{Kay_Loghin_Wathen}. This type of preconditioner shows only a mild dependence on the viscosity parameter, and is robust with respect to the discretization parameter. In \cite{Pearson_Navier_Stokes} the combination of saddle-point theory with a commutator argument has also been adapted to the control of the stationary Navier--Stokes equations; we note that a commutator argument of this type was first introduced in \cite{Pearson_Stokes} for the control of the stationary Stokes equations. In this work, we utilize a commutator argument for a block matrix in conjunction with saddle-point theory in order to derive robust preconditioners for the optimal control of the incompressible Navier--Stokes equations, in both the stationary and instationary settings. For instationary problems our approach leads to potent preconditioners when either the backward Euler or Crank--Nicolson scheme is used in the time variable, and also leads to a preconditioner for the instationary Stokes control problem using Crank--Nicolson.

This article is structured as follows. In Section \ref{2}, we define the problems that we examine, that is the stationary and instationary Navier--Stokes control problems along with instationary Stokes control; we then present the linearization adopted in this work, and outline the linear systems arising upon discretization of the forward problem. In Section \ref{3}, we introduce the saddle-point theory used to devise optimal preconditioners, giving as an example a preconditioner for the forward stationary Navier--Stokes equation in combination with the commutator argument presented in \cite{Kay_Loghin_Wathen}; the latter will then be generalized when multiple differential operators are involved in the system of equations. In Section \ref{4}, we derive the first-order optimality conditions of the control problems and their discretization. In Section \ref{5}, we present our suggested preconditioners, and in particular the commutator argument applied to the Schur complements arising from the control problems. Then, in Section \ref{6} we provide numerical results that show the robustness and efficiency of our approach.

\section{Problem Formulation}
\label{2}
In this work we derive fast and robust preconditioned iterative methods for the distributed control of incompressible viscous fluid flow; in this case, the physics is described by the (stationary or instationary) incompressible Navier--Stokes equations. The corresponding distributed control problem is defined as a minimization of a least-squares cost functional subject to the PDEs.

Specifically, given a spatial domain $\Omega \subset \mathbb{R}^{d}$, $d\in\{2,3\}$, the stationary Navier--Stokes control problem we consider is
\begin{equation}\label{Stationary_Navier_Stokes_control_functional}
\min_{\vec{v},\vec{u}}~~J_{\rm{S}}(\vec{v},\vec{u}) = \dfrac{1}{2} \int_{\Omega} |\vec{v}(\mathbf{x})- \vec{v}_{d}(\mathbf{x})|^2 \: {\rm d}x + \dfrac{\beta}{2} \int_{\Omega} |\vec{u}(\mathbf{x})|^2 \: {\rm d}x
\end{equation}
subject to
\begin{equation}\label{Stationary_Navier_Stokes_control_constraints}
\left\{
\begin{array}{rl}
\vspace{0.25ex}
- \nu \nabla^2 \vec{v} + \vec{v} \cdot \nabla \vec{v} + \nabla p = \vec{u} + \vec{f}(\mathbf{x}) & \quad \mathrm{in} \; \Omega, \\
\vspace{0.25ex}
- \nabla \cdot \vec{v} (\mathbf{x}) = 0 & \quad \mathrm{in} \; \Omega, \\
\vec{v}(\mathbf{x}) = \vec{g}(\mathbf{x}) & \quad \mathrm{on} \; \partial \Omega,
\end{array}
\right.
\end{equation}
where the \emph{state variables} $\vec{v}$ and $p$ denote velocity and pressure respectively, $\vec{v}_d$ is the \emph{desired state (velocity)}, and $\vec{u}$ is the \emph{control variable}. Further, $\beta>0$ is a \emph{regularization parameter}, and $\nu>0$ is the \emph{viscosity} of the fluid. The functions $\vec{f}$ and $\vec{g}$ are known.

Similarly, the control of the instationary Navier--Stokes equations is defined as
\begin{equation}\label{Instationary_Navier_Stokes_control_functional}
\min_{\vec{v},\vec{u}}~~J_{\rm{I}}(\vec{v},\vec{u}) = \dfrac{1}{2} \int_{0}^{t_{f}} \!\! \int_{\Omega} |\vec{v}(\mathbf{x},t)- \vec{v}_{d}(\mathbf{x},t)|^2 \: {\rm d}x \: {\rm d}t + \dfrac{\beta}{2} \int_{0}^{t_{f}} \!\! \int_{\Omega} |\vec{u}(\mathbf{x},t)|^2 \: {\rm d}x \: {\rm d}t,
\end{equation}
given also a final time $t_{f}>0$, subject to
\begin{equation}\label{Instationary_Navier_Stokes_control_constraints}
\left\{
\begin{array}{rl}
\vspace{0.25ex}
\frac{\partial \vec{v}}{\partial t} - \nu \nabla^2 \vec{v} + \vec{v} \cdot \nabla \vec{v} + \nabla p = \vec{u} + \vec{f}(\mathbf{x},t) & \quad \mathrm{in} \; \Omega \times (0,t_{f}), \\
\vspace{0.25ex}
- \nabla \cdot \vec{v} (\mathbf{x},t) = 0 & \quad \mathrm{in} \; \Omega \times (0,t_{f}), \\
\vspace{0.25ex}
\vec{v}(\mathbf{x},t) = \vec{g}(\mathbf{x},t) & \quad \mathrm{on} \; \partial \Omega \times (0,t_{f}),\\
\vec{v}(\mathbf{x},0) = \vec{v}_0(\mathbf{x}) & \quad \mathrm{in} \; \Omega,

\end{array}
\right.
\end{equation}
using the same notation as above. As for the stationary case, the functions $\vec{f}$ and $ \vec{g}$ are known; the initial condition $\vec{v}_0$ is also given. In the following, we assume that $\vec{v}_0$ is solenoidal, i.e. $\nabla \cdot \vec{v}_0 = 0$, adapting our strategy to the general case when possible.

The constraints \eqref{Stationary_Navier_Stokes_control_constraints} and \eqref{Instationary_Navier_Stokes_control_constraints} are a system of non-linear (stationary or instationary) PDEs. In order to obtain a solution of the corresponding control problem, we make use of the Oseen linearization of the non-linear term $\vec{v} \cdot \nabla \vec{v}$, as in \cite{Posta_Roubicek}.

If the non-linear term $\vec{v} \cdot \nabla \vec{v}$ is dropped in \eqref{Stationary_Navier_Stokes_control_constraints} or \eqref{Instationary_Navier_Stokes_control_constraints}, with $\nu = 1$, we obtain the corresponding distributed Stokes control problem. Many parameter-robust preconditioners for the stationary Stokes control problem have been derived in the literature \cite{Rees_Wathen, Zulehner}; however, less progress has been made towards the parameter-robust solution of instationary Stokes control problems, except in the time-periodic setting \cite{Axelsson_Farouq_Neytcheva, Krendl_Simoncini_Zulehner}. Below, we derive a preconditioner that will also result in a robust solver for the general formulation of the instationary distributed Stokes control problem, defined as the minimization of the functional \eqref{Instationary_Navier_Stokes_control_functional} subject to the instationary Stokes equations.

\subsection{Non-linear iteration and discretization matrices}
\label{2_1}
To introduce the linearization adopted for the control case as well as the discretization matrices, we consider the stationary Navier--Stokes equations:
\begin{equation}\label{forward_Navier_Stokes_control_constraints}
\left\{
\begin{array}{rl}
- \nu \nabla^2 \vec{v} + \vec{v} \cdot \nabla \vec{v} + \nabla p = \vec{u} + \vec{f}(\mathbf{x}) & \quad \mathrm{in} \; \Omega, \\
- \nabla \cdot \vec{v} (\mathbf{x}) = 0 & \quad \mathrm{in} \; \Omega,
\end{array}
\right.
\end{equation}
with $\vec{v} = \vec{g}$ on $\partial \Omega$. First, we introduce the weak formulation of \eqref{forward_Navier_Stokes_control_constraints} as follows. Let $V : = \{\vec{v} \in \mathcal{H}^1(\Omega)^d \: | \: \vec{v} = \vec{g} \; \mathrm{on} \; \partial \Omega \}$, $V_0 : = \{\vec{v} \in \mathcal{H}^1(\Omega)^d \: | \: \vec{v} = \vec{0} \; \mathrm{on} \; \partial \Omega \}$, and $Q:= L^2(\Omega)$, with $\mathcal{H}^1(\Omega)^d$ the Sobolev space of square-integrable functions in $\mathbb{R}^d$ with square-integrable weak derivatives; then, the weak formulation reads as:

Find $\vec{v} \in V$ and $p \in Q$ such that
\begin{equation}\label{weak_formulation_forward_Navier_Stokes_control_constraints}
\left\{
\begin{array}{rl}
 \nu (\nabla \vec{v}, \nabla \vec{w}) + (\vec{v} \cdot \nabla \vec{v}, \vec{w}) - ( p,\nabla \cdot \vec{w} ) = (\vec{u}, \vec{w}) + (\vec{f},\vec{w}) & \quad \mathrm{for \; all} \; \vec{w}\in V_0, \\
- (q, \nabla \cdot \vec{v} \: ) = 0 & \quad \mathrm{for \; all} \;q \in Q,
\end{array}
\right.
\end{equation}
where $(\cdot,\cdot)$ is the $L^2$-inner product on $\Omega$. The main issue in \eqref{weak_formulation_forward_Navier_Stokes_control_constraints} is how to deal with the non-linear term $(\vec{v} \cdot \nabla \vec{v}, \vec{w})$. A common strategy employs the \emph{Picard iteration}, which is described as follows.  Given the approximations $\vec{v}^{\: (k)} \in V$ and $p^{(k)} \in Q$ to $\vec{v}$ and $p$ respectively, we consider the non-linear residuals:
\begin{equation}\label{forward_non_linear_residuals}
\left\{
\begin{array}{r}
\vec{R}^{\: (k)} = (\vec{u}, \vec{w}) + (\vec{f},\vec{w}) - \nu ( \nabla \vec{v}^{\: (k)}, \nabla \vec{w}) - (\vec{v}^{\: (k)} \cdot \nabla \vec{v}^{\: (k)},\vec{w}) + ( p^{(k)}, \nabla\cdot \vec{w}), \\
r^{(k)} = (q,\nabla \cdot \vec{v}^{\: (k)}),
\end{array}
\right.
\end{equation}
for any $\vec{w} \in V_0$ and $q\in Q$. Then, the Picard iteration is defined as \cite[p.\,345--346]{Elman_Silvester_Wathen}:
\begin{equation}\label{forward_Picard_step}
\vec{v}^{\:(k+1)}=\vec{v}^{\:(k)}+ \vec{\delta v}^{(k)}, \qquad 
p^{(k+1)}=p^{(k)}+ \delta p^{(k)},
\end{equation}
where $\vec{\delta v}^{(k)}$ and $\delta p^{(k)}$ are the solutions of
\begin{equation}\label{forward_Picard_problems}
\left\{
\begin{array}{r}
 \nu (\nabla \vec{\delta v}^{\: (k)}, \nabla \vec{w}) + (\vec{v}^{\: (k)} \cdot \nabla \vec{\delta v}^{\: (k)}, \vec{w}) - ( \delta p^{(k)}, \nabla \cdot\vec{w})  = \vec{R}^{\: (k)}, \\
- (q,\nabla \cdot \vec{\delta v}^{\: (k)}) = r^{(k)},

\end{array}
\right.
\end{equation}
for any $\vec{w} \in V_0$ and $q\in Q$. Equations \eqref{forward_Picard_problems} are the Oseen equations for the forward stationary Navier--Stokes equations. These are posed on the continous level, so in order to find a solution to \eqref{forward_Navier_Stokes_control_constraints} we now discretize them. Before doing so, we note that \eqref{forward_Picard_problems} represents an incompressible convection--diffusion equation, with wind vector defined by $\vec{v}^{\: (k)}$, and it is clear that, for $\nu \ll 1$, the problem is convection-dominated. This requires us to make use of a stabilization procedure.

\color{black}Letting $\{\vec{\phi}_i\}_{i=1}^{n_v}$ and $\{\psi_i\}_{i=1}^{n_p}$ be an inf--sup stable finite element basis functions for the Navier--Stokes equations, we seek approximations $\vec{v}(\mathbf{x}) \approx \sum_{i=1}^{n_v} \mathbf{v}^{(k)}_i \vec{\phi}_i$, $\vec{u}(\mathbf{x}) \approx \sum_{i=1}^{n_v} \mathbf{u}_i \vec{\phi}_i$, $p(\mathbf{x}) \approx \sum_{i=1}^{n_p} p^{(k)}_i \psi_i$. Denoting the vectors $\mathbf{v}^{(k)} = \{\mathbf{v}^{(k)}_i\}_{i=1}^{n_v}$, $\mathbf{u} = \{\mathbf{u}_i\}_{i=1}^{n_v}$, $\mathbf{p}^{(k)} = \{ p^{(k)}_i \}_{i=1}^{n_p}$, a discretized version of \eqref{forward_non_linear_residuals} is:
\begin{displaymath}
\left\{
\begin{array}{rl}
\mathbf{R}^{ (k)}  =& \hspace{-0.5em} \mathbf{M} \mathbf{u} + \mathbf{f} - \mathbf{L}^{(k)} \: \mathbf{v}^{\: (k)} - B^\top \mathbf{p}^{(k)}  , \\
\mathbf{r}^{(k)}  = & \hspace{-0.5em} - B \: \mathbf{v}^{\: (k)}  ,
\end{array}
\right.
\end{displaymath}
where we set $\mathbf{L}^{(k)} =\nu \mathbf{K} + \mathbf{N}^{(k)} + \mathbf{W}^{(k)} $, with
\begin{align*}
\mathbf{K}=\{k_{ij}\}_{i,j=1}^{n_v}, & \;\; k_{ij} \! = \! \int_\Omega \nabla \vec{\phi}_i : \nabla \vec{\phi}_j , &
\mathbf{N}^{(k)}=\{n^{(k)}_{ij}\}_{i,j=1}^{n_v}, & \;\;  n^{(k)}_{ij} \! = \! \int_\Omega (  \vec{v}^{\: (k)} \cdot \nabla \vec{\phi}_j ) \cdot \vec{\phi}_i , \\
\mathbf{M}=\{m_{ij}\}_{i,j=1}^{n_v}, & \;\; m_{ij} \! = \! \int_\Omega \vec{\phi}_i \cdot \vec{\phi}_j , &
B = \{b_{ij} \}_{i = 1, \ldots, n_p}^{j = 1, \ldots, n_v}, & \;\; b_{ij} \! = \! - \int_\Omega \psi_i  \nabla \cdot \vec{\phi}_j , \\
\mathbf{f}=\{f_{i}\}_{i=1}^{n_v}, & \;\; f_i \! = \! \int_\Omega \vec{f} \cdot  \vec{\phi}_i ,
\end{align*}
and the matrix $\mathbf{W}^{(k)}$ denotes a possible stabilization matrix for the convection operator\color{black}. Then, the Picard iterate \eqref{forward_Picard_step} may be written in discrete form as
\begin{displaymath}
\mathbf{v}^{\:(k+1)}=\mathbf{v}^{\:(k)}+ \bm{\delta v}^{(k)}, \qquad 
\mathbf{p}^{(k+1)}=\mathbf{p}^{(k)}+ \bm{\delta p}^{(k)},
\end{displaymath}
with $\bm{\delta v}^{ (k)}$ and $\bm{\delta p}^{(k)}$ solutions of
\begin{equation}\label{forward_discrete_Picard_problems}
\left\{
\begin{array}{rl}
\mathbf{L}^{(k)} \: \bm{\delta v}^{ (k)} + B^\top \: \bm{\delta p}^{(k)} = & \hspace{-0.5em} \mathbf{R}^{(k)}, \\
B \: \bm{\delta v}^{ (k)} = & \hspace{-0.5em} \mathbf{r}^{(k)}.
\end{array}
\right.
\end{equation}
\color{black}The matrix $\mathbf{K}$ is generally referred to as a \emph{(vector-)stiffness matrix}, and the matrix $\mathbf{M}$ is referred to as a \emph{(vector-)mass matrix}; both the matrices are symmetric positive definite (s.p.d.). The matrix $\mathbf{N}^{(k)}$ is referred to as a \emph{(vector-)convection matrix}, and is skew-symmetric (i.e. $\mathbf{N}^{(k)}+(\mathbf{N}^{(k)})^\top=0$) if the incompressibility constraints $\nabla \cdot \vec{v}^{\: (k)}=0 $ are solved exactly\color{black}; finally, the matrix $B$ is the \emph{(negative) divergence matrix}.

Regarding the stabilization procedure applied, we note that the matrix $\mathbf{W}^{(k)}$ represents a differential operator that is not physical, and is introduced only to enhance coercivity (that is, increase the positivity of the real part of the eigenvalues) of the discretization, thereby allowing it to be stable. 
For the reasons discussed in \cite{Leveque_Pearson, Pearson_Wathen}, in the following we employ the Local Projection Stabilization (LPS) approach described in \cite{Becker_Braack, Becker_Vexler, Braack_Burman}. We point out \cite{Matthies_Tobiska}, where the authors derive the order of convergence of one- and two-level LPS applied to the Oseen problem. For other possible stabilizations applied to the Oseen problem, see \cite{Brooks_Hughes, Franca_Frey, Gelhard_Lube_Olshanskii_Starcke, Johnson_Saranen, Tobiska_Lube}.

\color{black}In the LPS formulation, the stabilization matrix $\mathbf{W}^{(k)}$ is defined as
\begin{equation}\label{LPS_stab_matrix}
\mathbf{W}^{(k)}=\{w^{(k)}_{ij}\}_{i,j=1}^{n_v}, \qquad w^{(k)}_{ij} = \delta^{(k)} \int_\Omega \kappa_h(\vec{v}^{\: (k)} \cdot \nabla \vec{\phi}_i) \cdot \kappa_h(\vec{v}^{\: (k)} \cdot \nabla \vec{\phi}_j) \; .
\end{equation}
Here, $\delta^{(k)}>0$ denotes a stabilization parameter, and $\kappa_h=\text{Id} +\pi_h$ is the fluctuation operator, with $\text{Id}$ the identity operator and $\pi_h$ an $L^2$-orthogonal (discontinuous) projection operator defined on patches of $\Omega$, where by a patch we mean a union of elements of our finite element discretization. In our implementation, the domain is divided into patches consisting of two elements in each dimension. Further, we define the stabilization parameter $\delta^{(k)}$ and the local projection $\pi_h$ as in \cite[Sec.\,2.1]{Leveque_Pearson}\color{black}.

\section{Saddle-Point Systems}
\label{3}
In this section we introduce saddle-point theory and the commutator argument derived in \cite{Kay_Loghin_Wathen} for the forward stationary Navier--Stokes equations. These will be the main ingredients for devising our preconditioners.

\color{black}Given an invertible system of the form
\begin{equation}\label{MatrixA}
\mathcal{A} =
\left[
\begin{array}{cc}
\Phi & \Psi_1\\
\Psi_2 & -\Theta
\end{array}
\right]
\end{equation}
with $\Phi$ invertible, a good candidate for a preconditioner is the block triangular matrix:
\begin{equation}\label{optimal_prec}
\mathcal{P}=
\left[
\begin{array}{cc}
\Phi & 0 \\
\Psi_2 & -S
\end{array}
\right],
\end{equation}
where $S$ denote the (negative) \emph{Schur complement}  $S=\Theta+\Psi_2 \Phi^{-1}\Psi_1$. Indeed, if $S$ is also invertible and denoting the set of eigenvalues of a matrix by $\lambda(\cdot)$, we have $\lambda (\mathcal{P}^{-1}\mathcal{A}) = \left\{1\right\}$, see \cite{Ipsen01,Murphy:1999:NPI:359189.359190}. Since the preconditioner is not symmetric, we need to employ a Krylov subspace method for non-symmetric matrices, such as GMRES \cite{Saad_Schultz}. Clearly, we do not want to apply the inverse of $\mathcal{P}$ as defined in \eqref{optimal_prec}, as the computational cost would be comparable to that of applying the inverse of $\mathcal{A}$. In particular, applying $S^{-1}$ would be problematic, as even when $\Phi$ and $\Theta$ are sparse $S$ is generally dense. For this reason, we consider a suitable approximation:
\begin{equation}\label{appr_preconditioner}
\widehat{\mathcal{P}} = 
\left[
\begin{array}{cc}
\widehat{\Phi} & 0 \\
\Psi_2 & -\widehat{S}
\end{array}
\right]
\end{equation}
of $\mathcal{P}$, or, more precisely, a cheap application of the effect of $\widehat{\mathcal{P}}^{-1}$ on a generic vector\color{black}. For instance, an efficient preconditioner for the matrix arising from \eqref{forward_discrete_Picard_problems} is given by \eqref{appr_preconditioner}, with $\widehat{\Phi}$ being the approximation of $\mathbf{L}^{(k)}$ using a multigrid routine, and $\widehat{S}$ the so called \emph{pressure convection--diffusion preconditioner} \cite[p.\,365--370]{Elman_Silvester_Wathen} (first derived in \cite{Kay_Loghin_Wathen}) for $S$. The latter is derived by mean of a commutator argument as follows. Consider the convection--diffusion operator $\mathcal{D}=-\nu \nabla^2 + \vec{v}^{\: (k)} \cdot \nabla$ defined on the velocity space as in \eqref{forward_Picard_problems}, and suppose the analogous operator $\mathcal{D}_p=(-\nu \nabla^2 + \vec{v}^{\: (k)} \cdot \nabla)_p$ on the pressure space is well defined. Suppose also that the commutator
\begin{equation}\label{forward_commutator}
\mathcal{E} = \mathcal{D} \nabla - \nabla \mathcal{D}_p
\end{equation}
is small in some sense. Then, discretizing \eqref{forward_commutator} with stable finite elements leads to
\begin{displaymath}
(\mathbf{M}^{-1} \mathbf{L}^{(k)}) \mathbf{M}^{-1} B^\top - \mathbf{M}^{-1} B^\top (M_p^{-1} L^{(k)}_p) \approx 0,
\end{displaymath}
where $L^{(k)}_p=\nu K_p + N_p^{(k)} + W_p^{(k)}$ is the discretization of $\mathcal{D}_p$ in the finite element basis for the pressure, with $M_p= \left[ ( \psi_i , \psi_j ) \right]$, $K_p= \left[ ( \nabla \psi_i , \nabla \psi_j ) \right]$, $N_p^{(k)}= \left[ (\vec{v}^{\: (k)} \:\! \cdot \:\! \nabla \psi_j, \psi_i ) \right]$, and $W_p^{(k)}= \left[ \delta^{(k)} (\kappa_h (\vec{v}^{\: (k)} \cdot \nabla \psi_i) , \kappa_h(\vec{v}^{\: (k)} \cdot \nabla \psi_j ) )\right]$ the (scalar) mass, stiffness, convection, and stabilization matrices, respectively, in the pressure finite element space. As above, $\kappa_h=\text{Id} + \pi_h$, and $\delta^{(k)}$ as well as $\pi_h$ are defined as in \eqref{LPS_stab_matrix}. Pre- and post-multiplying by $B (\mathbf{L}^{(k)})^{-1} \mathbf{M}$ and $(L_p^{(k)})^{-1} M_p$, the previous expression then gives
\begin{displaymath}
B \mathbf{M}^{-1} B^\top (L_p^{(k)})^{-1} M_p \approx B (\mathbf{L}^{(k)})^{-1} B^\top.
\end{displaymath}
We still have no practical preconditioner due to the matrix $B \mathbf{M}^{-1} B^\top$; however, it can be proved that $K_p \approx B \mathbf{M}^{-1} B^\top$ for problems with enclosed flow \cite[p.\,176--177]{Elman_Silvester_Wathen}. Finally, a good approximation of the Schur complement $S=B (\mathbf{L}^{(k)})^{-1} B^\top$ is
\begin{displaymath}
\widehat{S} = K_p (L_p^{(k)})^{-1} M_p \approx S.
\end{displaymath}
Note that in our derivation we have also included the stabilization matrices on the velocity and the pressure space, which was not done in \cite{Kay_Loghin_Wathen}.

In the following we present a generalization of the pressure convection--diffusion preconditioner, applying the commutator argument in \eqref{forward_commutator} to the case where the differential operators involved are to be considered vectorial differential operators, i.e.
\begin{equation}\label{general_commutator}
\mathcal{E}_n = \mathcal{D} \nabla_n - \nabla_n \mathcal{D}_p,
\end{equation}
where
\begin{displaymath}
\mathcal{D} = \left[
\begin{array}{ccc}
\mathcal{D}^{1,1} & \ldots & \mathcal{D}^{1,n}\\
\vdots & \ddots & \vdots\\
\mathcal{D}^{n,1} & \ldots & \mathcal{D}^{n,n}
\end{array}
\right], \qquad
\mathcal{D}_p = \left[
\begin{array}{ccc}
\mathcal{D}_p^{1,1} & \ldots & \mathcal{D}_p^{1,n}\\
\vdots & \ddots & \vdots\\
\mathcal{D}_p^{n,1} & \ldots & \mathcal{D}_p^{n,n}
\end{array}
\right].
\end{displaymath}
Here $\mathcal{D}^{i,j}$ is a differential operator on the velocity space with $D_p^{i,j}$ the corresponding differential operator on the pressure space, for $i,j=1,\ldots,n$, and $\nabla_n = I_n \otimes \nabla$, with $I_n\in\mathbb{R}^{n\times n}$ the identity matrix for some $n\in \mathbb{N}$. As above, we suppose that each $D_p^{i,j}, \: i,j=1,\ldots,n$, is well defined, and that the commutator $\mathcal{E}_n$ is small in some sense. Again, after discretizing with stable finite elements we can rewrite
\begin{equation}\label{discretized_general_commutator}
\left(\mathcal{M}^{-1} \mathbf{D}\right) \mathcal{M}^{-1} \vec{B}^{\: \top} - \mathcal{M}^{-1} \vec{B}^{\: \top} \left( \mathcal{M}_p^{-1} D_p \right) \approx 0,
\end{equation}
where $\mathcal{M}= I_n \otimes \mathbf{M}$, $\mathcal{M}_p= I_n \otimes M_p$, $\vec{B}= I_n \otimes B$, and
\begin{displaymath}
\mathbf{D} = \left[
\begin{array}{ccc}
\mathbf{D}^{1,1} & \ldots & \mathbf{D}^{1,n}\\
\vdots & \ddots & \vdots\\
\mathbf{D}^{n,1} & \ldots & \mathbf{D}^{n,n}
\end{array}
\right], \qquad
D_p = \left[
\begin{array}{ccc}
D_p^{1,1} & \ldots & D_p^{1,n}\\
\vdots & \ddots & \vdots\\
D_p^{n,1} & \ldots & D_p^{n,n}
\end{array}
\right],
\end{displaymath}
with $\mathbf{D}^{i,j}$ and $D^{i,j}_p$ the corresponding discretizations of $\mathcal{D}^{i,j}$ and $\mathcal{D}^{i,j}_p$, respectively. Pre-multiplying \eqref{discretized_general_commutator} by $\vec{B} \mathbf{D}^{-1} \mathcal{M}$, and post-multiplying by $D_p^{-1} \mathcal{M}_p$, gives that
\begin{displaymath}
\vec{B} \mathcal{M}^{-1} \vec{B}^{\: \top} D_p^{-1} \mathcal{M}_p \approx \vec{B} \mathbf{D}^{-1} \vec{B}^{\: \top}.
\end{displaymath}
Noting that $\vec{B} \mathcal{M}^{-1} \vec{B}^{\: \top} = I_n \otimes (B \mathbf{M}^{-1} B^\top)$ and recalling that $K_p \approx B \mathbf{M}^{-1} B^\top$, we derive the following approximation:
\begin{equation}\label{commutator_approx}
\mathcal{K}_pD_p^{-1} \mathcal{M}_p \approx \vec{B} \mathbf{D}^{-1} \vec{B}^{\: \top},
\end{equation}
where $\mathcal{K}_p = I_n \otimes K_p$. In Section \ref{5_2} we employ the approach outlined here to devise preconditioners for discrete optimality conditions of Navier--Stokes control problems.

\section{First-Order Optimality Conditions and Time-Stepping}
\label{4}
We now describe the strategy used for obtaining an approximate solution of \eqref{Stationary_Navier_Stokes_control_functional}--\eqref{Stationary_Navier_Stokes_control_constraints} and of \eqref{Instationary_Navier_Stokes_control_functional}--\eqref{Instationary_Navier_Stokes_control_constraints}. We introduce \emph{adjoint variables} (or \emph{Lagrange multipliers}) $\vec{\zeta}$ and $\mu$ and make use of an optimize-then-discretize scheme, stating the first-order optimality conditions. We then discretize the conditions so obtained, for both the stationary and instationary Navier--Stokes control problems, and derive the corresponding Oseen linearized problems. For the instationary problem \eqref{Instationary_Navier_Stokes_control_functional}--\eqref{Instationary_Navier_Stokes_control_constraints}, we consider employing both backward Euler and Crank--Nicolson schemes in time. \color{black} Both time-stepping schemes are A--stable, hence avoiding any constraints on the time step used. While only first-order accurate, backward Euler is also L--stable, and the technique presented below is easily generalized when the initial condition $\vec{v}_0$ is not solenoidal. On the other hand, Crank--Nicolson is not L--stable, but is second-order accurate. However, if $\vec{v}_0$ is not solenoidal, pre-processing is required in order to write the Oseen iteration\color{black}.

\subsection{Stationary Navier--Stokes control}
\label{4_1}
\color{black}Introducing the adjoint velocity $\vec{\zeta}$ and the adjoint pressure $\mu$, we may consider the Lagrangian associated with \eqref{Stationary_Navier_Stokes_control_functional}--\eqref{Stationary_Navier_Stokes_control_constraints} as in \cite{Posta_Roubicek}, and write the Karush--Kuhn--Tucker (KKT) conditions as:
\begin{equation}\label{Stationary_Navier_Stokes_control_optimality_conditions}
\left\{
\begin{array}{rl}
\vspace{0.25ex}
\left.
\begin{array}{rl}
 - \nu\nabla^2 \vec{v} + \vec{v} \cdot \nabla \vec{v} + \nabla p = \frac{1}{\beta}\: \vec{\zeta} + \vec{f} & \quad \mathrm{in} \; \Omega \\
- \nabla \cdot \vec{v}(\mathbf{x}) = 0 & \quad \mathrm{in} \; \Omega \\
\vec{v}(\mathbf{x}) = \vec{g}(\mathbf{x}) & \quad \mathrm{on} \; \partial \Omega
\end{array}
 \right\} & \!
\left.
\begin{array}{c}
\vspace{0.5ex}
\mathrm{state}\\
\mathrm{equations}
\end{array}
\right.

\\

\left.
\begin{array}{rl}
- \nu\nabla^2 \vec{\zeta} - \vec{v} \cdot \nabla \vec{\zeta} + (\: \nabla \vec{v} \:) ^\top \vec{\zeta} + \nabla \mu \: = \: \vec{v}_d - \vec{v}  & \quad \mathrm{in} \; \Omega  \\
- \nabla \cdot \vec{\zeta}(\mathbf{x}) = 0 & \quad \mathrm{in} \; \Omega \\
\vec{\zeta}(\mathbf{x}) = \vec{0} & \quad \mathrm{on} \; \partial \Omega
\end{array}
\right\} & \!
\left.
\begin{array}{c}
\vspace{0.5ex}
\mathrm{adjoint}\\
\mathrm{equations}
\end{array}
\right.

\end{array}
\right.
\end{equation}
where we have substituted the gradient equation $\beta \vec{u} - \vec{\zeta} = 0$ into the state equation\color{black}.

Problem \eqref{Stationary_Navier_Stokes_control_optimality_conditions} is a coupled system of non-linear, stationary PDEs. In order to find a numerical solution of \eqref{Stationary_Navier_Stokes_control_optimality_conditions}, we need to solve a sequence of linearizations of the system. As in \cite{Posta_Roubicek}, we solve at each step the Oseen approximation as follows. Letting $\vec{v}^{\:( k)}\in V, \: p^{(k)} \in Q, \: \vec{\zeta}^{\: (k)} \in V_0, \: \mu^{(k)} \in Q$ be the current approximations to $\vec{v}$, $p$, $\vec{\zeta}$, and $\mu$, respectively, with $V, \: V_0, \: Q$ defined as in Section \ref{2_1}, the Oseen iterate is defined as
\begin{equation}\label{vplusdv}
\begin{array}{c}
\vec{v}^{\: (k+1)}=\vec{v}^{\: (k)}+ \vec{\delta v}^{\: (k)}, \qquad p^{(k+1)} = p^{(k)} + \delta p ^{(k)}, \\
 \vec{\zeta}^{\: (k+1)}= \vec{\zeta}^{\: (k)} + \vec{\delta \zeta}^{\: (k)}, \qquad \mu^{(k+1)}=\mu^{(k)} + \delta \mu^{(k)},
\end{array}
\end{equation}
with $\vec{\delta v}^{\: (k)}$, $\delta p ^{(k)}$, $\vec{\delta \zeta}^{\: (k)}$, $\delta \mu^{(k)}$ the solutions of the following Oseen problem:
\begin{equation}\label{Stationary_Oseen_problem}
\!
\left\{
\begin{array}{rl}
\!\! \nu (\nabla \vec{\delta v}^{\: (k)}, \nabla \vec{w} ) + ( \vec{v}^{\: (k)} \cdot \nabla \vec{\delta v}^{\: (k)}, \vec{w}) - ( \delta p^{(k)}, \nabla \cdot \vec{w} ) - \frac{1}{\beta} ( \vec{\delta \zeta}^{\: (k)},\vec{w}) = &  \hspace{-0.6em} \vec{R}_1^{\: (k)}, \\
- (q,\nabla \cdot \vec{\delta v}^{\: (k)}) = & \hspace{-0.6em} r_1^{(k)}, \\

(\vec{\delta v}^{\: (k)},\vec{w}) + \nu(\nabla \vec{\delta \zeta}^{\: (k)},\nabla \vec{w}) - (\vec{v}^{\: (k)} \cdot \nabla \vec{\delta \zeta}^{\: (k)},\vec{w}) - ( \delta \mu^{(k)} , \nabla \cdot \vec{w})  = & \hspace{-0.6em} \vec{R}_2^{\: (k)},  \\
- (q,\nabla \cdot \vec{\delta \zeta}^{\: (k)}) = & \hspace{-0.6em} r_2^{(k)},
\end{array}
\right.
\end{equation}
for any $\vec{w} \in V_0$ and $q\in Q$. The residuals $\vec{R}_1^{\: (k)},\: r_1^{(k)},\: \vec{R}_2^{\: (k)},\: r_2^{(k)}$ are given by
\begin{equation*}\label{Stationary_Oseen_residuals}
\!\!
\left\{
\!\!
\begin{array}{rl}
\vspace{0.125ex}
\vec{R}_1^{\: (k)} \! = & \hspace{-0.8em} (\vec{f},\vec{w}) \! - \! \nu(\nabla \vec{v}^{\: (k)},\nabla \vec{w}) \! - \! (\vec{v}^{\: (k)} \cdot \nabla \vec{v}^{\: (k)}, \vec{w}) \! + \! ( p^{(k)}, \nabla \cdot \vec{w}) \! + \! \frac{1}{\beta}( \vec{ \zeta}^{\: (k)},\vec{w}), \\
\vspace{0.125ex}
r_1^{(k)} \! = & \hspace{-0.8em} (q,\nabla \cdot \vec{v}^{\: (k)}), \\
\vec{R}_2^{\: (k)} \! = & \hspace{-0.8em} (\vec{v}_d,\vec{w}) - (\vec{v}^{\: (k)},\vec{w}) - \nu (\nabla \vec{\zeta}^{\: (k)},\nabla \vec{w}) + (\vec{v}^{\: (k)} \cdot \nabla \vec{\zeta}^{\: (k)},\vec{w}) \\
\vspace{0.125ex}
& - ((\: \nabla \vec{v}^{\: (k)} \:) ^\top \vec{\zeta}^{\: (k)},\vec{w}) + ( \mu^{(k)}, \nabla \cdot \vec{w}), \\
r_2^{(k)} \! = & \hspace{-0.8em} (q, \nabla \cdot \vec{\zeta}^{\:( k)}).
\end{array}
\right.
\end{equation*}

The Oseen problem \eqref{Stationary_Oseen_problem} is posed on the continuous level, so we need to discretize it in order to obtain a numerical solution of \eqref{Stationary_Navier_Stokes_control_functional}--\eqref{Stationary_Navier_Stokes_control_constraints}. Let $\mathbf{v}^{(k)}=\{\mathbf{v}^{(k)}_i\}_{i=1}^{n_v}, \: \mathbf{p}^{(k)}=\{p^{(k)}_i\}_{i=1}^{n_p}, \: \bm{\zeta}^{(k)}=\{\bm{\zeta}^{(k)}_i\}_{i=1}^{n_v}, \: \bm{\mu}^{(k)}=\{\mu^{(k)}_i\}_{i=1}^{n_p}$ be the vectors containing the numerical solutions at the $k$-th iteration for $\vec{v}^{\: (k)}$, $p^{(k)}$, $\vec{\zeta}^{\; (k)}$, and $\mu^{(k)}$, respectively, that is, $\vec{v}^{\: (k)}\approx \sum_{i=1}^{n_v}\mathbf{v}^{(k)}_i \vec{\phi}_i$, $p^{\: (k)}\approx \sum_{i=1}^{n_p}\mathbf{p}^{(k)}_i \psi_i$, $\vec{\zeta}^{\: (k)}\approx \sum_{i=1}^{n_v}\bm{\zeta}^{(k)}_i \vec{\phi}_i$, $\mu^{\: (k)}\approx \sum_{i=1}^{n_p}\bm{\mu}^{(k)}_i \psi_i$. Then, the (discrete) Oseen iterate is defined as
\begin{displaymath}
\begin{array}{c}
\mathbf{v}^{\: (k+1)}=\mathbf{v}^{\: (k)}+ \bm{\delta v}^{\: (k)}, \qquad \bm{p}^{(k+1)} = \bm{p}^{(k)} + \bm{\delta p }^{(k)}, \\
 \bm{\zeta}^{\: (k+1)}= \bm{\zeta}^{\:( k)} + \bm{\delta \zeta}^{\:( k)}, \qquad \bm{\mu}^{(k+1)}=\bm{\mu}^{(k)} + \bm{\delta \mu}^{(k)},
\end{array}
\end{displaymath}
where
\begin{equation}\label{Stationary_discrete_Oseen_problem}
\left\{
\begin{array}{rl}
\mathbf{L}^{(k)} \: \bm{\delta v}^{\: (k)} + B^\top \: \bm{\delta p}^{(k)} -  \mathbf{M}_{\beta}\: \bm{\delta \zeta}^{\: (k)} = & \hspace{-0.6em} \bm{R}_1^{\: (k)}, \\
B \bm{\delta v}^{\: (k)} = & \hspace{-0.6em} \bm{r}_1^{(k)}, \\
\mathbf{M} \:\bm{\delta v}^{\: (k)} + \mathbf{L}^{(k)}_{\mathrm{adj}} \: \bm{\delta \zeta}^{\: (k)} + B^\top \bm{\delta \mu}^{(k)} \:  = & \hspace{-0.6em} \bm{R}_2^{\: (k)},\\
B \bm{\delta \zeta}^{\: (k)} = & \hspace{-0.6em} \bm{r}_2^{(k)},

\end{array}
\right.
\end{equation}
with $\mathbf{M}_\beta=\frac{1}{\beta}\mathbf{M}$, $\mathbf{L}^{(k)}_\mathrm{adj}= \nu \mathbf{K} - \mathbf{N}^{(k)} + \mathbf{W}^{(k)}$, and the discrete residuals given by
\begin{displaymath}
\left\{
\begin{array}{rll}
\bm{R}_1^{\: (k)} = & \hspace{-0.6em} \mathbf{f} - \mathbf{L}^{(k)} \: \mathbf{v}^{\: (k)} - B^\top \: \mathbf{p}^{(k)} + \mathbf{M}_\beta\: \bm{ \zeta}^{\: (k)}, \\
\bm{r}_1^{(k)} = & \hspace{-0.6em} - B \mathbf{v}^{\: (k)}, \\
\bm{R}_2^{\: (k)} = & \hspace{-0.6em} \mathbf{M} \: \mathbf{v}_d - \mathbf{M} \:  \mathbf{v}^{\: (k)} - \mathbf{L}^{(k)}_\mathrm{adj} \: \bm{\zeta}^{\: (k)} - B^\top \bm{\mu}^{(k)}  - \bm{\omega}^{(k)},\\
\bm{r}_2^{(k)} = & \hspace{-0.6em} - B \bm{\zeta}^{\: (k)}.
\end{array}
\right.
\end{displaymath}
Here $\mathbf{v}_d$ is the vector corresponding to the discretized desired state $\vec{v}_d$, and $\bm{\omega}^{(k)}=\{ \big( (\nabla \vec{v}^{\:(k)} )^\top   \vec{\zeta}^{\: (k)}, \vec{\phi}_i \big) \}_{i=1}^{n_v}$.
In our tests, the initial guesses $\boldsymbol{v}^{(1)}$ and $\boldsymbol{\zeta}^{(1)}$ for the non-linear process are the state and adjoint velocity solutions of the KKT conditions for the corresponding stationary Stokes control problem, with discretization given by \eqref{Stationary_discrete_Oseen_problem} with $\mathbf{L}^{(k)}=\mathbf{L}^{(k)}_\mathrm{adj}=  \mathbf{K} $, and residuals $\bm{R}_1^{\: (k)}=\mathbf{f}$, $ \bm{R}_2^{\: (k)}=\mathbf{M} \: \mathbf{v}_d$, $\bm{r}_1^{(k)}= \bm{r}_2^{(k)}=\mathbf{0}$. Note that the right-hand side may also take into account boundary conditions (as done in our implementation).

In matrix form, we rewrite system \eqref{Stationary_discrete_Oseen_problem} as
\begin{equation}\label{Stationary_discrete_Oseen_system}
\underbrace{
\left[
\begin{array}{cc}
\Phi^{(k)}_{\mathrm{S}} & \Psi_{\mathrm{S}}^\top\\
\Psi_{\mathrm{S}} & - \Theta_{\mathrm{S}}
\end{array}
\right]
}_{\mathcal{A}^{(k)}_{\mathrm{S}}}
\left[
\begin{array}{c}
\bm{\delta v}^{(k)}\\
\bm{\delta \zeta}^{(k)}\\
\bm{\delta \mu}^{(k)}\\
\bm{\delta p}^{(k)}
\end{array}
\right]=
\left[
\begin{array}{c}
\bm{R}_2^{(k)}\\
\bm{R}_1^{(k)}\\
\bm{r}_1^{(k)}\\
\bm{r}_2^{(k)}
\end{array}
\right],
\end{equation}
where
\begin{equation}\label{Phi_Psi_Theta_Stationary}
\Phi_{\mathrm{S}}^{(k)}=
\left[
\begin{array}{cc}
\mathbf{M} & \mathbf{L}^{(k)}_{\mathrm{adj}}\\
\mathbf{L}^{(k)} & -\mathbf{M}_\beta
\end{array}
\right], \quad
\Psi_{\mathrm{S}}=
\left[
\begin{array}{cc}
B & 0\\
0 & B
\end{array}
\right], \quad
\Theta_{\mathrm{S}}=
\left[
\begin{array}{cc}
0 & 0\\
0 & 0
\end{array}
\right].
\end{equation}
The matrix $\mathcal{A}_{\mathrm{S}}^{(k)}$ is of saddle-point type; however, since the incompressibility constraints $\nabla \cdot \vec{v}=0$ are not solved exactly, $\Phi^{(k)}_\mathrm{S}$ is not symmetric in general.

\subsection{Instationary Navier--Stokes control}
\label{4_2}
We now state the KKT conditions for the instationary problem \eqref{Instationary_Navier_Stokes_control_functional}--\eqref{Instationary_Navier_Stokes_control_constraints}. As before, introducing the adjoint variables $\vec{\zeta}$ and $\mu$, we consider the Lagrangian associated to \eqref{Instationary_Navier_Stokes_control_functional}--\eqref{Instationary_Navier_Stokes_control_constraints} as in \cite[p.\,318]{Troltzsch}. Then, by deriving the KKT conditions and substituting the gradient equation $\beta \vec{u} - \vec{\zeta} = 0$ into the state equation, the solution of \eqref{Instationary_Navier_Stokes_control_functional}--\eqref{Instationary_Navier_Stokes_control_constraints} satisfies:
\begin{equation}\label{Instationary_Navier_Stokes_control_optimality_conditions}
\left\{
\begin{array}{rl}

\vspace{0.125ex}
\frac{\partial \vec{v}}{\partial t} - \nu\nabla^2 \vec{v} + \vec{v} \cdot \nabla \vec{v} + \nabla p = \frac{1}{\beta}\: \vec{\zeta} + \vec{f} & \quad \mathrm{in} \; \Omega \times (0,t_{f}), \\
\vspace{0.125ex}
- \nabla \cdot \vec{v}(\mathbf{x},t) = 0 & \quad \mathrm{in} \; \Omega \times (0,t_{f}), \\
\vspace{0.125ex}
\vec{v}(\mathbf{x},t) = \vec{g}(\mathbf{x},t) & \quad \mathrm{on} \; \partial \Omega \times (0,t_{f}), \\
\vspace{0.25ex}
\vec{v}(\mathbf{x},0) = \vec{v}_0(\mathbf{x}) & \quad \mathrm{in} \; \Omega, \\

\vspace{0.125ex}
-\frac{\partial \vec{\zeta}}{\partial t} - \nu\nabla^2 \vec{\zeta} - \vec{v} \cdot \nabla \vec{\zeta} + (\: \nabla \vec{v} \: ) ^\top \vec{\zeta} + \nabla \mu \: = \: \vec{v}_d - \vec{v}  & \quad \mathrm{in} \; \Omega \times (0,t_{f}), \\
\vspace{0.125ex}
- \nabla \cdot \vec{\zeta}(\mathbf{x},t) = 0 & \quad \mathrm{in} \; \Omega \times (0,t_{f}), \\
\vspace{0.125ex}
\vec{\zeta}(\mathbf{x},t) = \vec{0} & \quad \mathrm{on} \; \partial \Omega \times (0,t_{f}), \\
\vec{\zeta}(\mathbf{x},t_{f})= \vec{0} & \quad \mathrm{in} \; \Omega.

\end{array}
\right.
\end{equation}

Problem \eqref{Instationary_Navier_Stokes_control_optimality_conditions} is a coupled system of non-linear, instationary PDEs. In order to find a numerical solution of \eqref{Instationary_Navier_Stokes_control_optimality_conditions}, as for the stationary case we take an Oseen linearization. Let $\vec{v}^{\:( k)}\in \bar{V}, \, p^{(k)} \in \bar{Q}, \, \vec{\zeta}^{\: (k)} \in \bar{V}_0, \, \mu^{(k)} \in \bar{Q}$ be the current approximation to $\vec{v},\, p, \, \vec{\zeta}, \, \mu$, with $\bar{V}:=\{ \vec{v} \in L^2( 0,t_f; \mathcal{H}^1(\Omega)^d) \, | \, \frac{\partial \vec{v}}{\partial t}(\cdot, t) \in L^2(0,t_f; \mathcal{H}^{-1}(\Omega)^d) \; \mathrm{for \; a.e.} \; t \in (0,t_f), \, \vec{v}=\vec{g} \; \mathrm{on} \; \partial \Omega, \, \vec{v}(\mathbf{x},0) = \vec{v}_0(\mathbf{x}) \}, $ $\bar{Q}= L^2( 0,t_f; L^2(\Omega))$, and $\bar{V}_0$ the corresponding space for the adjoint velocity (see \cite[p.\,315--321]{Troltzsch} and \cite[p.\,88--95]{Hinze}, for instance, for the case $d=2$). Then, the Oseen iterate is of the form \eqref{vplusdv}, with $\vec{\delta v}^{\: (k)}, \, \delta p ^{(k)}, \, \vec{\delta \zeta}^{\: (k)}, \, \delta \mu^{(k)}$ the solution of: 
\begin{equation}\label{Instationary_Oseen_problem}
\!\!
\left\{
\begin{array}{rl}
\!\!\! \frac{\partial }{\partial t}(\vec{\delta v}^{\! \:(k)} \! ,\vec{w}) \! + \! \nu (\nabla \vec{\delta v}^{\! \:(k)} \! ,\nabla \vec{w}) \! + \! (\vec{v}^{\! \:(k)} \! \cdot \! \nabla \vec{\delta v}^{\! \:(k)} \! ,\vec{w}) \! 
- \! ( \delta p^{(k)} \! ,\nabla \! \cdot \! \vec{w}) \! - \! \frac{1}{\beta} ( \vec{\delta \zeta}^{\! \:(k)} \! ,\vec{w}) \! = & \hspace{-0.6em}\! \vec{R}_1^{(k)}, \\
- \nabla \cdot \vec{\delta v}^{\! \:(k)} \! = & \hspace{-0.6em}\! r_1^{(k)}, \\

-\frac{\partial }{\partial t}(\vec{\delta \zeta}^{\! \:(k)} \! ,\vec{w}) \! + \! \nu (\nabla \vec{\delta \zeta}^{\! \:(k)} \! ,\vec{w}) \! - \! (\vec{v}^{\! \:(k)} \! \cdot \! \nabla \vec{\delta \zeta}^{\! \:(k)} \! ,\vec{w})
\! - \! ( \delta \mu^{(k)} \! ,\nabla \! \cdot \! \vec{w}) \! + \! (\vec{\delta v}^{\! \:(k)} \! ,\vec{w}) \! = & \hspace{-0.6em}\! \vec{R}_2^{(k)} ,\\

- \nabla \cdot \vec{\delta \zeta}^{\! \:(k)} \! = & \hspace{-0.6em}\! r_2^{(k)}, 
\end{array}
\right.
\end{equation}
for any $\vec{w} \in V_0$ and $q\in Q$. The residuals $\vec{R}_1^{(k)},\: r_1^{(k)},\: \vec{R}_2^{(k)}, \: r_2^{(k)}$ are given by
\begin{equation}\label{Instationary_Oseen_residuals}
\left\{
\begin{array}{rl}

\vec{R}_1^{(k)} =& \hspace{-0.6em} (\vec{f},\vec{w}) - \frac{\partial }{\partial t}(\vec{v}^{\:(k)},\vec{w}) - \nu (\nabla \vec{v}^{\:(k)},\nabla \vec{w}) - (\vec{v}^{\:(k)} \cdot \nabla \vec{v}^{\:(k)},\vec{w})\\
 & + ( p^{(k)},\nabla \cdot \vec{w}) + \frac{1}{\beta}\: (\vec{\zeta}^{\:(k)} , \vec{w}),  \\
r_1^{(k)} =& \hspace{-0.6em} (q,\nabla \cdot \vec{v}^{\:(k)}), \\
\vec{R}_2^{(k)} =& \hspace{-0.6em} (\vec{v}_d,\vec{w}) - (\vec{v}^{\:(k)},\vec{w}) + \frac{\partial }{\partial t}(\vec{\zeta}^{\:(k)},\vec{w}) + \nu (\nabla \vec{\zeta}^{\:(k)},\nabla \vec{w})\\
& + (\vec{v}^{\:(k)} \cdot \nabla \vec{\zeta}^{\:(k)},\vec{w}) - ((\nabla \vec{v}^{\:(k)}) ^\top \vec{\zeta}^{\:(k)},\vec{w}) + ( \mu^{(k)},\nabla \cdot \vec{w}), \\
r_2^{(k)} =& \hspace{-0.6em} (q,\nabla \cdot \vec{\zeta}^{\:(k)} ).
\end{array}
\right.
\end{equation}
Note that with this notation $\vec{\delta v}^{\:(k)}(\mathbf{x},0)=\vec{\delta \zeta}^{\:(k)}(\mathbf{x},t_f)=\vec{0}$ in $\Omega$, and $\vec{\delta v}^{\:(k)}(\mathbf{x},t)=\vec{\delta \zeta}^{\:(k)}(\mathbf{x},t)=\vec{0}$ on $\partial \Omega \times (0,t_f)$.

Equations \eqref{Instationary_Oseen_problem} are the Oseen equations for instationary Navier--Stokes control, involving a coupled system of instationary convection--diffusion equations and divergence-free conditions. Due to this structure, we present two discretized version of \eqref{Instationary_Oseen_problem}, one making use of backward Euler time-stepping, the other employing the Crank--Nicolson scheme. Further, in order to solve \eqref{Instationary_Oseen_problem} we need to choose an initial guess $\boldsymbol{v}^{(1)}$ and $\boldsymbol{\zeta}^{(1)}$ for the state and the adjoint velocities and then iteratively solve a sequence of linearized problems. In our tests, $\boldsymbol{v}^{(1)}$ and $\boldsymbol{\zeta}^{(1)}$ are again the (velocity) solutions of the KKT conditions for the corresponding Stokes control problem, for which the optimality conditions are:
\begin{equation*}
\left\{
\begin{array}{rl}

\frac{\partial \vec{v}}{\partial t} - \nabla^2 \vec{v} + \nabla p = \frac{1}{\beta}\: \vec{\zeta} + \vec{f} & \quad \mathrm{in} \; \Omega \times (0,t_{f}), \\

-\frac{\partial \vec{\zeta}}{\partial t} - \nabla^2 \vec{\zeta} + \nabla \mu \: = \: \vec{v}_d - \vec{v}  & \quad \mathrm{in} \; \Omega \times (0,t_{f}), \\

\end{array}
\right.
\end{equation*}
along with the incompressibility constraints together with initial, final, and boundary conditions on the state and adjoint velocities as in \eqref{Instationary_Navier_Stokes_control_optimality_conditions}.

We now derive the linear systems resulting from the time-stepping schemes. For the sake of exposition, we introduce the matrices $I_{n,1}$, $I_{n,2}$, $I_{n,3}$, $I_{n,4} \in \mathbb{R}^{n\times n}$, with
\begin{displaymath}
\begin{array}{rl}
\begin{array}{rl}
\vspace{1ex}
I_{n,1}= & \!\! \texttt{diag}(1,\ldots ,1,0),
\\
\vspace{1ex}
I_{n,2}= & \!\! \texttt{diag}(0,1,\ldots ,1),\\
I_{n,4}= & \!\! I_{n}+I_{n,3},
\end{array}
&
\quad I_{n,3}=
\left[
\begin{array}{cccc}
0 & 1 \\
& \ddots & \ddots \\
& & 0 & 1\\
& & & 0
\end{array}
\right].
\end{array}
\end{displaymath}
Here, \texttt{diag} denotes a diagonal matrix with the diagonal entries stated.

\subsubsection{Backward Euler for instationary Navier--Stokes control}\label{4_2_1}
In this section we introduce the backward Euler scheme for approximating \eqref{Instationary_Oseen_problem}--\eqref{Instationary_Oseen_residuals}, and then derive the resulting linear system. We discretize the interval $(0,t_f)$ into $n_t$ subintervals of length $\tau=\frac{t_f}{n_t}$, denoting the grid points as $t_n=n\tau$, for $n = 0, \ldots,n_t$. We approximate all the functions on this time grid, excluding the initial and final time points for the state and adjoint pressure, respectively. Specifically, our approximations of the solutions at the $k$-th step of the non-linear solver are given by $\boldsymbol{v}_n^{(k)} \approx \vec{v}(\mathbf{x},t_n)$, $\boldsymbol{\zeta}_n^{(k)} \approx \vec{\zeta}(\mathbf{x},t_n)$, for $n=0,\ldots,n_t$, and $\boldsymbol{p}^{(k)}_{n+1} \approx p(\mathbf{x},t_{n+1})$, $\boldsymbol{\mu}^{(k)}_{n} \approx \mu(\mathbf{x},t_n)$, for $n=0,\ldots,n_t-1$, for all $\mathbf{x} \in \Omega$. We also introduce the following finite element matrices:
\begin{equation*}
\begin{array}{c}
\mathbf{L}^{(k)}_n =  \tau (\nu \mathbf{K} + \mathbf{N}^{(k)}_n + \mathbf{W}_n^{(k)}) + \mathbf{M}, \qquad \mathbf{T}^{(k)}_n = \tau (\nu \mathbf{K} - \mathbf{N}^{(k)}_n + \mathbf{W}_n^{(k)}) + \mathbf{M},\\
\bar{\mathbf{M}}^\mathrm{BE}= \tau \mathbf{M}, \qquad \bar{\mathbf{M}}^\mathrm{BE}_\beta=\frac{\tau}{\beta} \mathbf{M}, \qquad \bar{B}= \tau B,
\end{array}
\end{equation*}
where $\mathbf{W}_n^{(k)}$ is the stabilization matrix related to $\vec{v}^{\:(k)}_n$, and $\mathbf{N}^{(k)}_n = [( \vec{v}^{\:(k)}_n \cdot \nabla \vec{\phi}_j, \vec{\phi}_i )]$,
with $\vec{v}^{\:(k)}_n$ the approximation to $\vec{v}$ at time $t_n$, at the $k$-th Oseen iteration. Note that the superscripts of $\mathbf{L}^{(k)}_0$, $\mathbf{T}^{(k)}_0$ are superfluous, as the initial condition on $\vec{v}$ is fixed; however we keep them for consistency. We then write the discrete Oseen iterate as
\begin{displaymath}
\begin{array}{lll}
\boldsymbol{v}_n^{(k+1)} = \boldsymbol{v}_n^{(k)} + \boldsymbol{\delta v}_n^{(k)}, & \boldsymbol{\zeta}_n^{(k+1)} = \boldsymbol{\zeta}_n^{(k)} + \boldsymbol{\delta \zeta}_n^{(k)}, & \quad n = 0 , \ldots, n_t,\\
\boldsymbol{p}_{n+1}^{(k+1)} = \boldsymbol{p}_{n+1}^{(k)} + \boldsymbol{\delta p}_{n+1}^{(k)}, & \boldsymbol{\mu}_n^{(k+1)} = \boldsymbol{\mu}_n^{(k)} + \boldsymbol{\delta \mu}_n^{(k)}, & \quad n = 0 , \ldots, n_t-1,
\end{array}
\end{displaymath}
with $\boldsymbol{\delta v}_n^{(k)}$, $\boldsymbol{\delta \zeta}_n^{(k)},$ $\boldsymbol{\delta p}_n^{(k)}$, $\boldsymbol{\delta \mu}_n^{(k)}$ the solutions of the following discretization of \eqref{Instationary_Oseen_problem}:
\begin{equation}\label{Instationary_discrete_Oseen_problem_BE}
\left\{
\begin{array}{rll}
\vspace{0.25ex}
\displaystyle 
\bar{\mathbf{M}}^\mathrm{E} \boldsymbol{\delta v}^{(k)}_n + \mathbf{T}^{(k)}_n \boldsymbol{\delta \zeta}^{(k)}_n - \mathbf{M}\: \boldsymbol{\delta \zeta}^{(k)}_{n+1} + \bar{B}^{\top} \boldsymbol{\delta \mu}^{(k)}_{n} = & \hspace{-0.6em} \boldsymbol{R}_{2,n}^{(k)} , \\
\vspace{0.25ex}
\displaystyle - \mathbf{M} \: \boldsymbol{\delta v}^{(k)}_n + \mathbf{L}^{(k)}_{n+1} \boldsymbol{\delta v}^{(k)}_{n+1} + \bar{B}^{\top} \boldsymbol{\delta p}^{(k)}_{n+1} - \bar{\mathbf{M}}^\mathrm{E}_\beta  \boldsymbol{\delta \zeta}^{(k)}_{n+1} = & \hspace{-0.6em} \boldsymbol{R}_{1,n}^{(k)},\\
\vspace{0.25ex}
B \: \boldsymbol{\delta v}^{(k)}_{n+1} = & \hspace{-0.6em} \boldsymbol{r}_{1,n+1}^{(k)}, \\
B \: \boldsymbol{\delta \zeta}^{(k)}_{n} = & \hspace{-0.6em} \boldsymbol{r}_{2,n}^{(k)},
\end{array}
\right.
\end{equation}
for $n=0,...,n_t-1$, with $\boldsymbol{\delta v}^{(k)}_0=\boldsymbol{0}$, $\boldsymbol{\delta \zeta}^{(k)}_{n_t}=\boldsymbol{0}$. The discretized residuals are given by
\begin{equation}\label{Instationary_discrete_Oseen_residual_BE}
\left\{
\begin{array}{rll}
\vspace{0.25ex}
\boldsymbol{R}^{(k)}_{1,n} = & \hspace{-0.6em} \tau \boldsymbol{f}^{n+1} + \mathbf{M} \boldsymbol{ v}^{(k)}_{n} - \mathbf{L}^{(k)}_{n+1} \boldsymbol{ v}^{(k)}_{n+1} - \bar{B}^{\top} \boldsymbol{ p}^{(k)}_{n+1} + \bar{\mathbf{M}}^\mathrm{E}_\beta \boldsymbol{ \zeta}^{(k)}_{n+1}, \\
\vspace{0.25ex}
\boldsymbol{r}^{(k)}_{1,n+1} = & \hspace{-0.6em} - B \: \boldsymbol{ v}^{(k)}_{n+1}, \\
\vspace{0.25ex}
\boldsymbol{R}^{(k)}_{2,n} = & \hspace{-0.6em} \bar{\mathbf{M}}^\mathrm{E} \boldsymbol{v}_d^n - \bar{\mathbf{M}}^\mathrm{E}\boldsymbol{ v}^{(k)}_n  - \mathbf{T}^{(k)}_n \boldsymbol{ \zeta}^{(k)}_n + \mathbf{M} \boldsymbol{\zeta}^{(k)}_{n+1} - \bar{B}^{\top} \boldsymbol{ \mu}^{(k)}_{n} - \tau  \boldsymbol{\omega}^{(k)}_n, \\
\boldsymbol{r}^{(k)}_{2,n} = & \hspace{-0.6em} - B \: \boldsymbol{ \zeta}^{(k)}_{n},
\end{array}
\right.
\end{equation}
where $\boldsymbol{f}^{n+1}=\{( \vec{f}(\mathbf{x},t_{n+1}), \vec{\phi}_i )\}_{i=1}^{n_v}$, and $\boldsymbol{\omega}^{(k)}_n=\{ \big( (\nabla \vec{v}^{\:(k)}_n )^\top  \vec{\zeta}^{\: (k)}_n ,\vec{\phi}_i \big) \}_{i=1}^{n_v}$, for $n=0,\ldots, n_t-1$.
Note that the non-linear residuals $\boldsymbol{R}^{(k)}_{1,0}$, $ \boldsymbol{R}^{(k)}_{1,n_t-1}$, $\boldsymbol{R}^{(k)}_{2,0}$, $\boldsymbol{R}^{(k)}_{2,n_t-1}$ in \eqref{Instationary_discrete_Oseen_residual_BE} take into account the initial and the final conditions on $\vec{v}$ and $\vec{\zeta}$.

Even if the incompressibility constraints $B \: \boldsymbol{\delta v}^{(k)}_{n+1} =\mathbf{0} $ are solved exactly, for $n=0,\ldots,n_t-1$, at each Oseen iteration the system described in \eqref{Instationary_discrete_Oseen_problem_BE} is not symmetric, due the conditions $\boldsymbol{\delta v}^{(k)}_0=\boldsymbol{0}$ and $\boldsymbol{\delta \zeta}^{(k)}_{n_t}=\boldsymbol{0}$. However it can be made symmetric in the Stokes control setting, using the following projections onto the space of divergence-free functions (\emph{solenoidal projection}), as done in \cite{Hinze_Koster_Turek}, for instance. Given a vector $\vec{b}$, its solenoidal projection is defined as $\boldsymbol{b}$, with
\begin{equation}\label{solenoidal_projection}
\left\{
\begin{array}{rl}
\mathbf{L}^{(k)}_0\: \boldsymbol{b} + \bar{B}^\top \boldsymbol{\bar{p}} = & \hspace{-0.6em}  \mathbf{L}_{\vec{b}} \: \vec{b},\\
B \:\boldsymbol{b} = & \hspace{-0.6em} \boldsymbol{0},
\end{array}
\right.
\end{equation}
with $\mathbf{L}_{\vec{b}} = \tau (\nu \mathbf{K} + \mathbf{N}_{\vec{b}} \, + \mathbf{W}_{\vec{b}} \,) + \mathbf{M}$, $\mathbf{N}_{\vec{b}}$ and $ \mathbf{W}_{\vec{b}}$ being the vector-convection and stabilization matrices related to $\vec{b}$. As the vector $\boldsymbol{0}$ is clearly divergence-free, the condition $\boldsymbol{\delta v}^{(k)}_0=\boldsymbol{0}$ is equivalent to
\begin{equation}\label{solenoidal_projection_BE}
\left\{
\begin{array}{rl}
\vspace{0.25ex}
\mathbf{L}^{(k)}_0\: \boldsymbol{\delta v}^{(k)}_0 + \bar{B}^\top \boldsymbol{\delta p}^{(k)}_0 = & \hspace{-0.6em} \boldsymbol{0},\\
B\: \boldsymbol{\delta v}^{(k)}_0 = & \hspace{-0.6em} \boldsymbol{0}.
\end{array}
\right.
\end{equation}
Analogously, the condition $\boldsymbol{\delta \zeta}^{(k)}_{n_t}=\boldsymbol{0}$ is equivalent to
\begin{displaymath}
\left\{
\begin{array}{rl}
\vspace{0.25ex}
 \mathbf{T}^{(k)}_{n_t}\: \boldsymbol{\delta \zeta}^{(k)}_{n_t} + \bar{B}^\top \boldsymbol{\delta \mu}^{(k)}_{n_t} = & \hspace{-0.6em} \boldsymbol{0},\\
B\: \boldsymbol{\delta \zeta}^{(k)}_{n_t} = & \hspace{-0.6em} \boldsymbol{0}.
\end{array}
\right.
\end{displaymath}

By imposing the previous projections and multiplying the incompressibility conditions by $\tau$, the linear system of \eqref{Instationary_discrete_Oseen_problem_BE} can be rewritten as
\begin{equation}\label{Instationary_Oseen_system}
\underbrace{\left[
\begin{array}{cc}
\Phi^{(k)}_{\mathrm{BE}} & \left(\Psi_{\mathrm{BE}} \right)^\top\\
\Psi_{\mathrm{BE}} & - \Theta_{\mathrm{BE}}
\end{array}
\right]}_{\mathcal{A}^{(k)}_{\mathrm{BE}}}
\left[
\begin{array}{c}
\boldsymbol{\delta v}^{(k)} \\
\boldsymbol{\delta \zeta}^{(k)}\\
\boldsymbol{\delta \mu}^{(k)}\\
\boldsymbol{\delta p}^{(k)}
\end{array}
\right]=
\left[
\begin{array}{c}
\boldsymbol{b}^{(k)}_1 \\
\boldsymbol{b}^{(k)}_2 \\
\boldsymbol{b}^{(k)}_3 \\
\boldsymbol{b}^{(k)}_4
\end{array}
\right],
\end{equation}
where the right-hand side accounts for the non-linear residual. Further,
\begin{equation}\label{Phi_Psi_Theta_BE}
\Phi^{(k)}_{\mathrm{BE}}=
\left[
\begin{array}{cccc}
\mathcal{M}^\mathrm{BE} & \mathcal{L}^{\mathrm{BE},(k)}_1 \\
\mathcal{L}^{\mathrm{BE},(k)}_{2} & -\mathcal{M}^\mathrm{BE}_\beta
\end{array}
\right], \quad
\Psi_{\mathrm{BE}}=
\left[
\begin{array}{cccc}
\mathcal{B}^\mathrm{BE} & 0 \\
0 & \mathcal{B}^\mathrm{BE}
\end{array}
\right], \quad
\Theta_{\mathrm{BE}}
\left[
\begin{array}{cccc}
0 & 0 \\
0 & 0
\end{array}
\right],
\end{equation}
with $\mathcal{M}^\mathrm{BE}= I_{n_t+1,1} \otimes \bar{\mathbf{M}}^\mathrm{BE}$, $\bar{\mathcal{M}}^\mathrm{BE}_\beta = I_{n_t+1,2} \otimes \bar{\mathbf{M}}_\beta^\mathrm{BE} $, $\mathcal{B}^\mathrm{BE}= I_{n_t+1} \otimes \bar{B}$, and
\begin{align}
\nonumber & \mathcal{L}^{\mathrm{BE},(k)}_{1} = 
\left[
\begin{array}{cccccc}
\mathbf{T}^{(k)}_0 & -\mathbf{M} &  \\
 & \ddots & \ddots &  \\
  &  & \mathbf{T}^{(k)}_{n_t-1} & -\mathbf{M} \\
  &  &  & \mathbf{T}^{(k)}_{n_t}
\end{array}
\right],
 \;\; \mathcal{L}^{\mathrm{BE},(k)}_{2}=
\left[
\begin{array}{cccccc}
\mathbf{L}^{(k)}_0 & \\
-\mathbf{M} & \mathbf{L}^{(k)}_1 & \\
 & \ddots & \ddots & \\
  &  & -\mathbf{M} & \mathbf{L}^{(k)}_{n_t}
\end{array}
\right].
\end{align}
Note that, in the case of the incompressibility conditions not being solved exactly, $\mathcal{L}_1^{\mathrm{BE},(k)} \neq (\mathcal{L}_2^{\mathrm{BE},(k)})^\top$; however, the system is symmetric if they are solved exactly.

We note that we can relax the incompressibility assumptions on $\vec{v}_0$ and modify the discretization of the Oseen problem \eqref{Instationary_Oseen_problem}--\eqref{Instationary_Oseen_residuals} as follows. Suppose that $\vec{v}_0$ is not solenoidal. Then, for the first backward Euler step in \eqref{Instationary_discrete_Oseen_problem_BE} we can rewrite \eqref{solenoidal_projection_BE} as
\begin{equation*}
\left\{
\begin{array}{rl}
\mathbf{L}^{(k)}_0\: \boldsymbol{\delta v}^{(k)}_0 + \bar{B}^\top \boldsymbol{\delta p}^{(k)}_0 = & \hspace{-0.6em} \boldsymbol{R}^{(k)}_{1,-1},\\
B\: \boldsymbol{\delta v}^{(k)}_0 = & \hspace{-0.6em} r^{(k)}_{1,0},
\end{array}
\right.
\end{equation*}
where, given $\bar{\boldsymbol{v}}_0$ as an appropriate discretization of $\vec{v}_0$,
\begin{displaymath}
\left\{
\begin{array}{rl}
\vspace{0.25ex}
\boldsymbol{R}^{(k)}_{1,-1} = & \hspace{-0.6em} \mathbf{\overline{L}}_0\:  \bar{\boldsymbol{v}}_0 - \mathbf{L}_0^{(k)} \boldsymbol{v}^{(k)}_0 ,\\
r^{(k)}_{1,0} = & \hspace{-0.6em} - B \: \boldsymbol{v}^{(k)}_0.
\end{array}
\right.
\end{displaymath}
Here, $\mathbf{\overline{L}}_0=\tau (\nu \mathbf{K} + \mathbf{\overline{N}}_0 + \mathbf{\overline{W}}_0) + \mathbf{M}$, $\mathbf{\overline{N}}_0$ and $ \mathbf{\overline{W}}_0$ are the vector-convection and stabilization matrices related to $\bar{\boldsymbol{v}}_0$, and the rest of the non-linear residuals are defined as in \eqref{Instationary_discrete_Oseen_residual_BE}. Note that in this case we cannot substitute $\mathbf{M}\boldsymbol{v}^{(k)}_0 = \mathbf{M}\bar{\boldsymbol{v}}_0$ into the non-linear residuals as $\bar{\boldsymbol{v}}_0$ is not incompressible; we note also that for $k=0$ (meaning $\boldsymbol{v}^{(0)}_0=\boldsymbol{0}$) and with $\nu = 1$ the above step gives the solenoidal projection \eqref{solenoidal_projection} for the instationary Stokes control problem, with $\mathbf{\overline{L}}_0=\mathbf{L}^{(0)}_0 = \tau  \mathbf{K} + \mathbf{M}$.

\subsubsection{Crank--Nicolson for instationary Navier--Stokes control}\label{4_2_2}
In this section we present the linear system arising upon employing Crank--Nicolson in time when solving \eqref{Instationary_Oseen_problem}--\eqref{Instationary_Oseen_residuals}. Again discretizing the interval $(0,t_{f})$ into $n_t$ subintervals of length $\tau=\frac{t_{f}}{n_t}$, we approximate $\vec{v}$ and $\vec{\zeta}$ at the time points $t_n=n\tau$, $n=0,\ldots, n_t$, and use a staggered grid for $p$ and $\mu$, as in \cite{Bell_Colella_Glaz}. Specifically, our approximations of the solutions at the $k$-th non-linear iteration are given by $\boldsymbol{v}^{(k)}_n \approx \vec{v}(\mathbf{x},t_n)$, $\boldsymbol{\zeta}^{(k)}_n \approx \vec{\zeta}(\mathbf{x},t_n)$, for $n=0,\ldots, n_t$, and $\boldsymbol{p}^{(k)}_{n+\frac{1}{2}} \approx p(\mathbf{x},t_n+\frac{1}{2}\tau)$, $\boldsymbol{\mu}^{(k)}_{n+\frac{1}{2}} \approx \mu(\mathbf{x},t_n+\frac{1}{2}\tau)$, for $n=0,\ldots, n_t-1$, for all $\mathbf{x} \in \Omega$. Let us introduce the following finite element matrices:
\begin{equation*}
\left.
\begin{array}{rl}
\mathbf{L}^{\pm,(k)}_n = \frac{ \tau}{2} (\nu \mathbf{K} + \mathbf{N}^{(k)}_n + \mathbf{W}^{(k)}_n) \pm \mathbf{M}, & \mathbf{T}^{\pm,(k)}_n = \frac{ \tau}{2} (\nu \mathbf{K} - \mathbf{N}^{(k)}_n + \mathbf{W}^{(k)}_n) \pm \mathbf{M},\\
\bar{\mathbf{M}}^\mathrm{CN}= \frac{\tau}{2}\mathbf{M}, &  \bar{\mathbf{M}}^\mathrm{CN}_\beta=\frac{\tau}{2 \beta}\mathbf{M},
\end{array}
\right.
\end{equation*}
with $\mathbf{W}_n^{(k)}$, $\mathbf{N}_n^{(k)}$ defined as for backward Euler. Then the discrete Oseen iterate is
\begin{displaymath}
\begin{array}{lll}
\boldsymbol{v}_n^{(k+1)} = \boldsymbol{v}_n^{(k)} + \boldsymbol{\delta v}_n^{(k)}, & \quad \boldsymbol{\zeta}_n^{(k+1)} = \boldsymbol{\zeta}_n^{(k)} + \boldsymbol{\delta \zeta}_n^{(k)},& \qquad n = 0 , \ldots, n_t, \\
\boldsymbol{p}_{n+\frac{1}{2}}^{(k+1)} = \boldsymbol{p}_{n+\frac{1}{2}}^{(k)} + \boldsymbol{\delta p}_{n+\frac{1}{2}}^{(k)}, & \quad \boldsymbol{\mu}_{n+\frac{1}{2}}^{(k+1)} = \boldsymbol{\mu}_{n+\frac{1}{2}}^{(k)} + \boldsymbol{\delta \mu}_{n+\frac{1}{2}}^{(k)},&\qquad n = 0 , \ldots, n_t-1,
\end{array}
\end{displaymath}
with $\boldsymbol{\delta v}_n^{(k)}$, $ \boldsymbol{\delta \zeta}_n^{(k)}$, $ \boldsymbol{\delta p}_{n+\frac{1}{2}}^{(k)}$, $ \boldsymbol{\delta \mu}_{n+\frac{1}{2}}^{(k)}$ solutions of the following discretized version of \eqref{Instationary_Oseen_problem}:
\begin{equation*}
\!\! \left\{
\begin{array}{rll}
\vspace{0.25ex}
\displaystyle 
\!\!\! \bar{\mathbf{M}}^\mathrm{CN} ( \boldsymbol{\delta v}^{(k)}_n + \boldsymbol{\delta v}^{(k)}_{n+1} )+ \mathbf{T}^{+,(k)}_n \boldsymbol{\delta \zeta}^{(k)}_n + \mathbf{T}^{-,(k)}_{n+1} \boldsymbol{\delta \zeta}^{(k)}_{n+1} + \bar{B}^{\top} \boldsymbol{\delta \mu}^{(k)}_{n+\frac{1}{2}} = & \hspace{-0.6em}  \mathbf{R}^{(k)}_{2,n},\\
\vspace{0.25ex}
\displaystyle \!\!\! \mathbf{L}^{-,(k)}_n \boldsymbol{\delta v}^{(k)}_n + \mathbf{L}^{+,(k)}_{n+1} \boldsymbol{\delta v}^{(k)}_{n+1} + \bar{B}^{\top} \boldsymbol{\delta p}^{(k)}_{n+\frac{1}{2}} - \bar{\mathbf{M}}^\mathrm{CN}_\beta (\boldsymbol{\delta \zeta}^{(k)}_n + \boldsymbol{\delta \zeta}^{(k)}_{n+1} ) = & \hspace{-0.6em} \mathbf{R}^{(k)}_{1,n},\\
\vspace{0.25ex}
B \boldsymbol{\delta v}^{(k)}_{n+1} = & \hspace{-0.6em} \mathbf{r}^{(k)}_{1,n+1}, \\
B \boldsymbol{\delta \zeta}^{(k)}_{n} = & \hspace{-0.6em} \mathbf{r}^{(k)}_{2,n},
\end{array}
\right.
\end{equation*}
for $n=0,...,n_t-1$, with $\boldsymbol{\delta v}^{(k)}_0=\boldsymbol{0}$, $\boldsymbol{\delta \zeta}^{(k)}_{n_t}=\boldsymbol{0}$. The discretized residuals are given by
\begin{equation}\label{Instationary_discrete_Oseen_residual_CN}
\left\{
\begin{array}{rll}
\boldsymbol{R}^{(k)}_{1,n} = & \hspace{-0.6em} \frac{\tau}{2}( \boldsymbol{f}^{n}+ \boldsymbol{f}^{n+1}) - \mathbf{L}^{-,(k)}_n \boldsymbol{ v}^{(k)}_{n} - \mathbf{L}^{+,(k)}_{n+1} \boldsymbol{ v}^{(k)}_{n+1} - \bar{B}^{\top} \boldsymbol{ p}^{(k)}_{n+\frac{1}{2}}\\
\vspace{0.25ex}
&  + \bar{\mathbf{M}}^\mathrm{CN}_\beta (\boldsymbol{ \zeta}^{(k)}_{n} + \boldsymbol{ \zeta}^{(k)}_{n+1} ) ,\\
\vspace{0.25ex}
\boldsymbol{r}^{(k)}_{1,n+1} = & \hspace{-0.6em} - B \: \boldsymbol{ v}^{(k)}_{n+1}, \\
\boldsymbol{R}^{(k)}_{2,n} = & \hspace{-0.6em} \bar{\mathbf{M}}^\mathrm{CN}( \boldsymbol{v}_d^n +\boldsymbol{v}_d^{n+1} ) - \bar{\mathbf{M}}^\mathrm{CN}(\boldsymbol{ v}^{(k)}_n+ \boldsymbol{ v}^{(k)}_{n+1})  - \mathbf{T}^{+,(k)}_n \boldsymbol{ \zeta}^{(k)}_n\\
\vspace{0.25ex}
& - \mathbf{T}^{-,(k)}_{n+1} \boldsymbol{\zeta}^{(k)}_{n+1} - \bar{B}^{\top} \boldsymbol{ \mu}^{(k)}_{n+\frac{1}{2}} - \frac{\tau}{2}( \boldsymbol{\omega}^{(k)}_n+ \boldsymbol{\omega}^{(k)}_{n+1} ),\\
\boldsymbol{r}^{(k)}_{2,n} = & \hspace{-0.6em} - B \: \boldsymbol{ \zeta}^{(k)}_{n},
\end{array}
\right.
\end{equation}
for $n=0,...,n_t-1$, with $\mathbf{f}^n$ and $\boldsymbol{\omega}^{(k)}_n$ defined as for backward Euler, for $n=0,\ldots,n_t$. Note also that here the non-linear residuals $\boldsymbol{R}^{(k)}_{1,0}$, $ \boldsymbol{R}^{(k)}_{1,n_t-1}$, $ \boldsymbol{R}^{(k)}_{2,0}$, and $\boldsymbol{R}^{(k)}_{2,n_t-1}$ in \eqref{Instationary_discrete_Oseen_residual_CN} take into account the initial and final conditions on $\vec{v}$ and $\vec{\zeta}$.

In matrix form, after multipling the incompressibility constraints by $\tau$, we write
\begin{equation}\label{original_system_Oseen_CN}
\left.
\begin{array}{c}
\left[
\begin{array}{cccc}
\bar{\mathcal{M}}^\mathrm{CN} & \bar{\mathcal{L}}^{\mathrm{CN},(k)}_1 & (\bar{\mathcal{B}}_2^\mathrm{CN})^{\top} & 0 \\
\bar{\mathcal{L}}^{\mathrm{CN},(k)}_{2} & -\bar{\mathcal{M}}^\mathrm{CN}_\beta & 0 & (\bar{\mathcal{B}}_1^\mathrm{CN})^{\top}\\
\bar{\mathcal{B}}^\mathrm{CN}_1 & 0 & 0 & 0\\
0 & \bar{\mathcal{B}}_2^\mathrm{CN} & 0 & 0
\end{array}
\right]

\left[
\begin{array}{c}
\bar{\boldsymbol{\delta v}}^{(k)} \\
\bar{\boldsymbol{\delta \zeta}}^{(k)}\\
\bar{\boldsymbol{\delta \mu}}^{(k)}\\
\bar{\boldsymbol{\delta p}}^{(k)}
\end{array}
\right]
=
\left[
\begin{array}{c}
\bar{\boldsymbol{b}}_1^{(k)} \\
\bar{\boldsymbol{b}}_2^{(k)} \\
\bar{\boldsymbol{b}}_3^{(k)} \\
\bar{\boldsymbol{b}}_4^{(k)}
\end{array}
\right],
\end{array}
\right.
\end{equation}
where $\bar{\boldsymbol{\delta v}}^{(k)}$, $\bar{\boldsymbol{\delta \zeta}}^{(k)}$, $\bar{\boldsymbol{\delta \mu}}^{(k)}$, $\bar{\boldsymbol{\delta p}}^{(k)}$ are the $k$-th Oseen iterates, and the right-hand side accounts for the non-linear residual. The blocks in the previous matrix are given by
\begin{displaymath}
\begin{array}{rlrl}
\vspace{0.25ex}
\hspace{-0.4em} \bar{\mathcal{L}}^{\mathrm{CN},(k)}_{1} \! = & \hspace{-1em}
\left[ \!\!
\begin{array}{cccc}
\mathbf{T}^{+,(k)}_0 & \!\!\! \mathbf{T}^{-,(k)}_1 &  \\
 & \!\!\! \ddots & \!\!\!\! \ddots & \\
  &  & \!\!\!\! \mathbf{T}^{+,(k)}_{n_t-1} & \!\!\!\! \mathbf{T}^{-,(k)}_{n_t} \\
  &  &  & \!\!\!\! \mathbf{M}
\end{array}
\!\!\!
\right] \! ,
 &
\!\!\!
\bar{\mathcal{L}}^{\mathrm{CN},(k)}_{2} \! = & \hspace{-1em}
\left[ \!\!
\begin{array}{cccc}
\! \mathbf{M} & \\
 \mathbf{L}^{-,(k)}_0 & \!\!\! \mathbf{L}^{+,(k)}_1 & \\
 & \!\!\! \ddots & \!\!\!\! \ddots & \\
  &  & \!\!\!\! \mathbf{L}^{-,(k)}_{n_t-1} & \!\!\!\! \mathbf{L}^{+,(k)}_{n_t}
\end{array}
\!\!\!
\right] \! ,
 \\
\hspace{-0.4em}
\bar{\mathcal{B}}^\mathrm{CN}_1\!  = & \hspace{-1em}
\left[
\begin{array}{cccc}
0 & \bar{B} &  \\
 &  & \ddots & \\
 &  & & \bar{B} \\
\end{array}
\right] \! , &
\!\!
\bar{\mathcal{B}}^\mathrm{CN}_2 \! = & \hspace{-1em}
\left[
\begin{array}{cccc}
 \bar{B} &  \\
 & \ddots & \\
 & & \bar{B} & 0 \\
\end{array}
\right] \! ,
\end{array}
\end{displaymath}
and $\bar{\mathcal{M}}^\mathrm{CN} = (I_{n_t+1,1}+I_{n_t+1,3}) \otimes \bar{\mathbf{M}}^\mathrm{CN} $, $\bar{\mathcal{M}}^\mathrm{CN}_{\beta} = (I_{n_t+1,2}+I_{n_t+1,3}^\top) \otimes \bar{\mathbf{M}}^\mathrm{CN}_\beta$.

The system \eqref{original_system_Oseen_CN} is clearly not symmetric; however, we work as in \cite{Leveque_Pearson} in order to transform the linear system above and make it as close to symmetric as possible. In fact, eliminating the initial and final-time conditions on $\vec{v}$ and $\vec{\zeta}$, we can rewrite
\begin{displaymath}
\left.
\begin{array}{c}
\left[
\begin{array}{cccc}
\widetilde{\mathcal{M}}^\mathrm{CN} & \widetilde{\mathcal{L}}^{\mathrm{\: CN},(k)}_1 & (\widetilde{\mathcal{B}}^\mathrm{\: CN})^{\top} & 0 \\
\widetilde{\mathcal{L}}^{\mathrm{\: CN},(k)}_2 & -\widetilde{\mathcal{M}}^\mathrm{CN}_\beta & 0 & (\widetilde{\mathcal{B}}^\mathrm{\: CN})^{\top}\\
\widetilde{\mathcal{B}}^\mathrm{\: CN} & 0 & 0 & 0\\
0 & \widetilde{\mathcal{B}}^\mathrm{\: CN} & 0 & 0
\end{array}
\right]

\left[
\begin{array}{c}
\boldsymbol{\delta v}^{(k)} \\
\boldsymbol{\delta \zeta}^{(k)}\\
\boldsymbol{\delta \mu}^{(k)}\\
\boldsymbol{\delta p}^{(k)}
\end{array}
\right]
=
\left[
\begin{array}{c}
\boldsymbol{b}^{(k)}_1 \\
\boldsymbol{b}^{(k)}_2 \\
\boldsymbol{b}^{(k)}_3 \\
\boldsymbol{b}^{(k)}_4
\end{array}
\right],
\end{array}
\right.
\end{displaymath}
with $\boldsymbol{\delta v}^{(k)}$, $\boldsymbol{\delta \zeta}^{(k)}$, $\boldsymbol{\delta \mu}^{(k)}$, $\boldsymbol{\delta p}^{(k)}$ as well as the right-hand side modified accordingly. The matrices $\widetilde{\mathcal{M}}^\mathrm{CN} = I_{n_t,4}^\top \otimes \bar{\mathbf{M}}^\mathrm{CN}$, $\widetilde{\mathcal{M}}_{\beta}^\mathrm{CN} = I_{n_t,4} \otimes \bar{\mathbf{M}}^\mathrm{CN}_{\beta}$, $\widetilde{\mathcal{B}}^\mathrm{\: CN}= I_{n_t} \otimes \bar{B}$, and
\begin{align}
\nonumber & \widetilde{\mathcal{L}}^{\mathrm{\: CN},(k)}_1= 
\!\!
\left[
\begin{array}{cccccc}
\!\! \mathbf{T}^{+,(k)}_0 & \!\!\!\! \mathbf{T}^{-,(k)}_1\\
 &  \ddots &  \ddots & \\
 &   & \!\!\!\! \mathbf{T}^{+,(k)}_{n_t-2} & \!\!\!\! \mathbf{T}^{-,(k)}_{n_t-1}\\
   &  &  & \!\!\!\! \mathbf{T}^{+,(k)}_{n_t-1}
\end{array}
\!\!
\right] \!, 
\;
\widetilde{\mathcal{L}}^{\mathrm{\: CN},(k)}_2= 
\!\!
\left[
\begin{array}{cccccc}
\!\! \mathbf{L}^{+,(k)}_1 & \\
\!\! \mathbf{L}^{-,(k)}_1 & \!\!\!\! \mathbf{L}^{+,(k)}_2 & \\
 & \ddots & \ddots & \\
  &  & \!\!\!\! \mathbf{L}^{-,(k)}_{n_t-1} & \!\!\!\! \mathbf{L}^{+,(k)}_{n_t}
\end{array}
\!\!
\right] .
\end{align}
Using \texttt{blkdiag} to define a block diagonal matrix, we apply the linear transformation
\begin{equation*}
T=\texttt{blkdiag}(T_1,T_2,T_3,T_4),
\end{equation*}
where $T_1=I_{n_t,4} \otimes I_{n_v}, \, T_2=T_1^\top=I_{n_t,4}^\top \otimes I_{n_v}, \, T_3=I_{n_t,4}^\top \otimes I_{n_p}, \, T_4=T_3^\top=I_{n_t,4} \otimes I_{n_p}$, we may equivalently consider the following linear system:
\begin{equation}\label{modified_system_Oseen_CN}
\underbrace{\left[
\begin{array}{cc}\vspace{1ex}
\Phi^{(k)}_{\mathrm{CN}} & \left(\Psi_{\mathrm{CN}} \right)^\top\\
\Psi_{\mathrm{CN}} & - \Theta_{\mathrm{CN}}
\end{array}
\right]}_{\mathcal{A}^{(k)}_{\mathrm{CN}}}
\left[
\begin{array}{c}
\boldsymbol{\delta v}^{(k)}\\
\boldsymbol{\delta \zeta}^{(k)}\\
\boldsymbol{\delta \mu}^{(k)}\\
\boldsymbol{\delta p}^{(k)}
\end{array}
\right]
=
T
\left[
\begin{array}{c}
\boldsymbol{b}_1^{(k)} \\
\boldsymbol{b}_2^{(k)}\\
\boldsymbol{b}_3^{(k)}\\
\boldsymbol{b}_4^{(k)}
\end{array}
\right].
\end{equation}
Here the matrix blocks are given by
\begin{equation}\label{Phi_Psi_Theta_CN}
\Phi^{(k)}_{\mathrm{CN}}=\left[
\begin{array}{cc}
\mathcal{M}^{\mathrm{CN}} & \mathcal{L}^{\mathrm{CN},(k)}_1\\
\mathcal{L}^{\mathrm{CN},(k)}_2 & -\mathcal{M}^{\mathrm{CN}}_\beta
\end{array}\right], \;\:
\Psi_{\mathrm{CN}} = \left[
\begin{array}{cc}
\mathcal{B}^{\mathrm{CN}}_1 & 0\\
0 & \mathcal{B}^{\mathrm{CN}}_2
\end{array}
\right], \;\:
\Theta_{\mathrm{CN}}=\left[
\begin{array}{cc}
0 & 0\\
0 & 0
\end{array}
\right],
\end{equation}
with
\begin{align}
\nonumber & \mathcal{M}^{\mathrm{CN}}  =  T_1 \widetilde{\mathcal{M}}^{\mathrm{CN}}  =  \left( I_{n_t,4} \; I_{n_t,4}^\top \right)  \otimes  \bar{\mathbf{M}}, & \mathcal{B}^{\mathrm{CN}}_1  =  T_3 \widetilde{\mathcal{B}}^{\mathrm{\: CN}}  =  I_{n_t,4}^\top  \otimes  \bar{B}, \\
\nonumber & \mathcal{M}^{\mathrm{CN}}_\beta  =  T_2 \widetilde{\mathcal{M}}^{\mathrm{CN}}_{\beta} =  \left( I_{n_t,4}^\top \; I_{n_t,4} \right)  \otimes  \bar{\mathbf{M}}_\beta, & \mathcal{B}^{\mathrm{CN}}_2  =  T_4 \widetilde{\mathcal{B}}^{\mathrm{\: CN}}  =  I_{n_t,4}  \otimes  \bar{B},\\
& \label{L_1_L_2_CN}
\mathcal{L}^{\mathrm{CN},(k)}_1=T_1 \widetilde{\mathcal{L}}^{\mathrm{\: CN},(k)}_1, & \mathcal{L}^{\mathrm{CN},(k)}_2=T_2 \widetilde{\mathcal{L}}^{\mathrm{\: CN},(k)}_2.
\end{align}
System \eqref{modified_system_Oseen_CN} is still not symmetric in general, as $\mathcal{L}^{\mathrm{CN},(k)}_1 \neq (\mathcal{L}^{\mathrm{CN},(k)}_2)^\top$ due to the mismatch of the indices for the convection terms, however it is now symmetric when the above strategy is applied to the instationary Stokes control problem, due to the absence of the convection terms. \color{black}We observe that the transformations $T_i, \; i=1,\ldots,4$, as well as their inverse operations are easy and cheap to apply, as they require only a sequence of block updates. Therefore, in particular, we may rewrite 
\begin{equation}\label{M_CN_Mbeta_CN_transform}
\mathcal{M}^{\mathrm{CN}}=T_1 \mathcal{M}^{\mathrm{CN}}_D T_1^\top, \qquad
\mathcal{M}^{\mathrm{CN}}_\beta=T_2 \mathcal{M}^{\mathrm{CN}}_{D,\beta} T_2^\top,
\end{equation}
where
\begin{equation}\label{bar_M_CN_Mbeta_CN}
\mathcal{M}^{\mathrm{CN}}_D
=  I_{n_t} \otimes \bar{\mathbf{M}},
\qquad
\mathcal{M}^{\mathrm{CN}}_{D,\beta}
= I_{n_t} \otimes \bar{\mathbf{M}}_\beta.
\end{equation}
We may therefore work efficiently with $\mathcal{M}^{\mathrm{CN}}$ and $\mathcal{M}^{\mathrm{CN}}_\beta$, using $T_1$, $T_2$, $\mathcal{M}^{\mathrm{CN}}_D$, $\mathcal{M}^{\mathrm{CN}}_{D,\beta}$. Further, since both $\mathcal{M}^{\mathrm{CN}}_D$ and $\mathcal{M}^{\mathrm{CN}}_{D,\beta}$ are s.p.d., the same holds for $\mathcal{M}^{\mathrm{CN}}$ and $\mathcal{M}^{\mathrm{CN}}_\beta$\color{black}.

We point out that it is not straightforward to generalize the Crank--Nicolson discretization to the case where $\vec{v}_0$ is not incompressible. In fact, in this case we must also solve an appropriate solenoidal projection; however, the projection cannot be solved along with the other equations, as our approach requires the elimination of the initial and final conditions on $\vec{v}$ and $\vec{\zeta}$. Therefore, before applying our solver we must solve the projection to a stricter tolerance than that required for the control problem.

\section{Preconditioning Approach}
\label{5}
As the discretizations \eqref{Stationary_discrete_Oseen_system}, \eqref{Instationary_Oseen_system}, and \eqref{modified_system_Oseen_CN} of the optimality conditions for the problems under examination lead to matrices of the structure \eqref{MatrixA}, we now devise preconditioners for each system by making use of saddle-point theory. We employ a preconditioner of the form \eqref{optimal_prec}: this requires us to (approximately) apply the inverse of the corresponding $(1,1)$-block of each matrix analysed; we accelerate this process by again employing an approximation of the form \eqref{optimal_prec}. In the following, subscripts refer to the corresponding matrix we are considering; to simplify the notation, we drop the superscript referring to the non-linear iterate $k$.

\subsection{Approximation of the \boldmath{$(1,1)$}-block}
\label{5_1}
We now describe suitable approximations of the inverses of the $(1,1)$-blocks for the systems \eqref{Stationary_discrete_Oseen_system}, \eqref{Instationary_Oseen_system}, and \eqref{modified_system_Oseen_CN}. As noted after the discretization of the optimality conditions, each of these matrices is not symmetric if we solve the incompressibility constraints inexactly (for a Crank--Nicolson discretization the block is not symmetric even if those constraints are solved exactly). We thus use a fixed number of GMRES iterations to approximate the $(1,1)$-block, accelerated with the preconditioners described below, as opposed to Uzawa iteration for example (see \cite{Uzawa}) which may be symmetrized. However, we note that, for Stokes control problems, the $(1,1)$-block of each system is symmetric, allowing the use of a fixed number of Uzawa iterations coupled with the preconditioners below. 

\subsubsection{Stationary Navier--Stokes control}
\label{5_1_1}
Consider the $(1,1)$-block $\Phi^{(k)}_\mathrm{S}$ defined in \eqref{Phi_Psi_Theta_Stationary}. This matrix can be considered as the discretization of the optimality conditions for a stationary convection--diffusion control problem. Using saddle-point theory, a suitable preconditioner is given by
\begin{displaymath}
\mathcal{P}_{\Phi , \mathrm{S}}=
\left[
\begin{array}{cc}
\mathbf{M} & 0 \\
\mathbf{L}^{(k)} & -\mathbf{S}_{\Phi , \mathrm{S}}
\end{array}
\right],
\end{displaymath}
with $\mathbf{S}_{\Phi , \mathrm{S}} = \mathbf{M}_\beta + \mathbf{L}^{(k)} \mathbf{M}^{-1} \mathbf{L}^{(k)}_{\mathrm{adj}}$ the corresponding Schur complement. As described in \cite{Pearson_Wathen}, a potent preconditioner for $\mathcal{P}_{\Phi , \mathrm{S}}$ (optimal in the symmetric case) is given by
\begin{displaymath}
\widehat{\mathcal{P}}_{\Phi , \mathrm{S}}=
\left[
\begin{array}{cc}
\mathbf{M}_c & 0 \\
\mathbf{L}^{(k)} & -\widehat{\mathbf{S}}_{\Phi , \mathrm{S}}
\end{array}
\right].
\end{displaymath}
Here, $\mathbf{M}_c$ represents a fixed number of steps of the Chebyshev semi-iterative method \cite{GolubVargaI,GolubVargaII,Wathen_Rees}, and $\widehat{\mathbf{S}}_{\Phi , \mathrm{S}} = \big( \mathbf{L}^{(k)} + \mathbf{M}_{\sqrt{\beta}} \big) \mathbf{M}^{-1} \big( \mathbf{L}^{(k)}_{\mathrm{adj}} + \mathbf{M}_{\sqrt{\beta}} \big)$ with $\mathbf{M}_{\sqrt{\beta}} = \frac{1}{\sqrt{\beta}} \mathbf{M}$ and the blocks $\mathbf{L}^{(k)} + \mathbf{M}_{\sqrt{\beta}}$ and $\mathbf{L}^{(k)}_{\mathrm{adj}} + \mathbf{M}_{\sqrt{\beta}}$ approximated by the action of a multigrid routine, for example. It is worth noting that, if the incompressibility constraints are solved exactly, $\Phi^{(k)}_\mathrm{S}$ is symmetric and the approximation $\widehat{\mathbf{S}}_{\Phi , \mathrm{S}}$ of the Schur complement $\mathbf{S}_{\Phi , \mathrm{S}}$ is optimal; in fact, it can be proved that $\lambda(\widehat{\mathbf{S}}_{\Phi , \mathrm{S}}^{-1} \; \mathbf{S}_{\Phi , \mathrm{S}}) \in \left[\frac{1}{2}, 1 \right]$ \cite{Pearson_Wathen}.

\subsubsection{Instationary Navier--Stokes control with backward Euler}
\label{5_1_2}
We now derive a preconditioner for the matrix $\Phi^{(k)}_\mathrm{BE}$ defined in \eqref{Phi_Psi_Theta_BE}. As in the stationary case, the matrix can be considered as the discretization of the optimality conditions for an instationary convection--diffusion control problem with backward Euler in time. As the matrix $\mathcal{M}^\mathrm{BE}$ is not invertible, we seek a preconditioner of the form:
\begin{displaymath}
\widetilde{\mathcal{P}}_{\Phi , \mathrm{BE}}=
\left[
\begin{array}{cc}
\widetilde{\mathcal{M}}^\mathrm{BE} & 0 \\
\mathcal{L}_2^{\mathrm{BE},(k)} & -\widetilde{\mathcal{S}}_{\Phi , \mathrm{BE}}
\end{array}
\right],
\end{displaymath}
with $\widetilde{\mathcal{M}}^\mathrm{BE}$ an invertible approximation of $\mathcal{M}^\mathrm{BE}$, and the perturbed Schur complement $\widetilde{\mathcal{S}}_{\Phi , \mathrm{BE}} = \mathcal{M}^\mathrm{BE}_\beta + \mathcal{L}_2^{\mathrm{BE},(k)} \big(\widetilde{\mathcal{M}}^\mathrm{BE}\big)^{-1} \mathcal{L}^{\mathrm{BE},(k)}_1$. In \cite{Pearson_Stoll_Wathen}, the authors found for the heat control problem that a suitable approximation of $\mathcal{M}^\mathrm{BE}$ is given by $\widetilde{\mathcal{M}}^\mathrm{BE} =\texttt{blkdiag}(\bar{\mathbf{M}}^\mathrm{BE}, \ldots, \bar{\mathbf{M}}^\mathrm{BE}, \epsilon \bar{\mathbf{M}}^\mathrm{BE})$, with $\epsilon \ll 1$. Following \cite{Pearson_Stoll_Wathen}, we can derive that a good approximation for $\widetilde{\mathcal{S}}_{\Phi , \mathrm{BE}}$ is $\widehat{\mathcal{S}}_{\Phi , \mathrm{BE}}= \big( \mathcal{L}_2^{\mathrm{BE},(k)} + \mathcal{M}_{\sqrt{\beta}}^\mathrm{BE} \big) \big(\widetilde{\mathcal{M}}^\mathrm{BE}\big)^{-1} \big( \mathcal{L}^{\mathrm{BE},(k)}_1 + \mathcal{M}_{\sqrt{\beta}}^\mathrm{BE} \big)$, with $\mathcal{M}_{\sqrt{\beta}}^\mathrm{BE} = \frac{\tau}{\sqrt{\beta}} \texttt{blkdiag}(0, \mathbf{M} \ldots, \sqrt{\epsilon} \, \mathbf{M})$. As above, we do not apply the inverse of the blocks $\mathcal{L}_2^{\mathrm{BE},(k)} + \mathcal{M}_{\sqrt{\beta}}^\mathrm{BE}$ and $\mathcal{L}_1^{\mathrm{BE},(k)} + \mathcal{M}_{\sqrt{\beta}}^\mathrm{BE}$ exactly, but rather we apply a block-forward and block-backward substitution respectively, with each block on the diagonal approximated by the action of multigrid process, for instance. Thus, a suitable approximation of the matrix $\widetilde{\mathcal{P}}_{\Phi , \mathrm{BE}}$ is given by
\begin{displaymath}
\widehat{\mathcal{P}}_{\Phi , \mathrm{BE}}=
\left[
\begin{array}{cc}
\widehat{\mathcal{M}}^\mathrm{\: BE}_c & 0 \\
\mathcal{L}_2^{\mathrm{BE},(k)} & -\widehat{\mathcal{S}}_{\Phi , \mathrm{BE}}
\end{array}
\right],\qquad
\widehat{\mathcal{M}}^\mathrm{\: BE}_c = \tau \, \texttt{blkdiag}(\mathbf{M}_c, \ldots, \mathbf{M}_c, \epsilon \mathbf{M}_c).
\end{displaymath}


\subsubsection{Instationary Navier--Stokes control with Crank--Nicolson}
\label{5_1_3}
We focus now on devising a preconditioner for the linear system $\Phi^{(k)}_\mathrm{CN}$ defined in \eqref{Phi_Psi_Theta_CN}, arising from a Crank--Nicolson discretization. Similarly to the backward Euler case, this matrix can be considered as the discretization of the optimality conditions for the control of the instationary convection--diffusion equation discretized using Crank--Nicolson in time. Again, we seek to use the block triangular matrix
\begin{displaymath}
\mathcal{P}_{\Phi , \mathrm{CN}}=
\left[
\begin{array}{cc}
\mathcal{M}^{\mathrm{CN}} & 0 \\
\mathcal{L}^{\mathrm{CN},(k)}_2 & -S_{\Phi , \mathrm{CN}}
\end{array}
\right]
\end{displaymath}
as a preconditioner, where $S_{\Phi , \mathrm{CN}}= \mathcal{M}^{\mathrm{CN}}_\beta +\mathcal{L}^{\mathrm{CN},(k)}_2 (\mathcal{M}^{\mathrm{CN}})^{-1} \mathcal{L}^{\mathrm{CN},(k)}_1$. In order to find an approximation of $\mathcal{P}_{\Phi , \mathrm{CN}}$, we adapt the strategy used in \cite{Leveque_Pearson} as follows.

From \eqref{M_CN_Mbeta_CN_transform}--\eqref{bar_M_CN_Mbeta_CN}, $\mathcal{M}^{\mathrm{CN}}$ can be written as $\mathcal{M}^{\mathrm{CN}}=T_1 \mathcal{M}^{\mathrm{CN}}_D T_1^\top$, with $\mathcal{M}^{\mathrm{CN}}_D$ a block diagonal matrix with each diagonal block a multiple of $\mathbf{M}$. Therefore, a good approximation of $\mathcal{M}^{\mathrm{CN}}$ is given by $\widehat{\mathcal{M}}^{\mathrm{CN}}=T_1 \widehat{\mathcal{M}}^{\mathrm{CN}}_D T_1^\top$, with $\widehat{\mathcal{M}}^{\mathrm{CN}}_D = \frac{\tau}{2} I_{n_t} \otimes \mathbf{M}_c$.

To derive an approximation of $S_{\Phi , \mathrm{CN}}$, we use \eqref{M_CN_Mbeta_CN_transform} together with \eqref{L_1_L_2_CN} to rewrite
\begin{equation}\label{S_Phi_k_CN}
S_{\Phi , \mathrm{CN}} = T_2 \underbrace{ \big[ \mathcal{M}^{\mathrm{CN}}_{D,\beta} + \big( \widetilde{\mathcal{L}}^{\, \mathrm{CN},(k)}_2 \big) \big(\mathcal{M}^{\mathrm{CN}} \big)^{-1} \big( T_1 \widetilde{\mathcal{L}}^{\, \mathrm{CN},(k)}_1 \: T_1^{-1} \big) \big] }_{S_{\Phi , \mathrm{CN}}^{\mathrm{int}}}  T_1,
\end{equation}
recalling that $T_1=T_2^\top$. We first seek an approximation $\widehat{S}_{\Phi , \mathrm{CN}}^{\: \mathrm{int}}$ for $S_{\Phi , \mathrm{CN}}^{\mathrm{int}}$ of the form
\begin{displaymath}
\widehat{S}_{\Phi , \mathrm{CN}}^{\: \mathrm{int}} = \big( \widetilde{\mathcal{L}}^{\, \mathrm{CN},(k)}_2 + \widehat{\mathcal{M}}_2 \big) \big(\mathcal{M}^{\mathrm{CN}} \big)^{-1} \big( T_1 \widetilde{\mathcal{L}}^{\, \mathrm{CN},(k)}_1 \: T_1^{-1} + \widehat{\mathcal{M}}_1 \big),
\end{displaymath}
such that
\begin{displaymath}
\widehat{\mathcal{M}}_2 \big(\mathcal{M}^{\mathrm{CN}} \big)^{-1} \widehat{\mathcal{M}}_1 = \big( \widehat{\mathcal{M}}_2 T_2^{-1} \big) \big(\mathcal{M}^{\mathrm{CN}}_D \big)^{-1} \big( T_1^{-1} \widehat{\mathcal{M}}_1 \big) = \mathcal{M}^{\mathrm{CN}}_{D,\beta}.
\end{displaymath}
The previous expression is clearly satisfied with the choice $\widehat{\mathcal{M}}_2 T_2^{-1} = T_1^{-1} \widehat{\mathcal{M}}_1 = \frac{\tau}{2 \sqrt{\beta}} I_{n_t} \otimes \mathbf{M}$. Then, our approximation of $S_{\Phi , \mathrm{CN}}^{\: \mathrm{int}}$ is given by
\begin{align*}
\widehat{S}_{\Phi , \mathrm{CN}}^{\: \mathrm{int}} & = \big( \widetilde{\mathcal{L}}^{\, \mathrm{CN},(k)}_2 + \widehat{\mathcal{M}} \big) \big(\mathcal{M}^{\mathrm{CN}} \big)^{-1} \big( T_1 \widetilde{\mathcal{L}}^{\, \mathrm{CN},(k)}_1 \: T_1^{-1} + \widehat{\mathcal{M}}^{\: \top} \big) \\
 & = \big( \widetilde{\mathcal{L}}^{\, \mathrm{CN},(k)}_2 + \widehat{\mathcal{M}} \big) T_2^{-1} \big(\mathcal{M}^{\mathrm{CN}}_D \big)^{-1} \big( \widetilde{\mathcal{L}}^{\, \mathrm{CN},(k)}_1 \: T_1^{-1} + T_1^{-1} \widehat{\mathcal{M}}^{\: \top} \big),
\end{align*}
with $\widehat{\mathcal{M}} = \frac{\tau}{2 \sqrt{\beta}} I_{n_t,4}^\top \otimes \mathbf{M}$. Finally, substituting $\widehat{S}_{\Phi , \mathrm{CN}}^{\: \mathrm{int}}$ into \eqref{S_Phi_k_CN} and observing that $\widehat{\mathcal{M}}$ and $T_1$ commute, we obtain that our approximation of $S_{\Phi , \mathrm{CN}}$ is given by
\begin{displaymath}
\widehat{S}_{\Phi , \mathrm{CN}} = T_2 \big( \widetilde{\mathcal{L}}^{\, \mathrm{CN},(k)}_2 + \widehat{\mathcal{M}} \big) T_2^{-1} \big(\mathcal{M}^{\mathrm{CN}}_D \big)^{-1} \big( \widetilde{\mathcal{L}}^{\, \mathrm{CN},(k)}_1  +  \widehat{\mathcal{M}}^{\: \top} \big).
\end{displaymath}
As for backward Euler, we approximate the blocks $\widetilde{\mathcal{L}}^{\, \mathrm{CN},(k)}_2 + \widehat{\mathcal{M}}$ and $\widetilde{\mathcal{L}}^{\, \mathrm{CN},(k)}_1  +  \widehat{\mathcal{M}}^{\: \top}$ using a block-forward and block-backward substitution, with the action of a multigrid process used to apply the inverse of each block diagonal entry inexactly.

It is not possible in general to prove bounds on eigenvalues for this preconditioner, derived for instationary Navier--Stokes control with Crank--Nicolson in time. However, for instationary Stokes control, the preconditioner derived here reduces to that derived in \cite{Leveque_Pearson} for the heat control problem, which was proved to be optimal, and such that the spectrum of the preconditioned Schur complement is contained in $[\frac{1}{2},1]$.

\subsection{Approximation of Schur complement}
\label{5_2}
We now derive efficient approximations for each Schur complement of the systems \eqref{Stationary_discrete_Oseen_system}, \eqref{Instationary_Oseen_system}, and \eqref{modified_system_Oseen_CN}. Since the $(1,2)$- and the $(2,1)$-block of these systems can be considered as a \emph{vector-divergence matrix}, we make use of the commutator argument presented in Section \ref{3}.

\subsubsection{Stationary Navier--Stokes control}
\label{5_2_1}
Let us consider the Schur complement $S_{\mathcal{A},\mathrm{S}}= \Psi_\mathrm{S}(\Phi^{(k)}_\mathrm{S})^{-1} \Psi_\mathrm{S}^\top$ of the system \eqref{Stationary_discrete_Oseen_system}, with $\Phi^{(k)}_\mathrm{S}$ and $\Psi_\mathrm{S}$ defined as in \eqref{Phi_Psi_Theta_Stationary}. We apply the commutator argument to $\mathcal{E}_n$ as defined in \eqref{general_commutator} with $n=2$, with the differential operator on the velocity space defined as
\begin{displaymath}
\mathcal{D} = \left[
\begin{array}{cc}
\text{Id} &  -\nu \nabla^2 - \vec{v}^{\: (k)} \cdot \nabla\\
-\nu \nabla^2 + \vec{v}^{\: (k)} \cdot \nabla & -\frac{1}{\beta} \text{Id}
\end{array}
\right],
\end{displaymath}
and $\mathcal{D}_p$ the corresponding differential operator on the pressure space; we recall from Section \ref{2_1} that $\text{Id}$ represents the identity operator. Employing stable finite elements and working as in Section \ref{3}, we obtain the following expression for \eqref{commutator_approx}:
\begin{displaymath}
\widehat{S}_{\mathcal{A},\mathrm{S}} = \left[
\begin{array}{cc}
K_p & 0 \\
0 & K_p
\end{array}
\right]
\left[
\begin{array}{cc}
M_p & L^{(k)}_{\mathrm{adj},p} \\
L^{(k)}_p & -M_{\beta,p}
\end{array}
\right]^{-1}
\left[
\begin{array}{cc}
M_p & 0 \\
0 & M_p
\end{array}
\right] \approx S_{\mathcal{A},\mathrm{S}},
\end{displaymath}
where we set $L^{(k)}_p= \nu K_p + N^{(k)}_p + W^{(k)}_p$, $L^{(k)}_{\mathrm{adj},p}= \nu K_p - N^{(k)}_p +W^{(k)}_p$, and $M_{\beta,p}= \frac{1}{\beta}M_p$.

\subsubsection{Instationary Navier--Stokes control with backward Euler}
\label{5_2_2}
We now derive an approximation to the Schur complement $S_{\mathcal{A}, \mathrm{BE}}= \Psi_\mathrm{BE}(\Phi^{(k)}_\mathrm{BE})^{-1} \Psi_\mathrm{BE}^\top$ of \eqref{Instationary_Oseen_system}. As above, we apply the commutator argument \eqref{general_commutator}; however, we do not consider the instationary Navier--Stokes equation as part of the differential operator $\mathcal{D}$, but rather employ an operator that ``mimics'' the blocks of $\Phi^{(k)}_\mathrm{BE}$ defined in \eqref{Phi_Psi_Theta_BE}. With this aim, we consider \eqref{general_commutator} with $n=2(n_t+1)$ and the differential operator:
\begin{displaymath}
\mathcal{D} = \left[
\begin{array}{cc}
\vspace{1ex}
\mathcal{D}^{1,1}_{\mathrm{BE}} &  \mathcal{D}^{1,2}_{\mathrm{BE}}\\
\mathcal{D}^{2,1}_{\mathrm{BE}} & \mathcal{D}^{2,2}_{\mathrm{BE}}
\end{array}
\right],
\end{displaymath}
where $\mathcal{D}^{1,1}_{\mathrm{BE}} = \tau I_{n_t+1,1} \otimes \text{Id}$, $\mathcal{D}^{2,2}_{\mathrm{BE}} = - \frac{\tau}{\beta} I_{n_t+1,2} \otimes \text{Id}$, and
\begin{displaymath}
\begin{array}{cc}
\mathcal{D}^{1,2}_{\mathrm{BE}} = \left[ \!
\begin{array}{ccccc}
\mathcal{D}_{0,\: \mathrm{adj}} & -\text{Id} \\
 & \ddots & \ddots\\
 & & \mathcal{D}_{n_t-1,\: \mathrm{adj}} & -\text{Id}\\
 & & & \mathcal{D}_{n_t,\: \mathrm{adj}}
\end{array} \!
\right] \! ,
\;\;
\mathcal{D}^{2,1}_{\mathrm{BE}} = \left[ \!
\begin{array}{ccccc}
\mathcal{D}_0 & \\
-\text{Id} & \mathcal{D}_{1} \\
 & \ddots & \ddots\\
 & & -\text{Id} & \mathcal{D}_{n_t}
\end{array} \!
\right] \! ,
\end{array}
\end{displaymath}
with $\mathcal{D}_i = \tau(-\nu \nabla^2 + \vec{v}_i^{\: (k)} \cdot \nabla) + \text{Id}$ and $\mathcal{D}_{i,\: \mathrm{adj}} = \tau(-\nu \nabla^2 - \vec{v}_i^{\: (k)} \cdot \nabla) + \text{Id}$. As above, we define $\mathcal{D}_p$ as the corresponding differential operator on the pressure space. Discretizing \eqref{general_commutator} and observing that $S_{\mathcal{A}, \mathrm{BE}}=\tau^2 \vec{B} \: \mathbf{D}^{-1} \vec{B}^{\: \top}$, with $\mathbf{D}$ the discretization of the differential operator $\mathcal{D}$ and $\vec{B} = I_{2(n_t+1)} \otimes B$, we obtain the following approximation:
\begin{displaymath}
\widehat{S}_{\mathcal{A},\mathrm{BE}} = \tau^2 \mathcal{K}_p^{\mathrm{BE}} \left[
\begin{array}{cc}
\vspace{1ex}
D^{1,1}_{p, \: \mathrm{BE}} & D^{1,2}_{p, \: \mathrm{BE}}\\
D^{2,1}_{p, \: \mathrm{BE}} & D^{2,2}_{p, \: \mathrm{BE}}
\end{array}
\right]^{-1}
\mathcal{M}_p^{\mathrm{BE}} \approx S_{\mathcal{A},\mathrm{BE}}.
\end{displaymath}
Here, we set $\mathcal{K}_p^{\mathrm{BE}} = I_{2(n_t+1)} \otimes K_p$, $\mathcal{M}_p^{\mathrm{BE}} = I_{2(n_t+1)} \otimes M_p$, $\mathcal{D}^{1,1}_{p, \:\mathrm{BE}} = \tau I_{n_t+1,1} \otimes M_p$, $\mathcal{D}^{2,2}_{p, \: \mathrm{BE}} = - \frac{\tau}{\beta} I_{n_t+1,2} \otimes M_p$, and
\begin{displaymath}
\begin{array}{cc}
D^{1,2}_{p, \: \mathrm{BE}} = \left[
\begin{array}{ccccc}
\!\! T^{(k)}_{0,p} & \!\! -M_p \\
 & \ddots & \ddots\\
 & & \!\! T^{(k)}_{n_t-1,p} & \!\! -M_p\\
 & & & \!\! T^{(k)}_{n_t,p}
\end{array}
\right] \! ,
\;\;
D^{2,1}_{p, \: \mathrm{BE}} = \left[
\begin{array}{ccccc}
\!\! L^{(k)}_{0,p} & \\
\!\! -M_p & \!\! L^{(k)}_{1,p} \\
 & \ddots & \ddots\\
 & & \!\! -M_p & \!\! L^{(k)}_{n_t,p}
\end{array}
\right],
\end{array}
\end{displaymath}
with $L^{(k)}_{i,p}= \tau ( \nu K_p + N_{i,p}^{(k)} + W_{i,p}^{(k)}) + M_p$ and $T^{(k)}_{i,p}= \tau ( \nu K_p - N_{i,p}^{(k)} + W_{i,p}^{(k)}) + M_p $.

\subsubsection{Instationary Navier--Stokes control with Crank--Nicolson}
\label{5_2_3}
As for the Schur complement arising from the backward Euler discretization, we apply the commutator argument \eqref{general_commutator}, employing a differential operator $\mathcal{D}$ that mimics the blocks of a suitable matrix. Before presenting $\mathcal{D}$, we note that the Schur complement $S_{\mathcal{A}, {\mathrm{CN}}}=\Psi_\mathrm{CN}(\Phi^{(k)}_\mathrm{CN})^{-1} \Psi_\mathrm{CN}^\top$ can be rewritten as
\begin{displaymath}
S_{\mathcal{A},{\mathrm{CN}}}=\left[
\begin{array}{cc}
T_3 & 0\\
0 & T_4
\end{array}
\right]
\left[
\begin{array}{cc}
\widetilde{\mathcal{B}}^\mathrm{\: CN} & 0\\
0 & \widetilde{\mathcal{B}}^\mathrm{\: CN}
\end{array}
\right]
\left[
\begin{array}{cc}
\widetilde{\mathcal{M}}^\mathrm{CN} & \widetilde{\mathcal{L}}^{\mathrm{\: CN},(k)}_1\\
\widetilde{\mathcal{L}}^{\mathrm{\: CN},(k)}_2 & -\widetilde{\mathcal{M}}^\mathrm{CN}_\beta
\end{array}
\right]^{-1}
\left[
\begin{array}{cc}
\widetilde{\mathcal{B}}^\mathrm{\: CN} & 0\\
0 & \widetilde{\mathcal{B}}^\mathrm{\: CN}
\end{array}
\right]^\top.
\end{displaymath}

We now consider \eqref{general_commutator} with $n=2 n_t$ and the differential operator
\begin{displaymath}
\mathcal{D} = \left[
\begin{array}{cc}
\vspace{1ex}
\mathcal{D}^{1,1}_{\mathrm{CN}} &  \mathcal{D}^{1,2}_{\mathrm{CN}}\\
\mathcal{D}^{2,1}_{\mathrm{CN}} & \mathcal{D}^{2,2}_{\mathrm{CN}}
\end{array}
\right],
\end{displaymath}
where $\mathcal{D}^{1,1}_{\mathrm{CN}} = \frac{\tau}{2} I_{n_t,4}^\top \otimes \text{Id} $, $\mathcal{D}^{2,2}_{\mathrm{CN}} = - \frac{\tau}{2\beta} I_{n_t,4} \otimes \text{Id}$, and
\begin{displaymath}
\begin{array}{cc}
\mathcal{D}^{1,2}_{\mathrm{CN}} = \!\! \left[
\begin{array}{ccccc}
\!\! \mathcal{D}^+_{0,\: \mathrm{adj}} & \!\!\! \mathcal{D}^-_{1,\: \mathrm{adj}} \\
 & \ddots & \ddots\\
 & & \!\!\! \mathcal{D}^+_{n_t-2,\: \mathrm{adj}} & \!\!\! \mathcal{D}^-_{n_t-1,\: \mathrm{adj}}\\
 & & & \!\! \mathcal{D}^+_{n_t-1,\: \mathrm{adj}}
\end{array}
\right] \! ,
\;\;
\mathcal{D}^{2,1}_{\mathrm{CN}} = \left[
\begin{array}{ccccc}
\!\! \mathcal{D}^+_1 & \\
\!\! \mathcal{D}^-_1 & \!\!\! \mathcal{D}^+_2 \\
 & \ddots & \ddots\\
 & & \!\!\! \mathcal{D}^-_{n_t-1} & \!\!\! \mathcal{D}^+_{n_t}
\end{array}
\right] \! ,
\end{array}
\end{displaymath}
with $\mathcal{D}^\pm_i = \frac{\tau}{2}(-\nu \nabla^2 + \vec{v}_i^{\: (k)} \cdot \nabla) \pm \text{Id}$ and
$\mathcal{D}^\pm_{i,\: \mathrm{adj}} = \frac{\tau}{2}(-\nu \nabla^2 - \vec{v}_i^{\: (k)} \cdot \nabla)\pm \text{Id}$.
Again, we define $\mathcal{D}_p$ as the corresponding differential operator on the pressure space. Proceeding as above, we then derive the following approximation:
\begin{displaymath}
\widehat{S}_{\mathcal{A},{\mathrm{CN}}} = \tau^2 \left[
\begin{array}{cc}
T_3 & 0\\
0 & T_4
\end{array}
\right] \mathcal{K}_p^{\mathrm{CN}}
\left[
\begin{array}{cc}
\vspace{1ex}
D^{1,1}_{p, \: \mathrm{CN}} & D^{1,2}_{p, \: \mathrm{CN}}\\
D^{2,1}_{p, \: \mathrm{CN}} & D^{2,2}_{p, \: \mathrm{CN}}
\end{array}
\right]^{-1}
\mathcal{M}_p^{\mathrm{CN}} \approx S_{\mathcal{A},{\mathrm{CN}}}.
\end{displaymath}
Here, we set $\mathcal{K}_p^{\mathrm{CN}} = I_{2n_t} \otimes K_p$, $\mathcal{M}_p^{\mathrm{CN}} = I_{2n_t} \otimes M_p$, $\mathcal{D}^{1,1}_{p, \: \mathrm{CN}} = \frac{\tau}{2} I_{n_t,4}^\top \otimes M_p $, $\mathcal{D}^{2,2}_{p,\: \mathrm{CN}} = - \frac{\tau}{2\beta} I_{n_t,4} \otimes M_p$, and
\begin{displaymath}
\begin{array}{cc}
D^{1,2}_{p, \: \mathrm{CN}} \! = \!\! \left[
\begin{array}{ccccc}
\!\!\! T^{+,(k)}_{0,p} & \!\!\!\! T^{-,(k)}_{1,p} \\
 & \!\!\!\! \ddots & \!\!\!\! \ddots\\
 & & \!\!\!\! T^{+,(k)}_{n_t-2,p} & \!\!\!\! T^{-,(k)}_{n_t-1,p}\\
 & & & \!\!\!\! T^{+,(k)}_{n_t-1,p}
\end{array}
\!\!
\right] \! ,
\;
D^{2,1}_{p, \: \mathrm{CN}} \! = \!\! \left[
\begin{array}{ccccc}
\!\!\! L^{+,(k)}_{1,p} & \\
\!\!\! L^{-,(k)}_{1,p} & \!\!\!\!\! L^{+,(k)}_{2,p} \\
 & \!\!\!\! \ddots & \!\!\!\! \ddots\\
 & & \!\!\!\! L^{-,(k)}_{n_t-1,p} & \!\!\!\! L^{+,(k)}_{n_t,p}
\end{array}
\!\!\!
\right] \! ,
\end{array}
\end{displaymath}
with
\begin{displaymath}
\begin{array}{ll}
L^{\pm,(k)}_{i,p} = \frac{\tau}{2}(\nu K_p + N^{(k)}_{i,p}+W^{(k)}_{i,p}) \pm M_p, & \quad T^{\pm,(k)}_{i,p} = \frac{\tau}{2}(\nu K_p - N^{(k)}_{i,p}+W^{(k)}_{i,p})\pm M_p.
\end{array}
\end{displaymath}

\begin{remark}
To summarize, aside from matrix--vector products, the main computational work for our Crank--Nicolson preconditioner involves $n_t$ applications of Chebyshev semi-iteration to $\mathbf{M}$ and $2n_t$ multigrid processes per Uzawa iteration, in addition to $2n_t$ applications of Chebyshev semi-iteration to $M_p$ and $2n_t$ multigrid processes for $K_p$ to approximate the Schur complement. This is a similar computational workload as for the backward Euler preconditioner, as the latter requires $n_t+1$ application of Chebyshev semi-iteration and $2(n_t+1)$ applications of a multigrid process per Uzawa iteration, and $2(n_t+1)$ approximations of $M_p$ and $K_p$ for the Schur complement approximation.
\end{remark}

\section{Numerical Results}
\label{6}
We now demonstrate the effectiveness of our preconditioners by presenting numerical results. In all our tests, $d=2$ (that is, $\mathbf{x}=(x_1,x_2)$), and $\Omega=(-1,1)^2$. All tests are run on MATLAB R2018b, using a 1.70GHz Intel quad-core i5 processor and 8 GB RAM on an Ubuntu 18.04.1 LTS operating system.

As our preconditioners are non-symmetric and require an inner solve for the $(1,1)$-block, for the outer solver we apply flexible GMRES \cite{Saad} restarted every 10 iterations, up to a tolerance $10^{-6}$ on the relative residual (unless otherwise stated). Our implementation is based on the flexible GMRES routine in the TT-Toolbox \cite{TT_Toolbox}. To apply the approximate inverse of the $(1,1)$-block, we take 5 iterations of the GMRES routine implemented in MATLAB. We apply 20 steps of Chebyshev semi-iteration to mass matrices (on the velocity or pressure space); we apply 4 V-cycles of the AGMG routine \cite{Napov_Notay,Notay1,Notay2,AGMG} for other matrices constructed on the velocity space, while employing 2 V-cycles (with 2 symmetric Gauss--Seidel iterations for pre-/post-smoothing) of the \texttt{HSL\_MI20} solver \cite{HSL_MI20} for stiffness matrices on the pressure space within our Schur complement approximation.

Regarding the non-linear iteration for solving the Navier--Stokes control problem, we allow 20 Oseen iterations, specifying as a stopping criteria a reduction of $10^{-5}$ on the (non-linear) relative residual; the initial residual is the right-hand side of the corresponding Stokes control problem, with $\nu = 1$ (for the instationary case with Crank--Nicolson, we evaluate the residual before applying $T$). For each problem below, the first Oseen iterate is employed for the Stokes control solve, whose solutions  $\boldsymbol{v}^{(1)}, \: \boldsymbol{\zeta}^{(1)}, \: \boldsymbol{p}^{(1)}, \: \boldsymbol{\mu}^{(1)}$ are then used as the initial guess. We use inf--sup stable Taylor--Hood $\mathbf{Q}_2$--$\mathbf{Q}_1$ finite elements in the spatial dimensions, with level of refinement $l$ representing a (spatial) uniform grid of mesh-size $h=2^{1-l}$ for $\mathbf{Q}_1$ basis functions, and $h=2^{-l}$ for $\mathbf{Q}_2$ elements, in each dimension. All CPU times below are reported in seconds.

\subsection{Stationary Navier--Stokes control}\label{6_1}
We first test our solver on the stationary Navier--Stokes control problem \eqref{Stationary_Navier_Stokes_control_functional}--\eqref{Stationary_Navier_Stokes_control_constraints}. We set $\vec{f}=\vec{0}$, $\vec{v}_d=\vec{0}$, and
\begin{displaymath}
\vec{g}_D = \left\{
\begin{array}{ll}
\left[1,0\right]^\top & \mathrm{on} \: \partial \Omega_1 := \left(-1,1 \right)\times \left\{1\right\},\\
\left[0,0\right]^\top & \mathrm{on} \:\partial \Omega \setminus \partial\Omega_1.
\end{array}
\right.
\end{displaymath}
We report the average number of GMRES iterations together with the average CPU time per GMRES solve in Tables \ref{Stationary_NS_test1_table1}--\ref{Stationary_NS_test1_table3}, and in Table \ref{Number_Oseen_Stationary_NS} we state the total degrees of freedom (DoF) together with the total number of Oseen iterations required. We provide results for different levels of refinement $l$, values of $\beta$, and viscosities $\nu$.

\begin{table}[!ht]
\caption{Average GMRES iterations and CPU times for stationary Navier--Stokes control problem, for $\nu=\frac{1}{20}$ and a range of $l$, $\beta$. In brackets are the results for the corresponding Stokes control problem.}\label{Stationary_NS_test1_table1}
\begin{small}
\begin{center}
\begin{tabular}{|c||c|c|c|c|c|c|c|c|c|c|c|c|c|c|}
\hline
\multicolumn{1}{|c||}{\!\!\!\!\!\!\!\!\!\!\!\!\!\!\!\! \phantom{$1$} \!\!\!\!\!\!\!\!\!\!\!\!\!\!} & \multicolumn{2}{c|}{$\beta=10^0$}& \multicolumn{2}{c|}{$\beta=10^{-1}$}& \multicolumn{2}{c|}{$\beta=10^{-2}$}& \multicolumn{2}{c|}{$\beta=10^{-3}$}& \multicolumn{2}{c|}{$\beta=10^{-4}$}& \multicolumn{2}{c|}{$\beta=10^{-5}$}& \multicolumn{2}{c|}{$\beta=10^{-6}$}\\
\cline{2-15}
\!\!\!\!\!\!\!\!\!\!\!\! $l$ \!\!\!\!\!\!\!\!\!\!\!\! & \!\!\!\!\! $\texttt{it}$ \!\!\!\!\! & \!\!\!\! CPU \!\!\!\! & \!\!\!\! $\texttt{it}$ \!\!\!\! & \!\!\!\! CPU \!\!\!\! & \!\!\!\! $\texttt{it}$ \!\!\!\! & \!\!\!\! CPU \!\!\!\! & \!\!\!\! $\texttt{it}$ \!\!\!\! & \!\!\!\! CPU \!\!\!\! & \!\!\!\! $\texttt{it}$ \!\!\!\! & \!\!\!\! CPU \!\!\!\! & \!\!\!\! $\texttt{it}$ \!\!\!\! & \!\!\!\! CPU \!\!\!\! & \!\!\!\! $\texttt{it}$ \!\!\!\! & \!\!\!\!\! CPU \!\!\!\!\! \\
\hline
\hline
\multirow{2}{*}{\!\!\!\! $3$ \!\!\!\!\!\!\!\!\!\!\!\! } & \!\!\!\!\! 21 \!\!\!\!\! & \!\!\!\! 0.40 \!\!\!\! & \!\!\!\! 19 \!\!\!\! & \!\!\!\! 0.37 \!\!\!\! & \!\!\!\! 15 \!\!\!\! & \!\!\!\! 0.24 \!\!\!\! & \!\!\!\! 12 \!\!\!\! & \!\!\!\! 0.12 \!\!\!\! & \!\!\!\! 11 \!\!\!\! & \!\!\!\! 0.17 \!\!\!\! & \!\!\!\! 10 \!\!\!\! & \!\!\!\! 0.13 \!\!\!\! & \!\!\!\! 9 \!\!\!\! & \!\!\!\!\! 0.13 \!\!\!\!\! \\
& \!\!\!\!\! (15) \!\!\!\!\! & \!\!\!\! (0.34) \!\!\!\! & \!\!\!\! (18) \!\!\!\! & \!\!\!\! (0.41) \!\!\!\! & \!\!\!\! (17) \!\!\!\! & \!\!\!\! (0.32) \!\!\!\! & \!\!\!\! (16) \!\!\!\! & \!\!\!\! (0.31) \!\!\!\! & \!\!\!\! (15) \!\!\!\! & \!\!\!\! (0.27) \!\!\!\! & \!\!\!\! (13) \!\!\!\! & \!\!\!\! (0.09) \!\!\!\! & \!\!\!\! (10) \!\!\!\! & \!\!\!\!\! (0.12) \!\!\!\!\! \\
\hline
\multirow{2}{*}{\!\!\!\! $4$ \!\!\!\!\!\!\!\!\!\!\!\! } & \!\!\!\!\! 22 \!\!\!\!\! & \!\!\!\! 1.14 \!\!\!\! & \!\!\!\! 20 \!\!\!\! & \!\!\!\! 1.19 \!\!\!\! & \!\!\!\! 18 \!\!\!\! & \!\!\!\! 0.92 \!\!\!\! & \!\!\!\! 15 \!\!\!\! & \!\!\!\! 0.74 \!\!\!\! & \!\!\!\! 12 \!\!\!\! & \!\!\!\! 0.44 \!\!\!\! & \!\!\!\! 11 \!\!\!\! & \!\!\!\! 0.57 \!\!\!\! & \!\!\!\! 10 \!\!\!\! & \!\!\!\!\! 0.52 \!\!\!\!\! \\
& \!\!\!\!\! (15) \!\!\!\!\! & \!\!\!\! (0.89) \!\!\!\! & \!\!\!\! (19) \!\!\!\! & \!\!\!\! (1.10) \!\!\!\! & \!\!\!\! (18) \!\!\!\! & \!\!\!\! (0.98) \!\!\!\! & \!\!\!\! (16) \!\!\!\! & \!\!\!\! (0.78) \!\!\!\! & \!\!\!\! (16) \!\!\!\! & \!\!\!\! (0.97) \!\!\!\! & \!\!\!\! (15) \!\!\!\! & \!\!\!\! (0.80) \!\!\!\! & \!\!\!\! (14) \!\!\!\! & \!\!\!\!\! (0.66) \!\!\!\!\! \\
\hline
\multirow{2}{*}{\!\!\!\! $5$ \!\!\!\!\!\!\!\!\!\!\!\! } & \!\!\!\!\! 24 \!\!\!\!\! & \!\!\!\! 4.81 \!\!\!\! & \!\!\!\! 21 \!\!\!\! & \!\!\!\! 4.06 \!\!\!\! & \!\!\!\! 20 \!\!\!\! & \!\!\!\! 3.62 \!\!\!\! & \!\!\!\! 17 \!\!\!\! & \!\!\!\! 2.73 \!\!\!\! & \!\!\!\! 15 \!\!\!\! & \!\!\!\! 2.37 \!\!\!\! & \!\!\!\! 12 \!\!\!\! & \!\!\!\! 1.55 \!\!\!\! & \!\!\!\! 12 \!\!\!\! & \!\!\!\!\! 1.59 \!\!\!\!\! \\
& \!\!\!\!\! (20) \!\!\!\!\! & \!\!\!\! (4.07) \!\!\!\! & \!\!\!\! (20) \!\!\!\! & \!\!\!\! (4.06) \!\!\!\! & \!\!\!\! (23) \!\!\!\! & \!\!\!\! (4.69) \!\!\!\! & \!\!\!\! (16) \!\!\!\! & \!\!\!\! (3.18) \!\!\!\! & \!\!\!\! (16) \!\!\!\! & \!\!\!\! (2.90) \!\!\!\! & \!\!\!\! (16) \!\!\!\! & \!\!\!\! (2.88) \!\!\!\! & \!\!\!\! (15) \!\!\!\! & \!\!\!\!\! (2.15) \!\!\!\!\! \\
\hline
\multirow{2}{*}{\!\!\!\! $6$ \!\!\!\!\!\!\!\!\!\!\!\! } & \!\!\!\!\! 26 \!\!\!\!\! & \!\!\!\! 24.2 \!\!\!\! & \!\!\!\! 25 \!\!\!\! & \!\!\!\! 22.8 \!\!\!\! & \!\!\!\! 20 \!\!\!\! & \!\!\!\! 18.0 \!\!\!\! & \!\!\!\! 18 \!\!\!\! & \!\!\!\! 16.1 \!\!\!\! & \!\!\!\! 17 \!\!\!\! & \!\!\!\! 14.1 \!\!\!\! & \!\!\!\! 16 \!\!\!\! & \!\!\!\! 12.4 \!\!\!\! & \!\!\!\! 13 \!\!\!\! & \!\!\!\!\! 8.43 \!\!\!\!\! \\
& \!\!\!\!\! (26) \!\!\!\!\! & \!\!\!\! (24.0) \!\!\!\! & \!\!\!\! (33) \!\!\!\! & \!\!\!\! (30.2) \!\!\!\! & \!\!\!\! (23) \!\!\!\! & \!\!\!\! (20.9) \!\!\!\! & \!\!\!\! (19) \!\!\!\! & \!\!\!\! (17.2) \!\!\!\! & \!\!\!\! (16) \!\!\!\! & \!\!\!\! (14.2) \!\!\!\! & \!\!\!\! (16) \!\!\!\! & \!\!\!\! (13.5) \!\!\!\! & \!\!\!\! (15) \!\!\!\! & \!\!\!\!\! (12.2) \!\!\!\!\! \\
\hline
\multirow{2}{*}{\!\!\!\! $7$ \!\!\!\!\!\!\!\!\!\!\!\! } & \!\!\!\!\! 31 \!\!\!\!\! & \!\!\!\! 112 \!\!\!\! & \!\!\!\! 25 \!\!\!\! & \!\!\!\! 89.6 \!\!\!\! & \!\!\!\! 23 \!\!\!\! & \!\!\!\! 83.3 \!\!\!\! & \!\!\!\! 20 \!\!\!\! & \!\!\!\! 68.6 \!\!\!\! & \!\!\!\! 17 \!\!\!\! & \!\!\!\! 57.9 \!\!\!\! & \!\!\!\! 16 \!\!\!\! & \!\!\!\! 52.0 \!\!\!\! & \!\!\!\! 16 \!\!\!\! & \!\!\!\!\! 50.8 \!\!\!\!\! \\
& \!\!\!\!\! (27) \!\!\!\!\! & \!\!\!\! (97.2) \!\!\!\! & \!\!\!\! (27) \!\!\!\! & \!\!\!\! (96.2) \!\!\!\! & \!\!\!\! (29) \!\!\!\! & \!\!\!\! (103) \!\!\!\! & \!\!\!\! (22) \!\!\!\! & \!\!\!\! (77.1) \!\!\!\! & \!\!\!\! (17) \!\!\!\! & \!\!\!\! (59.1) \!\!\!\! & \!\!\!\! (14) \!\!\!\! & \!\!\!\! (47.9) \!\!\!\! & \!\!\!\! (17) \!\!\!\! & \!\!\!\!\! (57.1) \!\!\!\!\! \\
\hline
\multirow{2}{*}{\!\!\!\! $8$ \!\!\!\!\!\!\!\!\!\!\!\! } & \!\!\!\!\! 40 \!\!\!\!\! & \!\!\!\! 665 \!\!\!\! & \!\!\!\! 32 \!\!\!\! & \!\!\!\! 526 \!\!\!\! & \!\!\!\! 28 \!\!\!\! & \!\!\!\! 457 \!\!\!\! & \!\!\!\! 22 \!\!\!\! & \!\!\!\! 360 \!\!\!\! & \!\!\!\! 19 \!\!\!\! & \!\!\!\! 304 \!\!\!\! & \!\!\!\! 18 \!\!\!\! & \!\!\!\! 281 \!\!\!\! & \!\!\!\! 16 \!\!\!\! & \!\!\!\!\! 257 \!\!\!\!\! \\
 & \!\!\!\!\! (36) \!\!\!\!\! & \!\!\!\! (594) \!\!\!\! & \!\!\!\! (37) \!\!\!\! & \!\!\!\! (612) \!\!\!\! & \!\!\!\! (36) \!\!\!\! & \!\!\!\! (594) \!\!\!\! & \!\!\!\! (26) \!\!\!\! & \!\!\!\! (428) \!\!\!\! & \!\!\!\! (20) \!\!\!\! & \!\!\!\! (330) \!\!\!\! & \!\!\!\! (18) \!\!\!\! & \!\!\!\! (296) \!\!\!\! & \!\!\!\! (15) \!\!\!\! & \!\!\!\!\! (246) \!\!\!\!\! \\
\hline
\end{tabular}
\end{center}
\end{small}
\end{table}

\begin{table}[!ht]
\caption{Average GMRES iterations and CPU times for stationary Navier--Stokes control problem, for $\nu=\frac{1}{100}$ and a range of $l$, $\beta$.}\label{Stationary_NS_test1_table2}
\begin{small}
\begin{center}
\begin{tabular}{|c||c|c|c|c|c|c|c|c|c|c|c|c|c|c|}
\hline
\multicolumn{1}{|c||}{\phantom{$1$}} & \multicolumn{2}{c|}{$\beta=10^0$}& \multicolumn{2}{c|}{$\beta=10^{-1}$}& \multicolumn{2}{c|}{$\beta=10^{-2}$}& \multicolumn{2}{c|}{$\beta=10^{-3}$}& \multicolumn{2}{c|}{$\beta=10^{-4}$}& \multicolumn{2}{c|}{$\beta=10^{-5}$}& \multicolumn{2}{c|}{$\beta=10^{-6}$}\\
\cline{2-15}
\!\!\!\! $l$ \!\!\!\! & \!\!\!\! $\texttt{it}$ \!\!\!\! & \!\!\!\! CPU \!\!\!\! & \!\!\!\! $\texttt{it}$ \!\!\!\! & \!\!\!\! CPU \!\!\!\! & \!\!\!\! $\texttt{it}$ \!\!\!\! & \!\!\!\! CPU \!\!\!\! & \!\!\!\! $\texttt{it}$ \!\!\!\! & \!\!\!\! CPU \!\!\!\! & \!\!\!\! $\texttt{it}$ \!\!\!\! & \!\!\!\! CPU \!\!\!\! & \!\!\!\! $\texttt{it}$ \!\!\!\! & \!\!\!\! CPU \!\!\!\! & \!\!\!\! $\texttt{it}$ \!\!\!\! & \!\!\!\! CPU \!\!\!\! \\
\hline
\hline
\!\!\!\! $3$ \!\!\!\! & \!\!\!\! 38 \!\!\!\! & \!\!\!\! 0.79 \!\!\!\! & \!\!\!\! 24 \!\!\!\! & \!\!\!\! 0.30 \!\!\!\! & \!\!\!\! 13 \!\!\!\! & \!\!\!\! 0.20 \!\!\!\! & \!\!\!\! 11 \!\!\!\! & \!\!\!\! 0.19 \!\!\!\! & \!\!\!\! 11 \!\!\!\! & \!\!\!\! 0.18 \!\!\!\! & \!\!\!\! 9 \!\!\!\! & \!\!\!\! 0.13 \!\!\!\! & \!\!\!\! 9 \!\!\!\! & \!\!\!\! 0.14 \!\!\!\! \\
\hline
\!\!\!\! $4$ \!\!\!\! & \!\!\!\! 31 \!\!\!\! & \!\!\!\! 1.84 \!\!\!\! & \!\!\!\! 24 \!\!\!\! & \!\!\!\! 1.28 \!\!\!\! & \!\!\!\! 18 \!\!\!\! & \!\!\!\! 0.96 \!\!\!\! & \!\!\!\! 12 \!\!\!\! & \!\!\!\! 0.43 \!\!\!\! & \!\!\!\! 11 \!\!\!\! & \!\!\!\! 0.61 \!\!\!\! & \!\!\!\! 11 \!\!\!\! & \!\!\!\! 0.59 \!\!\!\! & \!\!\!\! 10 \!\!\!\! & \!\!\!\! 0.50 \!\!\!\! \\
\hline
\!\!\!\! $5$ \!\!\!\! & \!\!\!\! 29 \!\!\!\! & \!\!\!\! 5.47 \!\!\!\! & \!\!\!\! 23 \!\!\!\! & \!\!\!\! 4.23 \!\!\!\! & \!\!\!\! 20 \!\!\!\! & \!\!\!\! 3.23 \!\!\!\! & \!\!\!\! 16 \!\!\!\! & \!\!\!\! 2.32 \!\!\!\! & \!\!\!\! 12 \!\!\!\! & \!\!\!\! 1.60 \!\!\!\! & \!\!\!\! 11 \!\!\!\! & \!\!\!\! 1.58 \!\!\!\! & \!\!\!\! 11 \!\!\!\! & \!\!\!\! 1.64 \!\!\!\! \\
\hline
\!\!\!\! $6$ \!\!\!\! & \!\!\!\! 31 \!\!\!\! & \!\!\!\! 28.0 \!\!\!\! & \!\!\!\! 27 \!\!\!\! & \!\!\!\! 23.6 \!\!\!\! & \!\!\!\! 22 \!\!\!\! & \!\!\!\! 18.5 \!\!\!\! & \!\!\!\! 18 \!\!\!\! & \!\!\!\! 14.6 \!\!\!\! & \!\!\!\! 15 \!\!\!\! & \!\!\!\! 11.1 \!\!\!\! & \!\!\!\! 12 \!\!\!\! & \!\!\!\! 8.05 \!\!\!\! & \!\!\!\! 11 \!\!\!\! & \!\!\!\! 8.14 \!\!\!\! \\
\hline
\!\!\!\! $7$ \!\!\!\! & \!\!\!\! 32 \!\!\!\! & \!\!\!\! 116 \!\!\!\! & \!\!\!\! 27 \!\!\!\! & \!\!\!\! 93.6 \!\!\!\! & \!\!\!\! 24 \!\!\!\! & \!\!\!\! 83.9 \!\!\!\! & \!\!\!\! 20 \!\!\!\! & \!\!\!\! 69.6 \!\!\!\! & \!\!\!\! 17 \!\!\!\! & \!\!\!\! 58.7 \!\!\!\! & \!\!\!\! 15 \!\!\!\! & \!\!\!\! 46.7 \!\!\!\! & \!\!\!\! 13 \!\!\!\! & \!\!\!\! 35.9 \!\!\!\! \\
\hline
\!\!\!\! $8$ \!\!\!\! & \!\!\!\! 38 \!\!\!\! & \!\!\!\! 627 \!\!\!\! & \!\!\!\! 32 \!\!\!\! & \!\!\!\! 528 \!\!\!\! & \!\!\!\! 27 \!\!\!\! & \!\!\!\! 437 \!\!\!\! & \!\!\!\! 22 \!\!\!\! & \!\!\!\! 351 \!\!\!\! & \!\!\!\! 19 \!\!\!\! & \!\!\!\! 298 \!\!\!\! & \!\!\!\! 17 \!\!\!\! & \!\!\!\! 276 \!\!\!\! & \!\!\!\! 15 \!\!\!\! & \!\!\!\! 236 \!\!\!\! \\
\hline
\end{tabular}
\end{center}
\end{small}
\end{table}

\begin{table}[!ht]
\caption{Average GMRES iterations and CPU times for stationary Navier--Stokes control problem, for $\nu=\frac{1}{500}$ and a range of $l$, $\beta$.}\label{Stationary_NS_test1_table3}
\begin{small}
\begin{center}
{\begin{tabular}{|c||c|c|c|c|c|c|c|c|c|c|c|c|c|c|}
\hline
\multicolumn{1}{|c||}{\phantom{$1$}} & \multicolumn{2}{c|}{$\beta=10^0$}& \multicolumn{2}{c|}{$\beta=10^{-1}$}& \multicolumn{2}{c|}{$\beta=10^{-2}$}& \multicolumn{2}{c|}{$\beta=10^{-3}$}& \multicolumn{2}{c|}{$\beta=10^{-4}$}& \multicolumn{2}{c|}{$\beta=10^{-5}$}& \multicolumn{2}{c|}{$\beta=10^{-6}$}\\
\cline{2-15}
\!\!\!\! $l$  \!\!\!\! & \!\!\!\! $\texttt{it}$ \!\!\!\! & \!\!\!\! CPU \!\!\!\! & \!\!\!\! $\texttt{it}$ \!\!\!\! & \!\!\!\! CPU \!\!\!\! & \!\!\!\! $\texttt{it}$ \!\!\!\! & \!\!\!\! CPU \!\!\!\! & \!\!\!\! $\texttt{it}$ \!\!\!\! & \!\!\!\! CPU \!\!\!\! & \!\!\!\! $\texttt{it}$ \!\!\!\! & \!\!\!\! CPU \!\!\!\! & \!\!\!\! $\texttt{it}$ \!\!\!\! & \!\!\!\! CPU \!\!\!\! & \!\!\!\! $\texttt{it}$ \!\!\!\! & \!\!\!\! CPU \!\!\!\! \\
\hline
\hline
\!\!\!\! $3$ \!\!\!\! & \!\!\!\! 77$\dagger$\footnotemark \!\!\!\! & \!\!\!\! 1.46$\dagger$ \!\!\!\! & \!\!\!\! 23 \!\!\!\! & \!\!\!\! 0.44 \!\!\!\! & \!\!\!\! 13 \!\!\!\! & \!\!\!\! 0.25 \!\!\!\! & \!\!\!\! 11 \!\!\!\! & \!\!\!\! 0.24 \!\!\!\! & \!\!\!\! 10 \!\!\!\! & \!\!\!\! 0.19 \!\!\!\! & \!\!\!\! 9 \!\!\!\! & \!\!\!\! 0.14 \!\!\!\! & \!\!\!\! 9 \!\!\!\! & \!\!\!\! 0.14 \!\!\!\! \\
\hline
\!\!\!\! $4$ \!\!\!\! & \!\!\!\! 86$\dagger$ \!\!\!\! & \!\!\!\! 4.18$\dagger$ \!\!\!\! & \!\!\!\! 48$\dagger$ \!\!\!\! & \!\!\!\! 2.02$\dagger$ \!\!\!\! & \!\!\!\! 18 \!\!\!\! & \!\!\!\! 1.01 \!\!\!\! & \!\!\!\! 12 \!\!\!\! & \!\!\!\! 0.68 \!\!\!\! & \!\!\!\! 11 \!\!\!\! & \!\!\!\! 0.67 \!\!\!\! & \!\!\!\! 10 \!\!\!\! & \!\!\!\! 0.57 \!\!\!\! & \!\!\!\! 10 \!\!\!\! & \!\!\!\! 0.55 \!\!\!\! \\
\hline
\!\!\!\! $5$ \!\!\!\! & \!\!\!\! 74 \!\!\!\! & \!\!\!\! 15.4 \!\!\!\! & \!\!\!\! 49 \!\!\!\! & \!\!\!\! 8.02 \!\!\!\! & \!\!\!\! 28 \!\!\!\! & \!\!\!\! 4.13 \!\!\!\! & \!\!\!\! 14 \!\!\!\! & \!\!\!\! 2.08 \!\!\!\! & \!\!\!\! 12 \!\!\!\! & \!\!\!\! 1.87 \!\!\!\! & \!\!\!\! 11 \!\!\!\! & \!\!\!\! 1.73 \!\!\!\! & \!\!\!\! 10 \!\!\!\! & \!\!\!\! 1.50 \!\!\!\! \\
\hline
\!\!\!\! $6$ \!\!\!\! & \!\!\!\! 57 \!\!\!\! & \!\!\!\! 55.1 \!\!\!\! & \!\!\!\! 39 \!\!\!\! & \!\!\!\! 34.2 \!\!\!\! & \!\!\!\! 27 \!\!\!\! & \!\!\!\! 22.0 \!\!\!\! & \!\!\!\! 20 \!\!\!\! & \!\!\!\! 14.6 \!\!\!\! & \!\!\!\! 13 \!\!\!\! & \!\!\!\! 8.93 \!\!\!\! & \!\!\!\! 11 \!\!\!\! & \!\!\!\! 8.86 \!\!\!\! & \!\!\!\! 11 \!\!\!\! & \!\!\!\! 8.63 \!\!\!\! \\
\hline
\!\!\!\! $7$ \!\!\!\! & \!\!\!\! 54 \!\!\!\! & \!\!\!\! 192 \!\!\!\! & \!\!\!\! 32 \!\!\!\! & \!\!\!\! 111 \!\!\!\! & \!\!\!\! 27 \!\!\!\! & \!\!\!\! 90.6 \!\!\!\! & \!\!\!\! 21 \!\!\!\! & \!\!\!\! 67.3 \!\!\!\! & \!\!\!\! 17 \!\!\!\! & \!\!\!\! 49.4 \!\!\!\! & \!\!\!\! 12 \!\!\!\! & \!\!\!\! 32.5 \!\!\!\! & \!\!\!\! 12 \!\!\!\! & \!\!\!\! 37.2 \!\!\!\! \\
\hline
\!\!\!\! $8$ \!\!\!\! & \!\!\!\! 53 \!\!\!\! & \!\!\!\! 878 \!\!\!\! & \!\!\!\! 34 \!\!\!\! & \!\!\!\! 561 \!\!\!\! & \!\!\!\! 29 \!\!\!\! & \!\!\!\! 472 \!\!\!\! & \!\!\!\! 23 \!\!\!\! & \!\!\!\! 369 \!\!\!\! & \!\!\!\! 19 \!\!\!\! & \!\!\!\! 292 \!\!\!\! & \!\!\!\! 16 \!\!\!\! & \!\!\!\! 232 \!\!\!\! & \!\!\!\! 13 \!\!\!\! & \!\!\!\! 172 \!\!\!\! \\
\hline
\end{tabular}}
\end{center}
\end{small}
\end{table}

\begin{table}[!ht]
\caption{Degrees of freedom (DoF) and number of Oseen iterations required for stationary Navier--Stokes control problem. In each cell are the Oseen iterations for the given $l$, $\nu$, and $\beta=10^{-j}$, $j=0,1,...,6$.}\label{Number_Oseen_Stationary_NS}
\begin{small}
\begin{center}
\begin{tabular}{|c|c||ccccccc|ccccccc|ccccccc|}
\hline
$l$ & DoF & \multicolumn{7}{c|}{$\nu=\frac{1}{20}$} & \multicolumn{7}{c|}{$\nu=\frac{1}{100}$} & \multicolumn{7}{c|}{$\nu=\frac{1}{500}$} \\
\hline
\hline
$3$ & 1062 & 5 & \!\!\!\!\! 5 & \!\!\!\!\! 5 & \!\!\!\!\! 5 & \!\!\!\!\! 4 & \!\!\!\!\! 4 & \!\!\!\!\! 3 & 13 & \!\!\!\!\! 9 & \!\!\!\!\! 7 & \!\!\!\!\! 5 & \!\!\!\!\! 4 & \!\!\!\!\! 4 & \!\!\!\!\! 3 & $\dagger$ & \!\!\!\!\! 20 & \!\!\!\!\! 7 & \!\!\!\!\! 5 & \!\!\!\!\! 4 & \!\!\!\!\! 4 & \!\!\!\!\! 3 \\
\hline
$4$ & 4422 & 5 & \!\!\!\!\! 5 & \!\!\!\!\! 4 & \!\!\!\!\! 4 & \!\!\!\!\! 4 & \!\!\!\!\! 4 & \!\!\!\!\! 4 & 8 & \!\!\!\!\! 6 & \!\!\!\!\! 6 & \!\!\!\!\! 5 & \!\!\!\!\! 4 & \!\!\!\!\! 4 & \!\!\!\!\! 4 & $\dagger$ & \!\!\!\!\! $\dagger$ & \!\!\!\!\! 9 & \!\!\!\!\! 5 & \!\!\!\!\! 4 & \!\!\!\!\! 4 & \!\!\!\!\! 4 \\
\hline
$5$ & 18,054 & 4 & \!\!\!\!\! 4 & \!\!\!\!\! 4 & \!\!\!\!\! 4 & \!\!\!\!\! 4 & \!\!\!\!\! 4 & \!\!\!\!\! 3 & 7 & \!\!\!\!\! 5 & \!\!\!\!\! 5 & \!\!\!\!\! 4 & \!\!\!\!\! 4 & \!\!\!\!\! 4 & \!\!\!\!\! 3 & 16 & \!\!\!\!\! 10 & \!\!\!\!\! 8 & \!\!\!\!\! 6 & \!\!\!\!\! 4 & \!\!\!\!\! 4 & \!\!\!\!\! 3 \\
\hline
$6$ & 72,966 & 4 & \!\!\!\!\! 4 & \!\!\!\!\! 4 & \!\!\!\!\! 3 & \!\!\!\!\! 3 & \!\!\!\!\! 3 & \!\!\!\!\! 3 & 6 & \!\!\!\!\! 4 & \!\!\!\!\! 4 & \!\!\!\!\! 4 & \!\!\!\!\! 4 & \!\!\!\!\! 4 & \!\!\!\!\! 3 & 11 & \!\!\!\!\! 6 & \!\!\!\!\! 5 & \!\!\!\!\! 5 & \!\!\!\!\! 4 & \!\!\!\!\! 4 & \!\!\!\!\! 3 \\
\hline
$7$ & 293,382 & 4 & \!\!\!\!\! 3 & \!\!\!\!\! 3 & \!\!\!\!\! 3 & \!\!\!\!\! 3 & \!\!\!\!\! 3 & \!\!\!\!\! 3 & 5 & \!\!\!\!\! 4 & \!\!\!\!\! 3 & \!\!\!\!\! 3 & \!\!\!\!\! 3 & \!\!\!\!\! 3 & \!\!\!\!\! 3 & 8 & \!\!\!\!\! 4 & \!\!\!\!\! 4 & \!\!\!\!\! 4 & \!\!\!\!\! 4 & \!\!\!\!\! 4 & \!\!\!\!\! 3 \\
\hline
$8$ & 1,176,582 & 3 & \!\!\!\!\! 3 & \!\!\!\!\! 3 & \!\!\!\!\! 3 & \!\!\!\!\! 3 & \!\!\!\!\! 3 & \!\!\!\!\! 3 & 4 & \!\!\!\!\! 3 & \!\!\!\!\! 3 & \!\!\!\!\! 3 & \!\!\!\!\! 3 & \!\!\!\!\! 3 & \!\!\!\!\! 3 & 5 & \!\!\!\!\! 3 & \!\!\!\!\! 3 & \!\!\!\!\! 3 & \!\!\!\!\! 3 & \!\!\!\!\! 3 & \!\!\!\!\! 3 \\
\hline
\end{tabular}
\end{center}
\end{small}
\end{table}

Tables \ref{Stationary_NS_test1_table1}--\ref{Stationary_NS_test1_table3} demonstrate the robustness of our proposed preconditioner. The numbers of iterations show only a mild dependence on the viscosity $\nu$, and a slight increase only for large values of $\beta$. The CPU time scales approximately linearly with respect to the dimension of the systems, with a marginal increase for very fine grids; in this case we observe that the AGMG multigrid routine does not scale exactly linearly. Table \ref{Number_Oseen_Stationary_NS} shows that the number of Oseen iterations strongly depends on the viscosity $\nu$, as expected as the non-linear term becoming more dominant for smaller $\nu$; however, as the grid is refined the number of outer iterations decreases. We also note that the number of non-linear iterations increases for larger values of $\beta$ and coarser grids.
\footnotetext{$\dagger$ means that the outer (Oseen) iteration did not converge in 20 iterations. The average number of GMRES iterations and CPU time is evaluated over the first 10 Oseen iterations.}

\subsection{Crank--Nicolson for instationary Stokes control}\label{6_2}
We now test our solver on an instationary Stokes control problem, which allows us to verify the predicted order of convergence of the Crank--Nicolson method. We take $t_f=2$ and
\begin{small}
\begin{align*}
\vec{v}_d(x_1,x_2,t) = {}&  4 \beta \big[x_2 \! \big(2 (3x_1^2- \! 1)(x_2^2-  1)\! + \! 3 (x_1^2-  1)^2\big) \!, -x_1 \! \big( 3 (x_2^2- \! 1)^2 \! + \! 2 (x_1^2- \! 1)(3x_2^2 \! - \! 1) \big) \! \big]^{\! \top} \\
& + e^{t_f - t} \big[20x_1 x_2^3 + 2 \beta x_2\big( (x_1^2-1)^2(x_2^2-7)-4(3x_1^2-1)(x_2^2-1)+2 \big),\\
& \phantom{+ e^{t_f - t} \big[\;\;} 5 (x_1^4 - x_2^4)-2\beta x_1 \big( (x_2^2-1)^2(x_1^2-7)-4(x_1^2-1)(3x_2^2-1)-2 \big)\big]^{\! \top},\\
\vec{f}(x_1,x_2,t) = {}& e^{t_f - t}\big[-20x_1x_2^3-2x_2(x_1^2-1)^2(x_2^2-1), 5(x_2^4-x_1^4) \! +2x_1(x_1^2-1)(x_2^2- 1)^2\big]^{\! \top}\\
& + \big[2x_2(x_1^2-1)^2(x_2^2-1),-2x_1(x_1^2-1)(x_2^2-1)^2\big]^{\! \top}.
\end{align*}
\end{small}
The analytic solutions for this problem are:
\begin{equation*}
\begin{array}{l}
\vspace{0.5ex}
\vec{v}(x_1,x_2,t) = e^{t_f - t} [20x_1 x_2^3, 5 x_1^4 - 5x_2^4]^\top,\\
\vspace{0.5ex}
p(x_1,x_2,t) =e^{t_f - t} \left( 60x_1^2x_2 - 20 x_2^3 \right)+\text{constant},\\
\vspace{0.5ex}
\vec{\zeta}(x_1,x_2,t) =\beta \left(e^{t_f - t}-1 \right) [2x_2(x_1^2-1)^2(x_2^2-1),-2x_1(x_1^2-1)(x_2^2-1)^2]^\top,\\
\mu(x_1,x_2,t) = \beta e^{t_f - t} (4 x_1 x_2)+\text{constant},
\end{array}
\end{equation*}
with initial and boundary condition obtained from this $\vec{v}$. In Table \ref{Instationary_Stokes_CN} we report the level of refinement $l$, the number of GMRES iterations\footnote{For this problem we run GMRES until a relative reduction on the residual of $10^{-9}$ is achieved.}, the CPU time, and the resulting errors for different values of $\beta$. For level of refinement $l$ we divide the time interval into subintervals of length $2^{1-l}$ and consider a spatial uniform grid of refinement level $l$. The error is evaluated in the $L^\infty(L^2)$ norm, approximated for $\vec{v}$ as
\begin{displaymath}
\vec{v}_{\text{err}} = \max_{n} \: \: \left[(\boldsymbol{v}_n - \boldsymbol{v}_{\text{sol}, n})^\top \mathbf{M} \: (\boldsymbol{v}_n - \boldsymbol{v}_{\text{sol}, n})\right]^{1/2},
\end{displaymath}
where $\boldsymbol{v}_{\text{sol},n}$ is the discretized exact solution for $\vec{v}$ at time $t_n$. In the same way we define the error for the adjoint velocity $\vec{\zeta}_{\text{err}}$. 
Further, for the total size of the systems solved we refer to Section \ref{6_3_2}.

\begin{table}[!ht]
\caption{Average GMRES iterations, average CPU times, and errors for instationary Stokes control problem, solved using Crank--Nicolson, for a range of $l$, $\beta$.}\label{Instationary_Stokes_CN}
\begin{small}
\begin{center}
{\begin{tabular}{|c||c|c|c|c||c|c|c|c||c|c|c|c|}
\hline
 & \multicolumn{4}{c||}{$\beta=10^0$}  & \multicolumn{4}{c||}{$\beta=10^{-2}$}  & \multicolumn{4}{c|}{$\beta=10^{-4}$} \\
\cline{2-13}
\!\!\!\! $l$ \!\!\!\! & \!\!\!\! $\texttt{it}$ \!\!\!\! & \!\!\!\! CPU \!\!\!\! & \!\!\!\! $\vec{v}_{\text{err}}$ \!\!\!\! & \!\!\!\! $\vec{\zeta}_{\text{err}}$ \!\!\!\! & \!\!\!\! $\texttt{it}$ \!\!\!\! & \!\!\!\! CPU \!\!\!\! & \!\!\!\! $\vec{v}_{\text{err}}$ \!\!\!\! & \!\!\!\! $\vec{\zeta}_{\text{err}}$ \!\!\!\! & \!\!\!\! $\texttt{it}$ \!\!\!\! & \!\!\!\! CPU \!\!\!\! & \!\!\!\! $\vec{v}_{\text{err}}$ \!\!\!\! & \!\!\!\! $\vec{\zeta}_{\text{err}}$ \!\!\!\!\!\!\!\!\! \phantom{\Large{I}} \\
\hline
\hline
\!\!\!\! $2$ \!\!\!\! & \!\!\!\! 22 \!\!\!\! & \!\!\!\! 0.85 \!\!\!\! & \!\!\!\! 4.76e-1 \!\!\!\! & \!\!\!\! 2.49e-1 \!\!\!\! & \!\!\!\! 22 \!\!\!\! & \!\!\!\! 0.79 \!\!\!\! & \!\!\!\! 5.66e-1 \!\!\!\! & \!\!\!\! 1.16e-1 \!\!\!\! & \!\!\!\! 16 \!\!\!\! & \!\!\!\! 0.73 \!\!\!\! & \!\!\!\! 8.63e0 \!\!\!\! & \!\!\!\! 5.45e-2 \!\!\!\! \\
\hline
\!\!\!\! $3$ \!\!\!\! & \!\!\!\! 22 \!\!\!\! & \!\!\!\! 4.28 \!\!\!\! & \!\!\!\! 3.34e-2 \!\!\!\! & \!\!\!\! 5.68e-2 \!\!\!\! & \!\!\!\! 22 \!\!\!\! & \!\!\!\! 4.24 \!\!\!\! & \!\!\!\! 7.07e-2 \!\!\!\! & \!\!\!\! 3.42e-2 \!\!\!\! & \!\!\!\! 19 \!\!\!\! & \!\!\!\! 3.39 \!\!\!\! & \!\!\!\! 2.47e0 \!\!\!\! & \!\!\!\! 2.67e-2 \!\!\!\! \\
\hline
\!\!\!\! $4$ \!\!\!\! & \!\!\!\! 23 \!\!\!\! & \!\!\!\! 23.9 \!\!\!\! & \!\!\!\! 2.25e-3 \!\!\!\! & \!\!\!\! 1.15e-2 \!\!\!\! & \!\!\!\! 24 \!\!\!\! & \!\!\!\! 24.1 \!\!\!\! & \!\!\!\! 7.35e-3 \!\!\!\! & \!\!\!\! 7.79e-3 \!\!\!\! & \!\!\!\! 20 \!\!\!\! & \!\!\!\! 23.2 \!\!\!\! & \!\!\!\! 3.73e-1 \!\!\!\! & \!\!\!\! 7.30e-3 \!\!\!\! \\
\hline
\!\!\!\! $5$ \!\!\!\! & \!\!\!\! 23 \!\!\!\! & \!\!\!\! 200 \!\!\!\! & \!\!\!\! 1.74e-4 \!\!\!\! & \!\!\!\! 2.15e-3 \!\!\!\! & \!\!\!\! 27 \!\!\!\! & \!\!\!\! 232 \!\!\!\! & \!\!\!\! 6.70e-4 \!\!\!\! & \!\!\!\! 1.59e-3 \!\!\!\! & \!\!\!\! 20 \!\!\!\! & \!\!\!\! 162 \!\!\!\! & \!\!\!\! 3.84e-2 \!\!\!\! & \!\!\!\! 1.55e-3 \!\!\!\! \\
\hline
\!\!\!\! $6$ \!\!\!\! & \!\!\!\! 26 \!\!\!\! & \!\!\!\! 2082 \!\!\!\! & \!\!\!\! 2.16e-5 \!\!\!\! & \!\!\!\! 4.00e-4 \!\!\!\! & \!\!\!\! 37 \!\!\!\! & \!\!\!\! 2960 \!\!\!\! & \!\!\!\! 5.97e-5 \!\!\!\! & \!\!\!\! 3.02e-4 \!\!\!\! & \!\!\!\! 23 \!\!\!\! & \!\!\!\! 1830 \!\!\!\! & \!\!\!\! 3.38e-3 \!\!\!\! & \!\!\!\! 3.00e-4 \!\!\!\! \\
\hline
\end{tabular}}
\end{center}
\end{small}
\end{table}

From the discretization errors reported in Table \ref{Instationary_Stokes_CN}, we first note that the method is converging at second-order, and we experience similar convergence behaviour for the pressure variables. Secondly, we note that the preconditioner behaves robustly with respect to the level of refinement $l$ and the regularization parameter $\beta$, with the number of iterations slightly increasing for very fine grids. The elapsed CPU time scales almost exactly linearly, aside from the multigrid routine for very fine grids.

\subsection{Instationary Navier--Stokes control}\label{6_3}
We now test our solver on the instationary Navier--Stokes control problem \eqref{Instationary_Navier_Stokes_control_functional}--\eqref{Instationary_Navier_Stokes_control_constraints}, where we set $t_f=2$, $\vec{f}(\mathbf{x},t)=\vec{0}$, the initial condition $\vec{v}_0(\mathbf{x}) = \vec{0}$, and boundary conditions
\begin{displaymath}
\vec{g}_D(\mathbf{x},t) = \left\{
\begin{array}{ll}\left[t,0\right]^\top & \mathrm{on} \: \partial \Omega_1 \times (0,1),\\
\left[1,0\right]^\top & \mathrm{on} \: \partial \Omega_1 \times [1,t_f),\\
\left[0,0\right]^\top & \mathrm{on} \: (\partial \Omega \setminus \partial\Omega_1) \times (0,t_f).
\end{array}
\right.
\end{displaymath}
We present results obtained by employing backward Euler and Crank--Nicolson discretizations in time. Setting $c_1 = 1 - \sqrt{ (\frac{100}{49}( x_1-  \frac{1}{2}))^2 + (\frac{100}{99}x_2)^2 }$ and $c_2 = 1 - \sqrt{ (\frac{100}{49}( x_1+ \frac{1}{2}))^2 + (\frac{100}{99}x_2)^2 }$, we seek the (divergence-free) desired state:
\begin{displaymath}
\vec{v}_d(\mathbf{x},t) = \left\{
\begin{array}{ll}
\vspace{0.5ex}
c_1 \cos(\frac{\pi t}{2}) \: [(\frac{100}{99})^2 x_2, -(\frac{100}{49})^2 (x_1-\frac{1}{2})]^\top & \mathrm{if} \: c_1 \geq 0,\\
\vspace{0.5ex}
c_2 \cos(\frac{\pi t}{2}) \: [-(\frac{100}{99})^2 x_2, (\frac{100}{49})^2 (x_1+\frac{1}{2})]^\top & \mathrm{if} \: c_2 \geq 0,\\
\left[0,0\right]^\top & \mathrm{otherwise}.
\end{array}
\right.
\end{displaymath}

\subsubsection{Backward Euler for instationary Navier--Stokes control}\label{6_3_1}
We first report the results obtained when employing the backward Euler scheme in time. We provide the average number of GMRES iterations together with the average elapsed CPU time in Tables \ref{Instationary_NS_with_BE_table1}--\ref{Instationary_NS_with_BE_table2}, and in Table \ref{Oseen_Instationary_NS_with_BE} the total dimensions of the systems solved and the Oseen iterations required, for different levels of refinements $l$, values of $\beta$, and viscosities $\nu$. Here, we choose the time-step $\tau=0.05$ (that is, $n_t=40$), while the level of refinement $l$ refers to a spatial uniform grid constructed as above.

\begin{table}[!ht]
\caption{Average GMRES iterations and CPU times for instationary Navier--Stokes control problem, with backward Euler in time ($\tau = 0.05$), for $\nu = \frac{1}{100}$ and a range of $l, \beta$.}\label{Instationary_NS_with_BE_table1}
\begin{small}
\begin{center}
\begin{tabular}{|c||c|c|c|c|c|c|c|c|c|c|c|c|c|c|}
\hline
\multicolumn{1}{|c||}{\phantom{$1$}} & \multicolumn{2}{c|}{$\beta=10^0$}& \multicolumn{2}{c|}{$\beta=10^{-1}$}& \multicolumn{2}{c|}{$\beta=10^{-2}$}& \multicolumn{2}{c|}{$\beta=10^{-3}$}& \multicolumn{2}{c|}{$\beta=10^{-4}$}& \multicolumn{2}{c|}{$\beta=10^{-5}$}& \multicolumn{2}{c|}{$\beta=10^{-6}$}\\
\cline{2-15}
\!\!\!\! $l$ \!\!\!\!  & \!\!\!\! $\texttt{it}$ \!\!\!\! & \!\!\!\! CPU \!\!\!\! & \!\!\!\! $\texttt{it}$ \!\!\!\! & \!\!\!\! CPU \!\!\!\! & \!\!\!\! $\texttt{it}$ \!\!\!\! & \!\!\!\! CPU \!\!\!\! & \!\!\!\! $\texttt{it}$ \!\!\!\! & \!\!\!\! CPU \!\!\!\! & \!\!\!\! $\texttt{it}$ \!\!\!\! & \!\!\!\! CPU \!\!\!\! & \!\!\!\! $\texttt{it}$ \!\!\!\! & \!\!\!\! CPU \!\!\!\! & \!\!\!\! $\texttt{it}$ \!\!\!\! & \!\!\!\! CPU \!\!\!\! \\
\hline
\hline
\!\!\!\! $2$ \!\!\!\! & \!\!\!\! 17 \!\!\!\! & \!\!\!\! 6.96 \!\!\!\! & \!\!\!\! 14 \!\!\!\! & \!\!\!\! 5.70 \!\!\!\! & \!\!\!\! 11 \!\!\!\! & \!\!\!\! 4.47 \!\!\!\! & \!\!\!\! 11 \!\!\!\! & \!\!\!\! 4.21 \!\!\!\! & \!\!\!\! 10 \!\!\!\! & \!\!\!\! 4.08 \!\!\!\! & \!\!\!\! 12 \!\!\!\! & \!\!\!\! 4.93 \!\!\!\! & \!\!\!\! 21 \!\!\!\! & \!\!\!\! 8.54 \!\!\!\! \\
\hline
\!\!\!\! $3$ \!\!\!\! & \!\!\!\! 22 \!\!\!\! & \!\!\!\! 18.8 \!\!\!\! & \!\!\!\! 19 \!\!\!\! & \!\!\!\! 17.7 \!\!\!\! & \!\!\!\! 14 \!\!\!\! & \!\!\!\! 13.2 \!\!\!\! & \!\!\!\! 11 \!\!\!\! & \!\!\!\! 10.9 \!\!\!\! & \!\!\!\! 11 \!\!\!\! & \!\!\!\! 10.6 \!\!\!\! & \!\!\!\! 13 \!\!\!\! & \!\!\!\! 11.8 \!\!\!\! & \!\!\!\! 21 \!\!\!\! & \!\!\!\! 19.7 \!\!\!\! \\
\hline
\!\!\!\! $4$ \!\!\!\! & \!\!\!\! 23 \!\!\!\! & \!\!\!\! 45.3 \!\!\!\! & \!\!\!\! 22 \!\!\!\! & \!\!\!\! 45.4 \!\!\!\! & \!\!\!\! 18 \!\!\!\! & \!\!\!\! 39.0 \!\!\!\! & \!\!\!\! 14 \!\!\!\! & \!\!\!\! 37.5 \!\!\!\! & \!\!\!\! 13 \!\!\!\! & \!\!\!\! 35.0 \!\!\!\! & \!\!\!\! 15 \!\!\!\! & \!\!\!\! 41.2 \!\!\!\! & \!\!\!\! 23 \!\!\!\! & \!\!\!\! 61.4 \!\!\!\! \\
\hline
\!\!\!\! $5$ \!\!\!\! & \!\!\!\! 22 \!\!\!\! & \!\!\!\! 190 \!\!\!\! & \!\!\!\! 22 \!\!\!\! & \!\!\!\! 196 \!\!\!\! & \!\!\!\! 19 \!\!\!\! & \!\!\!\! 169 \!\!\!\! & \!\!\!\! 17 \!\!\!\! & \!\!\!\! 148 \!\!\!\! & \!\!\!\! 15 \!\!\!\! & \!\!\!\! 126 \!\!\!\! &  \!\!\!\! 16 \!\!\!\! & \!\!\!\! 136 \!\!\!\! & \!\!\!\! 25 \!\!\!\! & \!\!\!\! 220 \!\!\!\! \\
\hline
\!\!\!\! $6$ \!\!\!\! & \!\!\!\! 25 \!\!\!\! & \!\!\!\! 1153 \!\!\!\! & \!\!\!\! 24 \!\!\!\! & \!\!\!\! 1099 \!\!\!\! & \!\!\!\! 21 \!\!\!\! & \!\!\!\! 979 \!\!\!\! & \!\!\!\! 18 \!\!\!\! & \!\!\!\! 809 \!\!\!\! & \!\!\!\! 17 \!\!\!\! & \!\!\!\! 729 \!\!\!\! & \!\!\!\! 17 \!\!\!\! & \!\!\!\! 685 \!\!\!\! & \!\!\!\! 25 \!\!\!\! & \!\!\!\! 1080 \!\!\!\! \\
\hline
\end{tabular}
\end{center}
\end{small}
\end{table}

\begin{table}[!ht]
\caption{Average GMRES iterations and CPU times for instationary Navier--Stokes control problem, with backward Euler in time ($\tau = 0.05$), for $\nu = \frac{1}{500}$ and a range of $l, \beta$.}\label{Instationary_NS_with_BE_table2}
\begin{small}
\begin{center}
\begin{tabular}{|c||c|c|c|c|c|c|c|c|c|c|c|c|c|c|}
\hline
\multicolumn{1}{|c||}{\phantom{$1$}} & \multicolumn{2}{c|}{$\beta=10^0$}& \multicolumn{2}{c|}{$\beta=10^{-1}$}& \multicolumn{2}{c|}{$\beta=10^{-2}$}& \multicolumn{2}{c|}{$\beta=10^{-3}$}& \multicolumn{2}{c|}{$\beta=10^{-4}$}& \multicolumn{2}{c|}{$\beta=10^{-5}$}& \multicolumn{2}{c|}{$\beta=10^{-6}$}\\
\cline{2-15}
\!\!\!\! $l$  \!\!\!\! & \!\!\!\! $\texttt{it}$ \!\!\!\! & \!\!\!\! CPU \!\!\!\! & \!\!\!\! $\texttt{it}$ \!\!\!\! & \!\!\!\! CPU \!\!\!\! & \!\!\!\! $\texttt{it}$ \!\!\!\! & \!\!\!\! CPU \!\!\!\! & \!\!\!\! $\texttt{it}$ \!\!\!\! & \!\!\!\! CPU \!\!\!\! & \!\!\!\! $\texttt{it}$ \!\!\!\! & \!\!\!\! CPU \!\!\!\! & \!\!\!\! $\texttt{it}$ \!\!\!\! & \!\!\!\! CPU \!\!\!\! & \!\!\!\! $\texttt{it}$ \!\!\!\! & \!\!\!\! CPU \!\!\!\! \\
\hline
\hline
\!\!\!\! $2$ \!\!\!\! & \!\!\!\! 16$\dagger$ \!\!\!\! & \!\!\!\! 6.66$\dagger$ \!\!\!\! & \!\!\!\! 14$\dagger$ \!\!\!\! & \!\!\!\! 5.55$\dagger$ \!\!\!\! & \!\!\!\! 11$\dagger$ \!\!\!\! & \!\!\!\! 4.43$\dagger$ \!\!\!\! & \!\!\!\! 10$\dagger$ \!\!\!\! & \!\!\!\! 3.92$\dagger$ \!\!\!\! & \!\!\!\! 10 \!\!\!\! & \!\!\!\! 3.84 \!\!\!\! & \!\!\!\! 12 \!\!\!\! & \!\!\!\! 4.91 \!\!\!\! & \!\!\!\! 22 \!\!\!\! & \!\!\!\! 8.87 \!\!\!\! \\
\hline
\!\!\!\! $3$ \!\!\!\! & \!\!\!\! 26$\dagger$ \!\!\!\! & \!\!\!\! 28.9$\dagger$ \!\!\!\! & \!\!\!\! 20$\dagger$ \!\!\!\! & \!\!\!\! 25.4$\dagger$ \!\!\!\! & \!\!\!\! 14 \!\!\!\! & \!\!\!\! 17.2 \!\!\!\! & \!\!\!\! 11 \!\!\!\! & \!\!\!\! 13.6 \!\!\!\! & \!\!\!\! 11 \!\!\!\! & \!\!\!\! 13.3 \!\!\!\! & \!\!\!\! 13 \!\!\!\! & \!\!\!\! 15.1 \!\!\!\! & \!\!\!\! 22 \!\!\!\! & \!\!\!\! 27.9 \!\!\!\! \\
\hline
\!\!\!\! $4$ \!\!\!\! & \!\!\!\! 43$\dagger$ \!\!\!\! & \!\!\!\! 130$\dagger$ \!\!\!\! & \!\!\!\! 34 \!\!\!\! & \!\!\!\! 120 \!\!\!\! & \!\!\!\! 17 \!\!\!\! & \!\!\!\! 64.5 \!\!\!\! & \!\!\!\! 13 \!\!\!\! & \!\!\!\! 50.6 \!\!\!\! & \!\!\!\! 12 \!\!\!\! & \!\!\!\! 47.8 \!\!\!\! & \!\!\!\! 14 \!\!\!\! & \!\!\!\! 53.6 \!\!\!\! & \!\!\!\! 24 \!\!\!\! & \!\!\!\! 82.7 \!\!\!\! \\
\hline
\!\!\!\! $5$ \!\!\!\! & \!\!\!\! 54 \!\!\!\! & \!\!\!\! 467 \!\!\!\! & \!\!\!\! 43 \!\!\!\! & \!\!\!\! 417 \!\!\!\! & \!\!\!\! 28 \!\!\!\! & \!\!\!\! 294 \!\!\!\! & \!\!\!\! 16 \!\!\!\! & \!\!\!\! 160 \!\!\!\! & \!\!\!\! 15 \!\!\!\! & \!\!\!\! 164 \!\!\!\! & \!\!\!\! 16 \!\!\!\! & \!\!\!\! 187 \!\!\!\! & \!\!\!\! 24 \!\!\!\! & \!\!\!\! 263 \!\!\!\! \\
\hline
\!\!\!\! $6$ \!\!\!\! & \!\!\!\! 39 \!\!\!\! & \!\!\!\! 1723 \!\!\!\! & \!\!\!\! 35 \!\!\!\! & \!\!\!\! 1552 \!\!\!\! & \!\!\!\! 27 \!\!\!\! & \!\!\!\! 1209 \!\!\!\! & \!\!\!\! 20 \!\!\!\! & \!\!\!\! 842 \!\!\!\! & \!\!\!\! 17 \!\!\!\! & \!\!\!\! 742 \!\!\!\! & \!\!\!\! 18 \!\!\!\! & \!\!\!\! 830 \!\!\!\! & \!\!\!\! 26 \!\!\!\! & \!\!\!\! 1165 \!\!\!\! \\
\hline
\end{tabular}
\end{center}
\end{small}
\end{table}

As for the stationary case, Tables \ref{Instationary_NS_with_BE_table1}--\ref{Instationary_NS_with_BE_table2} show robustness of the proposed preconditioner with respect to all the parameters involved. We note that the number of iterations increases slightly for small viscosities and large values of $\beta$. The elapsed CPU time scales almost linearly with the dimension of the system, except for very fine grids. We see from Table \ref{Oseen_Instationary_NS_with_BE} that the number of Oseen iterations increases for small values of $\nu$ and large values of $\beta$ when employing a coarse grid; however, as the grid is refined, the number of non-linear iterations decreases.

\begin{table}[!ht]
\caption{Degrees of freedom (DoF) and number of Oseen iterations required for instationary Navier--Stokes control problem, with backward Euler in time ($\tau \hspace{-0.2em} = \hspace{-0.2em} 0.05$). In each cell are the Oseen iterations for the given $l, \nu$, and $\beta \hspace{-0.2em} = \hspace{-0.2em} 10^{-j}, j \hspace{-0.2em} = \hspace{-0.2em} 0,1,...,6$.}\label{Oseen_Instationary_NS_with_BE}
\begin{small}
\begin{center}
\begin{tabular}{|c|c||ccccccc|ccccccc|}
\hline
$l$ & DoF & \multicolumn{7}{c|}{$\nu=\frac{1}{100}$} & \multicolumn{7}{c|}{$\nu=\frac{1}{500}$} \\
\hline
\hline
$2$ & 10,086 & 15 & 8 & 6 & 5 & 5 & 5 & 5 & $\dagger$ & $\dagger$ & $\dagger$ & $\dagger$ & 10 & 8 & 8 \\
\hline
$3$ & 43,542 & 9 & 8 & 6 & 5 & 5 & 5 & 5 & $\dagger$ & $\dagger$ & 7 & 6 & 6 & 6 & 6 \\
\hline
$4$ & 181,302 & 6 & 6 & 6 & 5 & 5 & 5 & 5 & $\dagger$ & 14 & 8 & 6 & 6 & 6 & 6 \\
\hline
$5$ & 740,214 & 5 & 5 & 5 & 4 & 4 & 4 & 4 & 8 & 8 & 7 & 5 & 5 & 5 & 5 \\
\hline
$6$ & 2,991,606 & 4 & 4 & 4 & 4 & 4 & 3 & 3 & 6 & 5 & 5 & 5 & 4 & 4 & 4 \\
\hline
\end{tabular}
\end{center}
\end{small}
\end{table}

\subsubsection{Crank--Nicolson for instationary Navier--Stokes control}\label{6_3_2}
\begin{figure}[!htb]
\centering
\caption{Solution plots for the instationary Navier--Stokes control problem, for $\nu=\frac{1}{100}$, $\beta=10^{-1}$, and $l=4$. Top left: velocity $\vec{v}$ at $t=1$. Top right: pressure $p$ at $t=1.0625$. Bottom: adjoint velocity $\vec{\zeta}$ at $t=1$.}
\label{fig:state_adjoint_pressure}
\vspace{-1.15em}
\begin{subfigure}[b]{0.34\textwidth}
\centering
\includegraphics[width=1.1\linewidth]{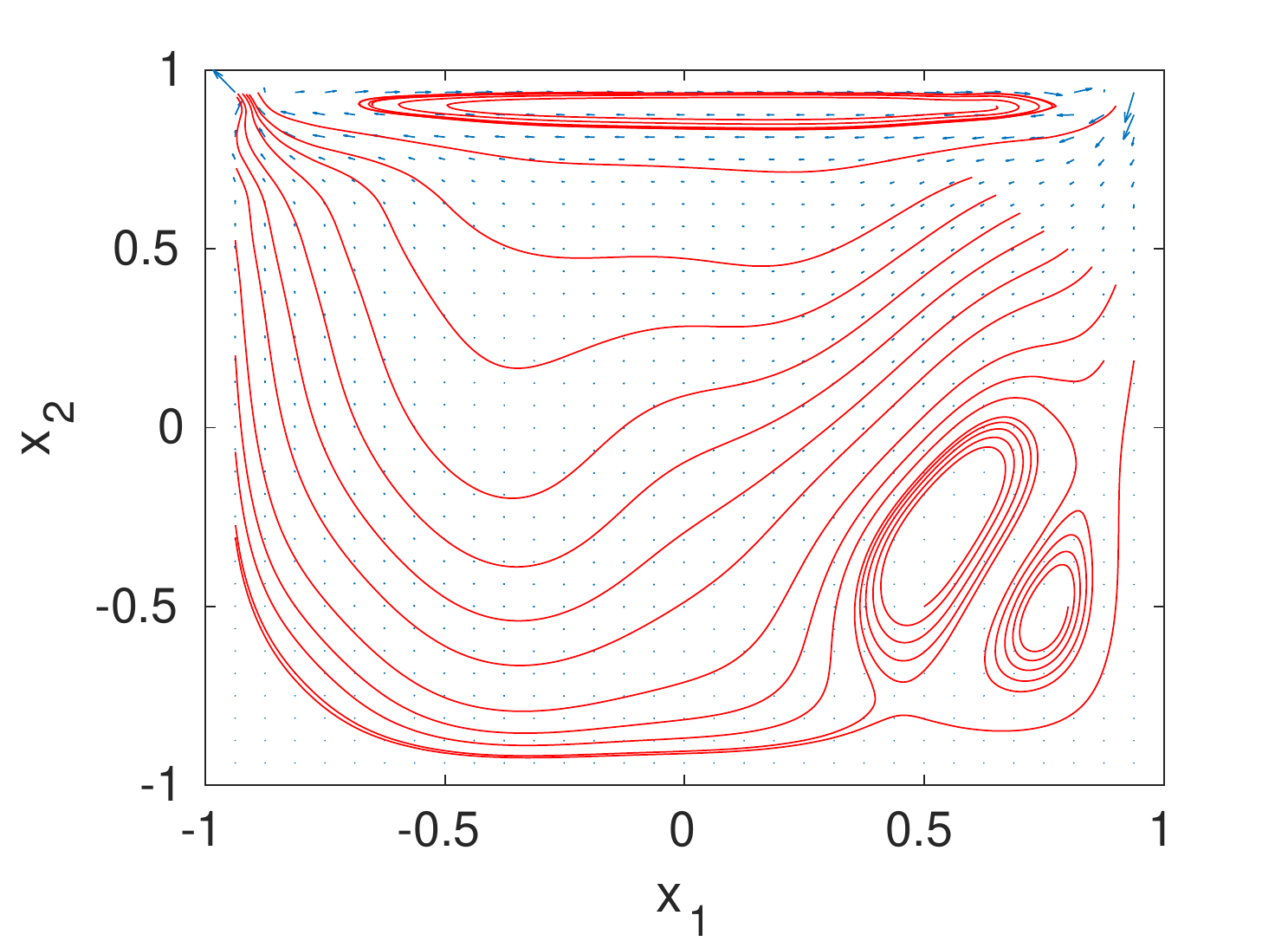}
\end{subfigure}\hspace{1.5em}
\begin{subfigure}[b]{0.34\textwidth}
\centering
\includegraphics[width=1.1\linewidth]{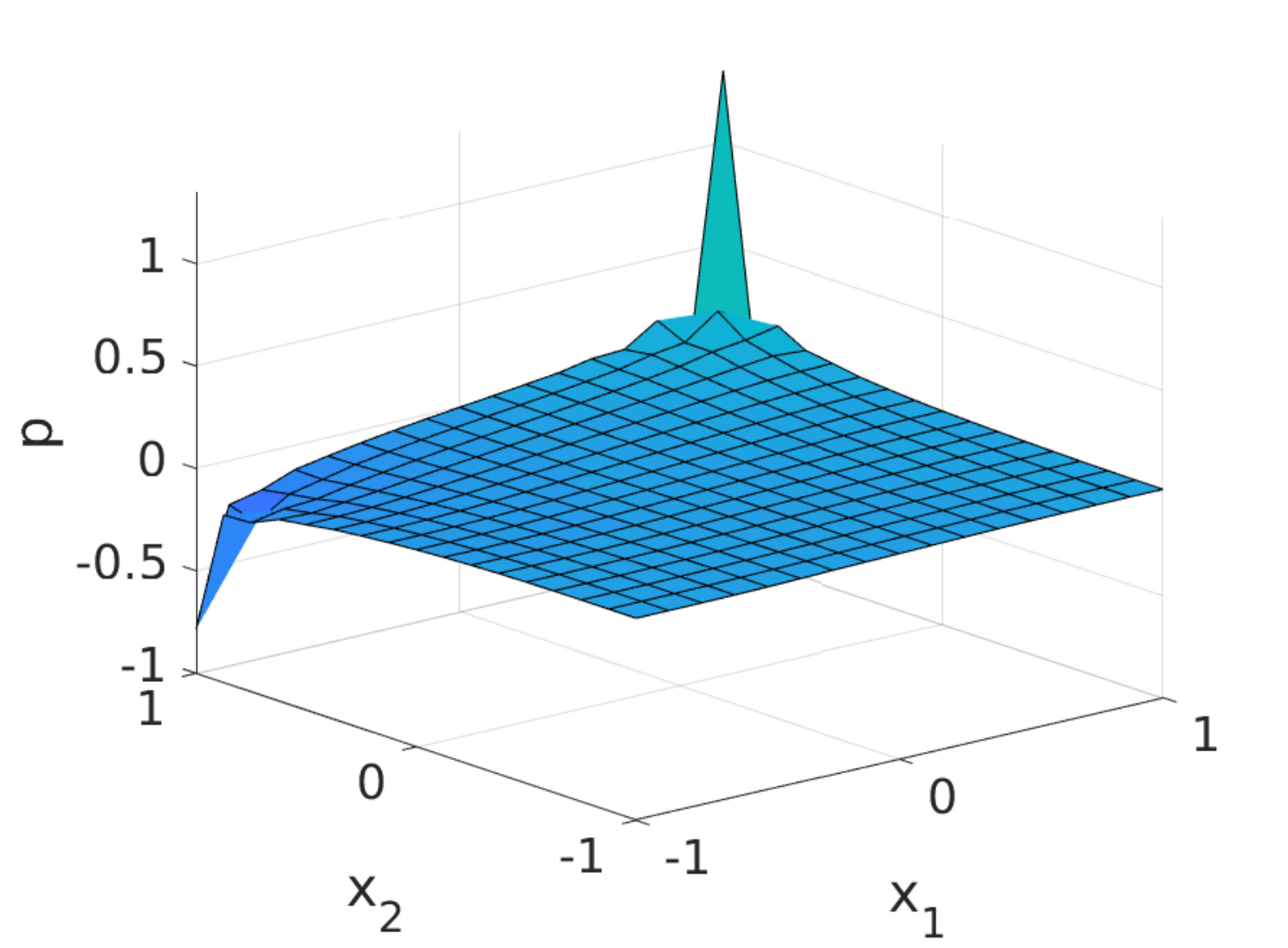}
\end{subfigure}
\begin{subfigure}[b]{0.34\textwidth}
\centering
\includegraphics[width=1.1\linewidth]{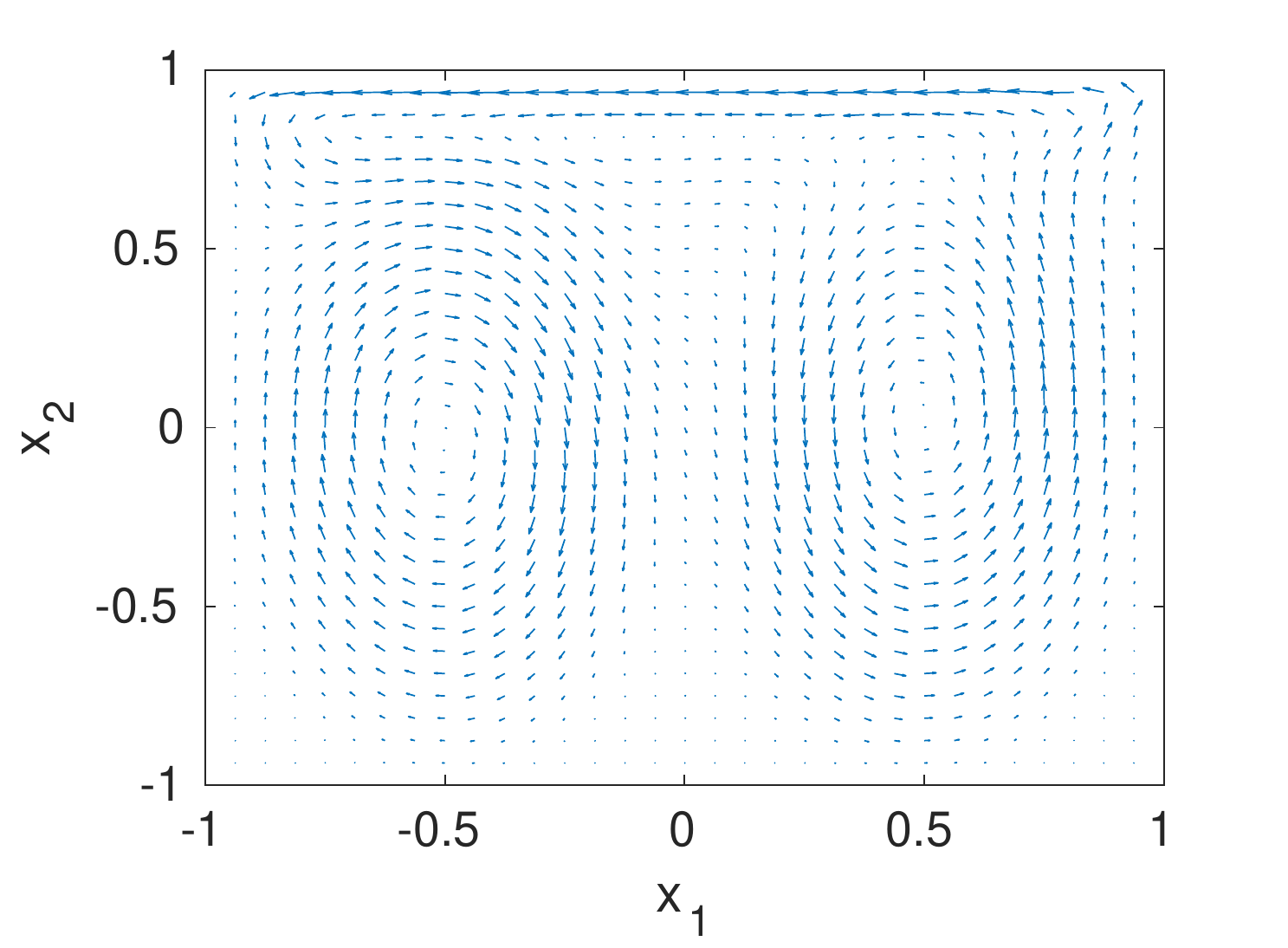}
\end{subfigure}
\vspace{-2em}
\end{figure}

We now report the results obtained when applying Crank--Nicolson in time. We report the average number of GMRES iterations together with the average elapsed time in Tables \ref{Instationary_NS_with_CN_table1}--\ref{Instationary_NS_with_CN_table3}, and in Table \ref{Oseen_Instationary_NS_with_CN} the total dimensions of the systems solved and the numbers of Oseen iterations, for different levels of refinements $l$, values of $\beta$, and viscosities $\nu$. As in Section \ref{6_2}, for level of refinement $l$ we divide the time interval into subintervals of length $2^{1-l}$ and consider a spatial uniform grid of refinement level $l$. In Figure \ref{fig:state_adjoint_pressure} we show the numerical solutions of the state and adjoint velocities $\vec{v}$ and $ \vec{\zeta}$, at time $t=1$, and of the pressure $p$, at time $t=1.0625$, for $\nu=\frac{1}{100}$, $\beta=10^{-1}$, and $l=4$.

\begin{table}[!ht]
\caption{Average GMRES iterations and CPU times for instationary Navier--Stokes control problem, with Crank--Nicolson in time ($\tau = h$), for $\nu = \frac{1}{20}$ and a range of $l, \beta$. In brackets are the results for the corresponding Stokes control problem.}\label{Instationary_NS_with_CN_table1}
\begin{small}
\begin{center}
\begin{tabular}{|c||c|c|c|c|c|c|c|c|c|c|c|c|c|c|}
\hline
\multicolumn{1}{|c||}{ \!\!\!\!\!\!\!\!\!\!\!\! \phantom{$1$}} & \multicolumn{2}{c|}{$\beta=10^0$}& \multicolumn{2}{c|}{$\beta=10^{-1}$}& \multicolumn{2}{c|}{$\beta=10^{-2}$}& \multicolumn{2}{c|}{$\beta=10^{-3}$}& \multicolumn{2}{c|}{$\beta=10^{-4}$}& \multicolumn{2}{c|}{$\beta=10^{-5}$}& \multicolumn{2}{c|}{$\beta=10^{-6}$}\\
\cline{2-15}
\!\!\!\!\!\!\! $l$  \!\!\!\!\!\!\!\! & \!\!\!\!\!\! $\texttt{it}$ \!\!\!\!\!\! & \!\!\!\!\!\! CPU \!\!\!\!\!\! & \!\!\!\!\!\! $\texttt{it}$ \!\!\!\!\!\! & \!\!\!\!\!\! CPU \!\!\!\!\!\! & \!\!\!\!\!\! $\texttt{it}$ \!\!\!\!\!\! & \!\!\!\!\!\! CPU \!\!\!\!\!\! & \!\!\!\!\!\! $\texttt{it}$ \!\!\!\!\!\! & \!\!\!\!\!\! CPU \!\!\!\!\!\! & \!\!\!\!\!\! $\texttt{it}$ \!\!\!\!\!\! & \!\!\!\!\!\! CPU \!\!\!\!\!\! & \!\!\!\!\!\! $\texttt{it}$ \!\!\!\!\!\! & \!\!\!\!\!\! CPU \!\!\!\!\!\! & \!\!\!\!\!\! $\texttt{it}$ \!\!\!\!\!\! & \!\!\!\!\!\! CPU \!\!\!\!\!\! \\
\hline
\hline
\multirow{2}{*}{ \!\!\!\!\! $2$ \!\!\!\!\!\!\!\! } & \!\!\!\!\!\! 16 \!\!\!\!\!\! & \!\!\!\!\!\! 0.73 \!\!\!\!\!\! & \!\!\!\!\!\! 15 \!\!\!\!\!\! & \!\!\!\!\!\! 0.68 \!\!\!\!\!\! & \!\!\!\!\!\! 12 \!\!\!\!\!\! & \!\!\!\!\!\! 0.53 \!\!\!\!\!\! & \!\!\!\!\!\! 10 \!\!\!\!\!\! & \!\!\!\!\!\! 0.44 \!\!\!\!\!\! & \!\!\!\!\!\! 9 \!\!\!\!\!\! & \!\!\!\!\!\! 0.39 \!\!\!\!\!\! & \!\!\!\!\!\! 9 \!\!\!\!\!\! & \!\!\!\!\!\! 0.37 \!\!\!\!\!\! & \!\!\!\!\!\! 8 \!\!\!\!\!\! & \!\!\!\!\!\! 0.36 \!\!\!\!\!\! \\
& \!\!\!\!\!\! (14) \!\!\!\!\!\! & \!\!\!\!\!\! (0.54) \!\!\!\!\!\! & \!\!\!\!\!\! (15) \!\!\!\!\!\! & \!\!\!\!\!\! (0.68) \!\!\!\!\!\! & \!\!\!\!\!\! (16) \!\!\!\!\!\! & \!\!\!\!\!\! (0.67) \!\!\!\!\!\! & \!\!\!\!\!\! (15) \!\!\!\!\!\! & \!\!\!\!\!\! (0.62) \!\!\!\!\!\! & \!\!\!\!\!\! (12) \!\!\!\!\!\! & \!\!\!\!\!\! (0.50) \!\!\!\!\!\! & \!\!\!\!\!\! (10) \!\!\!\!\!\! & \!\!\!\!\!\! (0.42) \!\!\!\!\!\! & \!\!\!\!\!\! (9) \!\!\!\!\!\! & \!\!\!\!\!\! (0.37) \!\!\!\!\!\! \\
\hline
\multirow{2}{*}{ \!\!\!\!\! $3$ \!\!\!\!\!\!\!\! } & \!\!\!\!\!\! 18 \!\!\!\!\!\! & \!\!\!\!\!\! 3.40 \!\!\!\!\!\! & \!\!\!\!\!\! 17 \!\!\!\!\!\! & \!\!\!\!\!\! 3.23 \!\!\!\!\!\! & \!\!\!\!\!\! 15 \!\!\!\!\!\! & \!\!\!\!\!\! 2.14 \!\!\!\!\!\! & \!\!\!\!\!\! 12 \!\!\!\!\!\! & \!\!\!\!\!\! 1.68 \!\!\!\!\!\! & \!\!\!\!\!\! 10 \!\!\!\!\!\! & \!\!\!\!\!\! 1.93 \!\!\!\!\!\! & \!\!\!\!\!\! 10 \!\!\!\!\!\! & \!\!\!\!\!\! 1.56 \!\!\!\!\!\! & \!\!\!\!\!\! 9 \!\!\!\!\!\! & \!\!\!\!\!\! 1.55 \!\!\!\!\!\! \\
& \!\!\!\!\!\! (15) \!\!\!\!\!\! & \!\!\!\!\!\! (2.89) \!\!\!\!\!\! & \!\!\!\!\!\! (16) \!\!\!\!\!\! & \!\!\!\!\!\! (3.11) \!\!\!\!\!\! & \!\!\!\!\!\! (17) \!\!\!\!\!\! & \!\!\!\!\!\! (3.28) \!\!\!\!\!\! & \!\!\!\!\!\! (16) \!\!\!\!\!\! & \!\!\!\!\!\! (2.88) \!\!\!\!\!\! & \!\!\!\!\!\! (15) \!\!\!\!\!\! & \!\!\!\!\!\! (2.69) \!\!\!\!\!\! & \!\!\!\!\!\! (13) \!\!\!\!\!\! & \!\!\!\!\!\! (1.23) \!\!\!\!\!\! & \!\!\!\!\!\! (10) \!\!\!\!\!\! & \!\!\!\!\!\! (1.36) \!\!\!\!\!\! \\
\hline
\multirow{2}{*}{ \!\!\!\!\! $4$ \!\!\!\!\!\!\!\! } & \!\!\!\!\!\! 18 \!\!\!\!\!\! & \!\!\!\!\!\! 22.7 \!\!\!\!\!\! & \!\!\!\!\!\! 19 \!\!\!\!\!\! & \!\!\!\!\!\! 22.9 \!\!\!\!\!\! & \!\!\!\!\!\! 18 \!\!\!\!\!\! & \!\!\!\!\!\! 21.2 \!\!\!\!\!\! & \!\!\!\!\!\! 15 \!\!\!\!\!\! & \!\!\!\!\!\! 17.4 \!\!\!\!\!\! & \!\!\!\!\!\! 12 \!\!\!\!\!\! & \!\!\!\!\!\! 11.9 \!\!\!\!\!\! & \!\!\!\!\!\! 11 \!\!\!\!\!\! & \!\!\!\!\!\! 12.6 \!\!\!\!\!\! & \!\!\!\!\!\! 10 \!\!\!\!\!\! & \!\!\!\!\!\! 11.8 \!\!\!\!\!\! \\
& \!\!\!\!\!\! (16) \!\!\!\!\!\! & \!\!\!\!\!\! (16.5) \!\!\!\!\!\! & \!\!\!\!\!\! (18) \!\!\!\!\!\! & \!\!\!\!\!\! (18.3) \!\!\!\!\!\! & \!\!\!\!\!\! (19) \!\!\!\!\!\! & \!\!\!\!\!\! (18.7) \!\!\!\!\!\! & \!\!\!\!\!\! (16) \!\!\!\!\!\! & \!\!\!\!\!\! (18.3) \!\!\!\!\!\! & \!\!\!\!\!\! (16) \!\!\!\!\!\! & \!\!\!\!\!\! (18.4) \!\!\!\!\!\! & \!\!\!\!\!\! (15) \!\!\!\!\!\! & \!\!\!\!\!\! (16.0) \!\!\!\!\!\! & \!\!\!\!\!\! (13) \!\!\!\!\!\! & \!\!\!\!\!\! (12.6) \!\!\!\!\!\! \\
\hline
\multirow{2}{*}{ \!\!\!\!\! $5$ \!\!\!\!\!\!\!\! } & \!\!\!\!\!\! 19 \!\!\!\!\!\! & \!\!\!\!\!\! 170 \!\!\!\!\!\! & \!\!\!\!\!\! 19 \!\!\!\!\!\! & \!\!\!\!\!\! 173 \!\!\!\!\!\! & \!\!\!\!\!\! 18 \!\!\!\!\!\! & \!\!\!\!\!\! 162 \!\!\!\!\!\! & \!\!\!\!\!\! 17 \!\!\!\!\!\! & \!\!\!\!\!\! 151 \!\!\!\!\!\! & \!\!\!\!\!\! 15 \!\!\!\!\!\! & \!\!\!\!\!\! 122 \!\!\!\!\!\! & \!\!\!\!\!\! 13 \!\!\!\!\!\! & \!\!\!\!\!\! 98.4 \!\!\!\!\!\! & \!\!\!\!\!\! 11 \!\!\!\!\!\! & \!\!\!\!\!\! 85.0 \!\!\!\!\!\! \\
& \!\!\!\!\!\! (16) \!\!\!\!\!\! & \!\!\!\!\!\! (139) \!\!\!\!\!\! & \!\!\!\!\!\! (19) \!\!\!\!\!\! & \!\!\!\!\!\! (163) \!\!\!\!\!\! & \!\!\!\!\!\! (20) \!\!\!\!\!\! & \!\!\!\!\!\! (171) \!\!\!\!\!\! & \!\!\!\!\!\! (19) \!\!\!\!\!\! & \!\!\!\!\!\! (158) \!\!\!\!\!\! & \!\!\!\!\!\! (16) \!\!\!\!\!\! & \!\!\!\!\!\! (128) \!\!\!\!\!\! & \!\!\!\!\!\! (15) \!\!\!\!\!\! & \!\!\!\!\!\! (119) \!\!\!\!\!\! & \!\!\!\!\!\! (15) \!\!\!\!\!\! & \!\!\!\!\!\! (103) \!\!\!\!\!\! \\
\hline
\multirow{2}{*}{ \!\!\!\!\! $6$ \!\!\!\!\!\!\!\! } & \!\!\!\!\!\! 21 \!\!\!\!\!\! & \!\!\!\!\!\! 1948 \!\!\!\!\!\! & \!\!\!\!\!\! 21 \!\!\!\!\!\! & \!\!\!\!\!\! 1898 \!\!\!\!\!\! & \!\!\!\!\!\! 21 \!\!\!\!\!\! & \!\!\!\!\!\! 1848 \!\!\!\!\!\! & \!\!\!\!\!\! 18 \!\!\!\!\!\! & \!\!\!\!\!\! 1587 \!\!\!\!\!\! & \!\!\!\!\!\! 17 \!\!\!\!\!\! & \!\!\!\!\!\! 1448 \!\!\!\!\!\! & \!\!\!\!\!\! 15 \!\!\!\!\!\! & \!\!\!\!\!\! 1295 \!\!\!\!\!\! & \!\!\!\!\!\! 13 \!\!\!\!\!\! & \!\!\!\!\!\! 1022 \!\!\!\!\!\! \\
& \!\!\!\!\! (22) \!\!\!\!\! & \!\!\!\!\! (1758) \!\!\!\!\! & \!\!\!\!\! (24) \!\!\!\!\! & \!\!\!\!\! (1915) \!\!\!\!\! & \!\!\!\!\! (26) \!\!\!\!\! & \!\!\!\!\! (2087) \!\!\!\!\! & \!\!\!\!\! (18) \!\!\!\!\! & \!\!\!\!\! (1437) \!\!\!\!\! & \!\!\!\!\! (17) \!\!\!\!\! & \!\!\!\!\! (1344) \!\!\!\!\! & \!\!\!\!\! (15) \!\!\!\!\! & \!\!\!\!\! (1149) \!\!\!\!\! & \!\!\!\!\! (15) \!\!\!\!\! & \!\!\!\!\! (1155) \!\!\!\!\! \\
\hline
\end{tabular}
\end{center}
\end{small}
\end{table}

\begin{table}[!ht]
\caption{Average GMRES iterations and CPU times for instationary Navier--Stokes control problem, with Crank--Nicolson in time ($\tau \! = \! h$), for $\nu \! = \! \frac{1}{100}$ and a range of $l, \, \beta$.}\label{Instationary_NS_with_CN_table2}
\begin{small}
\begin{center}
\begin{tabular}{|c||c|c|c|c|c|c|c|c|c|c|c|c|c|c|}
\hline
\multicolumn{1}{|c||}{\phantom{$1$}} & \multicolumn{2}{c|}{$\beta=10^0$}& \multicolumn{2}{c|}{$\beta=10^{-1}$}& \multicolumn{2}{c|}{$\beta=10^{-2}$}& \multicolumn{2}{c|}{$\beta=10^{-3}$}& \multicolumn{2}{c|}{$\beta=10^{-4}$}& \multicolumn{2}{c|}{$\beta=10^{-5}$}& \multicolumn{2}{c|}{$\beta=10^{-6}$}\\
\cline{2-15}
$l$  & \!\!\!\! $\texttt{it}$ \!\!\!\! & \!\!\!\! CPU \!\!\!\! & \!\!\!\! $\texttt{it}$ \!\!\!\! & \!\!\!\! CPU \!\!\!\! & \!\!\!\! $\texttt{it}$ \!\!\!\! & \!\!\!\! CPU \!\!\!\! & \!\!\!\! $\texttt{it}$ \!\!\!\! & \!\!\!\! CPU \!\!\!\! & \!\!\!\! $\texttt{it}$ \!\!\!\! & \!\!\!\! CPU \!\!\!\! & \!\!\!\! $\texttt{it}$ \!\!\!\! & \!\!\!\! CPU \!\!\!\! & \!\!\!\! $\texttt{it}$ \!\!\!\! & \!\!\!\! CPU \!\!\!\! \\
\hline
\hline
$2$ & \!\!\!\! 16 \!\!\!\! & \!\!\!\! 0.74 \!\!\!\! & \!\!\!\! 13 \!\!\!\! & \!\!\!\! 0.64 \!\!\!\! & \!\!\!\! 11 \!\!\!\! & \!\!\!\! 0.51 \!\!\!\! & \!\!\!\! 10 \!\!\!\! & \!\!\!\! 0.45 \!\!\!\! & \!\!\!\! 9 \!\!\!\! & \!\!\!\! 0.40 \!\!\!\! & \!\!\!\! 9 \!\!\!\! & \!\!\!\! 0.38 \!\!\!\! & \!\!\!\! 8 \!\!\!\! & \!\!\!\! 0.38 \!\!\!\! \\
\hline
$3$ & \!\!\!\! 21 \!\!\!\! & \!\!\!\! 3.80 \!\!\!\! & \!\!\!\! 19 \!\!\!\! & \!\!\!\! 3.72 \!\!\!\! & \!\!\!\! 13 \!\!\!\! & \!\!\!\! 2.80 \!\!\!\! & \!\!\!\! 10 \!\!\!\! & \!\!\!\! 2.20 \!\!\!\! & \!\!\!\! 10 \!\!\!\! & \!\!\!\! 2.20 \!\!\!\! & \!\!\!\! 9 \!\!\!\! & \!\!\!\! 1.77 \!\!\!\! & \!\!\!\! 9 \!\!\!\! & \!\!\!\! 1.87 \!\!\!\! \\
\hline
$4$ & \!\!\!\! 23 \!\!\!\! & \!\!\!\! 22.6 \!\!\!\! & \!\!\!\! 22 \!\!\!\! & \!\!\!\! 22.2 \!\!\!\! & \!\!\!\! 18 \!\!\!\! & \!\!\!\! 18.5 \!\!\!\! & \!\!\!\! 12 \!\!\!\! & \!\!\!\! 15.4 \!\!\!\! & \!\!\!\! 11 \!\!\!\! & \!\!\!\! 13.6 \!\!\!\! & \!\!\!\! 10 \!\!\!\! & \!\!\!\! 13.3 \!\!\!\! & \!\!\!\! 10 \!\!\!\! & \!\!\!\! 12.2 \!\!\!\! \\
\hline
$5$ & \!\!\!\! 22 \!\!\!\! & \!\!\!\! 187 \!\!\!\! & \!\!\!\! 21 \!\!\!\! & \!\!\!\! 184 \!\!\!\! & \!\!\!\! 19 \!\!\!\! & \!\!\!\! 166 \!\!\!\! & \!\!\!\! 16 \!\!\!\! & \!\!\!\! 135 \!\!\!\! & \!\!\!\! 12 \!\!\!\! & \!\!\!\! 103 \!\!\!\! & \!\!\!\! 11 \!\!\!\! & \!\!\!\! 87.7 \!\!\!\! & \!\!\!\! 11 \!\!\!\! & \!\!\!\! 89.7 \!\!\!\! \\
\hline
$6$ & \!\!\!\! 24 \!\!\!\! & \!\!\!\! 2141 \!\!\!\! & \!\!\!\! 24 \!\!\!\! & \!\!\!\! 2087 \!\!\!\! & \!\!\!\! 22 \!\!\!\! & \!\!\!\! 1922 \!\!\!\! & \!\!\!\! 18 \!\!\!\! & \!\!\!\! 1507 \!\!\!\! & \!\!\!\! 15 \!\!\!\! & \!\!\!\! 1272 \!\!\!\! & \!\!\!\! 12 \!\!\!\! & \!\!\!\! 973 \!\!\!\! & \!\!\!\! 11 \!\!\!\! & \!\!\!\! 913 \!\!\!\! \\
\hline
\end{tabular}
\end{center}
\end{small}
\end{table}

\begin{table}[!ht]
\caption{Average GMRES iterations and CPU times for instationary Navier--Stokes control problem, with Crank--Nicolson in time ($\tau \! = \! h$), for $\nu \! = \! \frac{1}{500}$ and a range of $l, \, \beta$.}\label{Instationary_NS_with_CN_table3}
\begin{small}
\begin{center}
\begin{tabular}{|c||c|c|c|c|c|c|c|c|c|c|c|c|c|c|}
\hline
\multicolumn{1}{|c||}{\phantom{$1$}} & \multicolumn{2}{c|}{$\beta=10^0$}& \multicolumn{2}{c|}{$\beta=10^{-1}$}& \multicolumn{2}{c|}{$\beta=10^{-2}$}& \multicolumn{2}{c|}{$\beta=10^{-3}$}& \multicolumn{2}{c|}{$\beta=10^{-4}$}& \multicolumn{2}{c|}{$\beta=10^{-5}$}& \multicolumn{2}{c|}{$\beta=10^{-6}$}\\
\cline{2-15}
$l$ & \!\!\!\! $\texttt{it}$ \!\!\!\! & \!\!\!\! CPU \!\!\!\! & \!\!\!\! $\texttt{it}$ \!\!\!\! & \!\!\!\! CPU \!\!\!\! & \!\!\!\! $\texttt{it}$ \!\!\!\! & \!\!\!\! CPU \!\!\!\! & \!\!\!\! $\texttt{it}$ \!\!\!\! & \!\!\!\! CPU \!\!\!\! & \!\!\!\! $\texttt{it}$ \!\!\!\! & \!\!\!\! CPU \!\!\!\! & \!\!\!\! $\texttt{it}$ \!\!\!\! & \!\!\!\! CPU \!\!\!\! & \!\!\!\! $\texttt{it}$ \!\!\!\! & \!\!\!\! CPU \!\!\!\! \\
\hline
\hline
$2$ & \!\!\!\! 16$\dagger$ \!\!\!\! & \!\!\!\! 0.75$\dagger$ \!\!\!\! & \!\!\!\! 14$\dagger$ \!\!\!\! & \!\!\!\! 0.64$\dagger$ \!\!\!\! & \!\!\!\! 10$\dagger$ \!\!\!\! & \!\!\!\! 0.46$\dagger$ \!\!\!\! & \!\!\!\! 10 \!\!\!\! & \!\!\!\! 0.43 \!\!\!\! & \!\!\!\! 9 \!\!\!\! & \!\!\!\! 0.40 \!\!\!\! & \!\!\!\! 8 \!\!\!\! & \!\!\!\! 0.38 \!\!\!\! & \!\!\!\! 8 \!\!\!\! & \!\!\!\! 0.37 \!\!\!\! \\
\hline
$3$ & \!\!\!\! 26$\dagger$ \!\!\!\! & \!\!\!\! 6.49$\dagger$ \!\!\!\! & \!\!\!\! 20$\dagger$ \!\!\!\! & \!\!\!\! 5.70$\dagger$ \!\!\!\! & \!\!\!\! 13 \!\!\!\! & \!\!\!\! 3.78 \!\!\!\! & \!\!\!\! 10 \!\!\!\! & \!\!\!\! 2.90 \!\!\!\! & \!\!\!\! 9 \!\!\!\! & \!\!\!\! 2.60 \!\!\!\! & \!\!\!\! 9 \!\!\!\! & \!\!\!\! 2.19 \!\!\!\! & \!\!\!\! 9 \!\!\!\! & \!\!\!\! 2.19 \!\!\!\! \\
\hline
$4$ & \!\!\!\! 47 \!\!\!\! & \!\!\!\! 64.6 \!\!\!\! & \!\!\!\! 33 \!\!\!\! & \!\!\!\! 53.5 \!\!\!\! & \!\!\!\! 17 \!\!\!\! & \!\!\!\! 28.2 \!\!\!\! & \!\!\!\! 11 \!\!\!\! & \!\!\!\! 20.3 \!\!\!\! & \!\!\!\! 10 \!\!\!\! & \!\!\!\! 17.7 \!\!\!\! & \!\!\!\! 10 \!\!\!\! & \!\!\!\! 15.8 \!\!\!\! & \!\!\!\! 9 \!\!\!\! & \!\!\!\! 14.6 \!\!\!\! \\
\hline
$5$ & \!\!\!\! 54 \!\!\!\! & \!\!\!\! 476 \!\!\!\! & \!\!\!\! 44 \!\!\!\! & \!\!\!\! 417 \!\!\!\! & \!\!\!\! 28 \!\!\!\! & \!\!\!\! 291 \!\!\!\! & \!\!\!\! 15 \!\!\!\! & \!\!\!\! 148 \!\!\!\! & \!\!\!\! 11 \!\!\!\! & \!\!\!\! 119 \!\!\!\! &  \!\!\!\! 11 \!\!\!\! & \!\!\!\! 109 \!\!\!\! & \!\!\!\! 10 \!\!\!\! & \!\!\!\! 98.2 \!\!\!\! \\
\hline
$6$ & \!\!\!\! 44 \!\!\!\! & \!\!\!\! 3728 \!\!\!\! & \!\!\!\! 37 \!\!\!\! & \!\!\!\! 3140 \!\!\!\! & \!\!\!\! 29 \!\!\!\! & \!\!\!\! 2472 \!\!\!\! & \!\!\!\! 20 \!\!\!\! & \!\!\!\! 1701 \!\!\!\! & \!\!\!\! 14 \!\!\!\! & \!\!\!\! 1157 \!\!\!\! & \!\!\!\! 11 \!\!\!\! & \!\!\!\! 983 \!\!\!\! & \!\!\!\! 11 \!\!\!\! & \!\!\!\! 968 \!\!\!\! \\
\hline
\end{tabular}
\end{center}
\end{small}
\end{table}

\begin{table}[!ht]
\caption{Degrees of freedom (DoF) and number of Oseen iterations required for instationary Navier--Stokes control problem, with Crank--Nicolson in time ($\tau = h$). In each cell are the Oseen iterations for the given $l, \, \nu$, and $\beta=10^{-j}, \, j=0,1,...,6$.}\label{Oseen_Instationary_NS_with_CN}
\begin{small}
\begin{center}
\begin{tabular}{|c|c||ccccccc|ccccccc|ccccccc|}
\hline
$l$ & DoF & \multicolumn{7}{c|}{$\nu=\frac{1}{20}$} & \multicolumn{7}{c|}{$\nu=\frac{1}{100}$} & \multicolumn{7}{c|}{$\nu=\frac{1}{500}$} \\
\hline
\hline
$2$ & 984 & 6 & \!\!\!\!\! 6 & \!\!\!\!\! 6 & \!\!\!\!\! 6 & \!\!\!\!\! 5 & \!\!\!\!\! 4 & \!\!\!\!\! 3 & 12 & \!\!\!\!\! 7 & \!\!\!\!\! 7 & \!\!\!\!\! 7 & \!\!\!\!\! 5 & \!\!\!\!\! 4 & \!\!\!\!\! 4 & $\dagger$ & \!\!\!\!\! $\dagger$ & \!\!\!\!\! $\dagger$ & \!\!\!\!\! 10 & \!\!\!\!\! 7 & \!\!\!\!\! 5 & \!\!\!\!\! 4 \\
\hline
$3$ & 8496 & 5 & \!\!\!\!\! 5 & \!\!\!\!\! 5 & \!\!\!\!\! 6 & \!\!\!\!\! 5 & \!\!\!\!\! 4 & \!\!\!\!\! 4 & 8 & \!\!\!\!\! 8 & \!\!\!\!\! 6 & \!\!\!\!\! 7 & \!\!\!\!\! 6 & \!\!\!\!\! 4 & \!\!\!\!\! 4 & $\dagger$ & \!\!\!\!\! $\dagger$ & \!\!\!\!\! 7 & \!\!\!\!\! 7 & \!\!\!\!\! 7 & \!\!\!\!\! 5 & \!\!\!\!\! 4 \\
\hline
$4$ & 70,752 & 4 & \!\!\!\!\! 4 & \!\!\!\!\! 4 & \!\!\!\!\! 4 & \!\!\!\!\! 4 & \!\!\!\!\! 4 & \!\!\!\!\! 3 & 6 & \!\!\!\!\! 6 & \!\!\!\!\! 5 & \!\!\!\!\! 5 & \!\!\!\!\! 4 & \!\!\!\!\! 4 & \!\!\!\!\! 3 & 18 & \!\!\!\!\! 14 & \!\!\!\!\! 6 & \!\!\!\!\! 5 & \!\!\!\!\! 6 & \!\!\!\!\! 5 & \!\!\!\!\! 4 \\
\hline
$5$ & 577,728 & 4 & \!\!\!\!\! 4 & \!\!\!\!\! 4 & \!\!\!\!\! 3 & \!\!\!\!\! 3 & \!\!\!\!\! 3 & \!\!\!\!\! 3 & 5 & \!\!\!\!\! 4 & \!\!\!\!\! 4 & \!\!\!\!\! 4 & \!\!\!\!\! 4 & \!\!\!\!\! 4 & \!\!\!\!\! 3 & 9 & \!\!\!\!\! 8 & \!\!\!\!\! 6 & \!\!\!\!\! 5 & \!\!\!\!\! 4 & \!\!\!\!\! 4 & \!\!\!\!\! 3 \\
\hline
$6$ & 4,669,824 & 3 & \!\!\!\!\! 3 & \!\!\!\!\! 3 & \!\!\!\!\! 3 & \!\!\!\!\! 3 & \!\!\!\!\! 3 & \!\!\!\!\! 3 & 4 & \!\!\!\!\! 4 & \!\!\!\!\! 4 & \!\!\!\!\! 3 & \!\!\!\!\! 3 & \!\!\!\!\! 3 & \!\!\!\!\! 3 & 5 & \!\!\!\!\! 5 & \!\!\!\!\! 5 & \!\!\!\!\! 4 & \!\!\!\!\! 4 & \!\!\!\!\! 3 & \!\!\!\!\! 3 \\
\hline
\end{tabular}
\end{center}
\end{small}
\end{table}

From Tables \ref{Instationary_NS_with_CN_table1}--\ref{Instationary_NS_with_CN_table3} we observe that the number of iterations required for reaching a prescribed accuracy is, again, roughly constant, increasing only for small $\nu$ and large $\beta$. As experienced above, the CPU time scales approximately linearly with the size of the system, except for very fine grids. Regarding the non-linear iteration, as above we note in Table \ref{Oseen_Instationary_NS_with_CN} that the number of Oseen iterations is decreasing as the grid is refined, while it is increasing for small values of $\nu$ and large values of $\beta$.

\section{Concluding Remarks}
\label{7}
We presented mesh- and parameter-robust preconditioners for distributed (Stokes and) Navier--Stokes control problems, of both stationary and instationary type, coupled with an Oseen linearization. The preconditioners were applied within the flexible GMRES algorithm, and in the instationary setting to backward Euler and Crank--Nicolson discretizations in time. Numerical results demonstrated the versatility and effectiveness of this approach when solving a range of huge-scale linear systems. Future work involves adapting this solver to problems involving more complicated PDEs from fluid dynamics, boundary control problems, and problems with additional algebraic constraints on the state and/or control variables.

\section*{Acknowledgements}
SL acknowledges a University of Edinburgh School of Mathematics PhD studentship. JWP acknowledges the EPSRC (UK) grant EP/S027785/1.

\vspace{-1.2em}

%
%
%

\end{document}